\documentclass[12pt]{article}
\usepackage{amsfonts,amscd,a4}
\usepackage{color}
\usepackage{graphicx}
\usepackage{epsfig}
\usepackage{epic,eepic,pictex}
\usepackage{theorem}
\newtheorem{theo}{Theorem}[section]
\newtheorem{prop}[theo]{Proposition}
\newtheorem{lemma}[theo]{Lemma}
\newtheorem{coro}[theo]{Corollary}
\newtheorem{defi}[theo]{Definition}
{\theorembodyfont{\rm}
\newtheorem{example}[theo]{Example}

}
\newcommand{\bA}{{\bf A}}
\newcommand{\bB}{{\bf B}}
\newcommand{\bC}{{\bf C}}

\newcommand{\bG}{{\bf G}}

\newcommand{\bI}{{\bf I}}
\newcommand{\bJ}{{\bf J}}
\newcommand{\bK}{{\bf K}}
\newcommand{\bL}{{\bf L}}

\newcommand{\cA}{{\mathcal A}}
\newcommand{\cB}{{\mathcal B}}
\newcommand{\cC}{{\mathcal C}}

\newcommand{\cF}{{\mathcal F}}
\newcommand{\cG}{{\mathcal G}}

\newcommand{\cI}{{\mathcal I}}
\newcommand{\cJ}{{\mathcal J}}
\newcommand{\cK}{{\mathcal K}}
\newcommand{\cL}{{\mathcal L}}

\newcommand{\cO}{{\mathcal O}}
\newcommand{\cP}{{\mathcal P}}

\newcommand{\cS}{{\mathcal S}}

\newcommand{\eB}{{\sf B}}

\newcommand{\eD}{{\sf D}}

\newcommand{\eO}{{\sf O}}

\newcommand{\eT}{{\sf T}}

\newcommand{\sC}{{\mathbb C}}

\newcommand{\sM}{{\mathbb M}}
\newcommand{\sN}{{\mathbb N}}

\newcommand{\sR}{{\mathbb R}}

\newcommand{\sT}{{\mathbb T}}

\newcommand{\sZ}{{\mathbb Z}}
\newcommand{\za}{{\vec{\alpha}}}
\newcommand{\zb}{{\vec{\beta}}}

\newcommand{\zs}{{\vec{\sigma}}}
\newcommand{\qed}{\rule{1ex}{1ex}}

\newcommand{\diag}{\mbox{\rm diag} \,}
\newcommand{\dist}{\mbox{\rm dist} \,}
\newcommand{\ess}{\mbox{\rm ess} \,}

\newcommand{\im}{{\rm im} \,}
\newcommand{\ind}{\mbox{\rm ind} \,}
\newcommand{\op}{\mbox{\rm Op} \,}

\newcommand{\rank}{\mbox{\rm rank} \,}

\makeindex
\begin{document}
\title{Fractal algebras of discretization sequences}
\author{Steffen Roch\footnote{Address: Steffen Roch, Technische Universit\"at Darmstadt, Fachbereich Mathematik, Schlossgartenstra{\ss}e 7, 64289
Darmstadt, Germany.}}
\date{Accompanying material to lectures at the\\ Summer School \\ on \\ Applied Analysis \\[1.5mm] Chemnitz, September 2011}
\maketitle
%
%\begin{abstract}
%\end{abstract}
%
%
\newpage
\tableofcontents
\newpage
\section{Introduction} \label{s0}
First a warning: Fractality, in the sense of these lectures, has nothing to do with fractal geometries or broken dimensions or other involved things. Rather, the notion {\em fractal algebra} had been chosen in order to emphasize
an important property of many discretization sequences, namely their {\em self-similarity}, in the sense that each subsequence has the same properties as the full sequence. (But note that self-similarity is also a characteristic
aspect of many fractal sets. I guess that everyone is fascinated by zooming into the Mandelbrot set, which reveals the same details at finer and finer levels.)

We start with a precise definition of the concept of fractality and show that the fractal property is enormously useful for several spectral approximation problems. These results will be illustrated by sequences in the algebra of the
finite sections method for Toeplitz operators. ({\em What else?} one might ask: these algebras (first) played the prominent role in the development of the use of algebraic techniques in numerical analysis, and they were (second) a main object of study in Silbermann's school; so one can hardly think of a lecture on this topic in Chemnitz, which does not come across with these algebras.)

Then we discuss some structural consequences of fractality, which are related with the notion of a compact sequence. Discretized Cuntz algebras will show that idea of fractality is also a very helpful guide in order to analyze concrete  algebras of approximation sequences, which illustrates the importance of the idea of {\em fractal restriction}. Our final example is the algebra of the finite sections method for band operators. This algebra is not fractal, but
has a related property which we call {\em essential fractality} and which is related with the approximation of points in the essential spectrum.

I suppose that the participants have some (really) basic knowledge on $C^*$-algebras and their representations. A short script will be available during the Summer School. The textbooks and review papers \cite{Boe3,BSi2,HRS1,HRS2,Roc10,Sil9} provide both an introduction to the field and suggestions for further reading.
\section{Stability} \label{s1}
\subsection{Algebras of matrix sequences} \label{ss11}
Let $(A_n)$ be a sequence of squared matrices of increasing size. We think of $A_n$ as the $n$th approximant of a bounded linear operator $A$ on a Hilbert space $H$. A basic question in numerical analysis asks if the method $(A_n)$ is applicable to $A$ in the sense that the equations $A_n x_n = f_n$ (with $f_n$ a suitable approximant of an element $f \in H$) are uniquely solvable for all sufficiently large $n$ and all right hand sides and if their solutions converge to a solution of the equation $Ax=f$. Typically, one can answer this question in the affirmative if one is able to decide the stability question for the sequence $(A_n)$. The sequence $(A_n)$ is called stable if there is an $n_0$ such that the matrices $A_n$ are invertible for $n \ge n_0$ and if the norms of their inverses are uniformly bounded. It turns out that the stability of a sequence is equivalent to the invertibility of a certain element (related with $(A_n)$ in a suitably constructed $C^*$-algebra. This is the point where the story begins.

Given a sequence $\delta : \sN \to \sN$ tending to infinity, let $\cF^\delta$ denote the set of all bounded sequences $\bA = (A_n)$ of matrices $A_n \in \sC^{\delta(n) \times \delta(n)}$. Equipped with the operations
\[
(A_n) + (B_n) := (A_n + B_n), \quad (A_n) (B_n) := (A_n B_n), \quad (A_n)^* := (A_n^*)
\]
and the norm
\[
\|\bA\|_\cF := \sup \|A_n\|,
\]
the set $\cF^\delta$ becomes a $C^*$-algebra, and the set $\cG^\delta$ of all sequences $(A_n) \in \cF^\delta$ with $\lim \|A_n\| = 0$ forms a closed ideal of $\cF^\delta$. We call $\cF^\delta$ the algebra of matrix sequences with dimension function $\delta$ and $\cG^\delta$ the associated ideal of zero sequences. When the concrete choice of $\delta$ is irrelevant or evident from the context, we will simply write $\cF$ and $\cG$ in place of $\cF^\delta$ and $\cG^\delta$.

The relevance of the algebra $\cF$ and its ideal $\cG$ in our context stems from the fact (following via a simple Neumann series argument which is left as an exercise) that a sequence $(A_n) \in \cF$ is stable if, and only if, the coset $(A_n) + \cG$ is invertible in the quotient algebra $\cF/\cG$. This equivalence is also known as Kozak's theorem. Thus, every stability problem is equivalent to an invertibility problem in a suitably chosen $C^*$-algebra, and to understand stability means to understand subalgebras of the quotient algebra $\cF/\cG$. Note in this connection that
\begin{equation} \label{e170309.3}
\limsup \|A_n\| = \|(A_n) + \cG \|_{\cF/\cG}
\end{equation}
for each sequence $(A_n)$ in $\cF$ (a simple exercise again).

It will sometimes be desirable to identify the entries of a sequence $(A_n)$ with operators acting on a common Hilbert  space. The general setting is as follows. Let $H$ be a separable infinite-dimensional Hilbert space and $\cP =
(P_n)$ a sequence of orthogonal projections of finite rank on $H$ which converges strongly to the identity operator $I$ on $H$, i.e., $\|P_n x - x\| \to 0$ for every $x \in H$. A sequence $\cP$ with these properties is also called a filtration on $H$. A typical filtration is that of the finite sections method, where one fixes an orthonormal basis $\{e_i\}_{i \in \sN}$ of $H$ and defines $P_n$ as the orthogonal projection from $H$ onto the linear span of $e_1, \, \ldots, \, e_n$.

Given a filtration $\cP = (P_n)$, we let $\cF^\cP$ stand for the set of all sequences $\bA = (A_n)$ of operators $A_n : \im P_n \to \im P_n$ with the property that the sequences $(A_nP_n)$ and $(A_n^*P_n)$ converge strongly. The set of all sequences $(A_n) \in \cF^\cP$ with $\|A_n P_n\| \to 0$ is denoted by $\cG^\cP$. By the uniform boundedness principle, the quantity $\sup \|A_n P_n\|$ is finite for every sequence $\bA$ in $\cF^\cP$. Thus, if we fix a basis in $\im P_n$ and identify each operator $A_n$ on $\im P_n$ with its matrix representation with respect to this basis, then we can think of $\cF^\cP$ as a $C^*$-subalgebra of the algebra of matrix functions with dimension function $\delta(n) = \rank P_n$. Then $\cG^\cP$ can be identified with $\cG^\delta$. Note that the mapping
\begin{equation} \label{e91.5}
W : \cF^\cP \to L(H), \quad (A_n) \mapsto \mbox{s-lim} \, A_n P_n
\end{equation}
is a $^*$-homomorphism, which we call the {\em consistency map}. Again we simply write $\cF$ and $\cG$ in place of $\cF^\cP$ and $\cG^\cP$ if the concrete choice of the filtration is irrelevant or evident from the context.

The attentive reader will note that many of the concepts and facts presented in this paper make sense and hold in the more general context when $\cF$ is a direct product of a countable family of $C^*$-algebras. There are also Banach algebraic version of parts of his story.
\subsection{Discretization of the Toeplitz algebra} \label{ss13}
Now we introduce a concrete $C^*$-algebra of matrix sequences which will serve as a running example throughout these lectures. This algebra is generated by finite sections sequences of Toeplitz operators or, slightly more general, of operators in the Toeplitz algebra. Below we present only some basic facts about Toeplitz operators and their finite sections, as far as we will need them in this text. Much more on this fascinating topic can be found, for example, in \cite{BGr5,BSi1,BSi2}.

There are several characterizations of the Toeplitz algebra. From the view point of abstract $C^*$-algebra theory, it can be defined as the universal algebra $C^*(s)$ generated by one isometry, i.e., by an element $s$ such that $s^*s$ is the identity element. The universal property of $C^*(s)$ implies that whenever $S$ is an isometry in a $C^*$-algebra $\cA$, then there is a $^*$-homomorphism from $C^*(s)$ onto the smallest $C^*$-subalgebra of $\cA$ containing $S$ which sends $s$ to $S$. Coburn \cite{Cob1} showed that the algebra $C^*(s)$ is $^*$-isomorphic to the smallest closed $^*$-subalgebra $\eT (C)$ of $L(l^2(\sZ^+))$ which contains the (isometric, but not unitary) operator
\[
V : l^2(\sZ^+) \to l^2(\sZ^+), \quad (x_k)_{k \ge 0} \mapsto (0, \,  x_0, \, x_1, \, \ldots)
\]
of forward shift. The algebra $\eT(C)$ is known as the {\em Toeplitz algebra}, since each of its elements is of the form $T(c) + K$ where $T(c)$ is a Toeplitz operator and $K$ a compact operator. To recall the definition of a Toeplitz operator, let $a$ be a function in $L^\infty (\sT)$ with $k$th Fourier coefficient
\[
a_k := \frac{1}{2 \pi} \int_0^{2 \pi} a(e^{i \theta}) e^{-ik
\theta} \, d \theta, \quad k \in \sZ.
\]
Then the {\em Laurent operator} \index{operator!Laurent} $L(a)$ on $l^2(\sZ)$, the {\em Toeplitz operator} \index{operator!Toeplitz} $T(a)$ on $l^2(\sZ^+)$, and the {\em Hankel operator} \index{operator!Hankel} $H(a)$ on $l^2(\sZ^+)$ {\em with generating function} $a$ are defined via their matrix representations with respect to the standard bases of $l^2(\sZ)$ and $l^2(\sZ^+)$ by
\[
L(a) = (a_{i-j})_{i, j = - \infty}^\infty, \quad
T(a) = (a_{i-j})_{i, j = 0}^\infty, \quad \mbox{and} \quad
H(a) = (a_{i+j+1})_{i, j = 0}^\infty.
\]
These operators are bounded, and
\[
\|H(a)\| \le \|T(a)\| = \|L(a)\| = \|a\|_\infty.
\]
It is also useful to know that $L(a)$ and $T(a)$ are compact if and only if $a$ is the zero function, whereas $H(a)$ is compact for each continuous function $a$. With these notations, one has
\begin{theo}\label{t14.10}
$\quad \eT(C) = \{ T(a) + K : a \in C(\sT) \; \mbox{and} \; K \in K(l^2(\sZ^+)) \}$.
\end{theo}
To discretize the Toeplitz algebra $\eT(C)$, consider the orthogonal projections
\[
P_n : l^2(\sZ^+) \to l^2(\sZ^+), \quad (x_0, \, x_1, \, x_2, \,
\dots, \,) \mapsto (x_0, \, x_1, \, \dots, \, x_{n-1}, \, 0, \, 0, \, \dots)
\]
which converge strongly to the identity operator. Hence, $\cP := (P_n)_{n \ge 1}$ is a filtration. We let $\cS(\eT(C))$ stand for the $C^*$-algebra of the finite sections discretization of the Toeplitz algebra, i.e., for the smallest closed subalgebra of $\cF^\cP$ which contains all sequences $(P_n A P_n)$ with $A \in \eT(C)$. By the way, one can show that already the sequences $(P_n T(a) P_n)$ with $a \in C(\sT)$ generate $\cS(\eT(C))$.

It is a lucky circumstance that, similarly to the Toeplitz algebra $\eT(C)$, all elements of $\cS(\eT(C))$ are known explicitly. This makes the algebra $\cS(\eT(C))$ to an ideal model in numerical analysis, and this is the reason why this algebra will serve as an illustrative example in this text. For the description of $\cS(\eT(C))$, we will need the {\em reflection operators} \index{operator!reflection}
\[
R_n : l^2(\sZ^+) \to l^2(\sZ^+), \quad (x_0, \, x_1, \, \dots ) \mapsto (x_{n-1}, \, x_{n-2}, \, \dots , \, x_1, \, x_0, \, 0, \, 0, \, \dots).
\]
\begin{theo}[B\"ottcher/Silbermann] \label{t14.21}
The algebra $\cS(\eT(C))$ coincides with the set of all sequences
\begin{equation} \label{e290409.1}
(P_nT(a)P_n + P_nKP_n + R_nLR_n + G_n)
\end{equation}
where $a \in C(\sT)$, $K$ and $L$ are compact on $l^2(\sZ^+)$, and $(G_n) \in \cG^\cP$.
\end{theo}
{\bf Proof.} Denote the set of all sequences of the form (\ref{e290409.1}) by $\cS_1$ for a moment. In a first step we show that $\cS_1$ is a symmetric algebra. This follows essentially from Widom's identity
\begin{equation} \label{e170309.4}
P_n T(ab) P_n = P_n T(a) P_n T(b) P_n + P_n H(a) H(\tilde{b}) P_n + R_n H(\tilde{a}) H(b) R_n,
\end{equation}
where $\tilde{a}(t) := a(t^{-1})$, and from the compactness of Hankel operators with continuous generating function.

The proof that $\cS_1$ is closed (hence, a $C^*$-algebra) proceeds in the standard way if one employs the fact that the strong limits $W(\bA) := \mbox{s-lim} \, A_n P_n$ and $\widetilde{W}(\bA) := \mbox{s-lim} \, R_nAR_n$ exist for each sequence $\bA := (A_n) \in \cS_1$ and that
\begin{equation} \label{e14.23}
W((P_n T(a) P_n + P_n K P_n + R_n L R_n + G_n)) = T(a) + K
\end{equation}
and
\begin{equation} \label{e14.24}
\widetilde{W}((P_n T(a) P_n + P_n K P_n + R_n L R_n + G_n)) =
T(\tilde{a}) + L.
\end{equation}
Since the generating sequences of $\cS(\eT(C))$ belong to $\cS_1$ and $\cS_1$ is a closed algebra, we conclude that $\cS(\eT(C)) \subseteq \cS_1$.

For the reverse inclusion we have to show that the sequence $(R_n L R_n)$ belongs to $\cS(\eT(C))$ for every compact operator $L$ and that $\cG \subseteq \cS(\eT(C))$. Note that $V^*$ is the operator of backward shift and that all non-negative powers of $V$ and $V^*$ are Toeplitz operators with polynomial generating function. Hence, the identities
\[
(R_n V^i P_1 (V^*)^j R_n) = (P_n (V^*)^i P_n) (R_n P_1 R_n) (P_n V^j P_n)
\]
and $(R_n P_1 R_n) = (P_n) - (P_n V P_n) (P_n V^* P_n)$ imply that $\cS(\eT(C))$ contains all sequences $(R_n L R_n)$ with $L$ a finite linear combination of operators of the form $V_i P_1 V_{-j}$ with $i, \, j \ge 0$. Since these operators form a dense subset of $K(l^2(\sZ^+))$, the first claim follows. The inclusion $\cG^\cP \subseteq \cS(\eT(C))$ is a consequence of a more general result which we formulate as a separate proposition. \hfill \qed
\\[3mm]
The following proposition shows a close symbiosis between sequences of the form $(P_n K P_n)$ with $K$ compact and sequences which tend to zero in the norm: each algebra which contains all sequences $(P_n K P_n)$ also contains all sequences tending to zero. The only (evidently necessary) obstruction is that no two of the $P_n$ coincide.
\begin{prop} \label{p12.20}
Let $\cP = (P_n)$ be a filtration on a Hilbert space $H$ and suppose that $P_m \neq P_n$ whenever $m \neq n$. Then the ideal $\cG^\cP$ of the zero sequences is contained in the smallest closed subalgebra $\cJ$ of $\cF^\cP$ which contains all sequences $(P_n K P_n)$ with $K$ compact.
\end{prop}
{\bf Proof.} It is sufficient to show that, for each $n_0 \in \sN$, there is a sequence $(G_n)$ in $\cJ$ such that $G_{n_0}$ is a projection of rank 1 and $G_n = 0$ for all $n \neq n_0$. Since the matrix algebras $\sC^{k \times k}$ are simple, this fact already implies that each sequence $(G_n)$ with arbitrarily prescribed $G_{n_0} \in L(\im P_{n_0})$ and $G_n = 0$ for $n \neq n_0$ belongs to $\cJ$. Since $\cG^\cP$ is generated (as a Banach space) by sequences of this special form, the assertion follows.

Let $n_0 \in \sN$, put
\[
\sN_< := \{ n \in \sN : \im P_n \cap \im P_{n_0} \; \mbox{is a proper subspace of} \; \im P_{n_0} \},
\]
and set $\sN_> := \sN \setminus (\{n_0\} \cup \sN_< )$. The set $\sN_<$ is at most countable. If $n \in \sN_<$, then none of the closed linear spaces
$\im P_n \cap \im P_{n_0}$ has interior points relative to $\im P_{n_0}$. By the Baire category theorem, $\cup_{n \in \sN_<} (\im P_n \cap \im P_{n_0})$ is a proper subset of $\im P_{n_0}$. Choose a unit vector
\[
f \in \im P_{n_0} \setminus \cup_{n \in \sN_<} (\im P_n \cap \im P_{n_0}).
\]
Then $\|P_n f\| < 1$ for all $n \in \sN_<$ by the Pythagoras theorem. (Indeed, otherwise $\|P_n f\| = 1$, and the equality $1 = \|f\|^2 = \|P_n f\|^2 + \|f - P_n f\|^2$ implies $f = P_n f$, whence $f \in \im P_n$.)

Let $Q_n:= I - P_n$. If $n \in \sN_>$, then $\im P_n \cap \im P_{n_0} = \im P_{n_0}$ by the definition of $\sN_>$. Thus, $\im P_{n_0} \subseteq \im P_n$, and since no two of the projections $P_n$ coincide, this implies that $\im P_{n_0}$ is a proper subspace of $\im P_n$ and $\im Q_n$ is a proper subspace of $\im Q_{n_0}$ for $n \in \sN_>$. Again by the Baire category theorem, $\cup_{n \in \sN_>} \im Q_n$ is a proper subset of $\im Q_{n_0}$. Choose a unit vector
\[
g \in \im Q_{n_0} \setminus \cup_{n \in \sN_>} \im Q_n.
\]
Then, as above,  $\|Q_n g\| < 1$ for all $n \in \sN_>$. Consider the operator $K : x \mapsto \langle x, \, g \rangle f$ on $H$. Its adjoint is $K^* : x \mapsto \langle x, \, f \rangle g$, and
\[
P_n K Q_n K^* P_n x = \langle P_n x, \, f \rangle  \, \langle Q_n g, \, g \rangle \, P_n f = \langle x, \, P_n f \rangle \, \|Q_n g\|^2 P_n f.
\]
If $n \in \sN_<$, then $\|P_n f\| < 1$, and if $n \in \sN_>$, then $\|Q_n g\| < 1$ by construction. In both cases, $\|P_n K Q_n K^*
P_n\| < 1$. In case $n = n_0$,
\[
P_n K Q_n K^* P_n x = \langle x, \, f \rangle f
\]
is an orthogonal projection of rank 1, which we call $P$. The sequence $\bK := (P_n K Q_n K^* P_n)$ belongs to the algebra $\cJ$ since
\[
(P_n K Q_n K^* P_n) = (P_n KK^*P_n) - (P_n K P_n) \, (P_n K^*
P_n).
\]
As $r \to \infty$, the powers $\bK^r$ converge in the norm of $\cF^\cP$ to the sequence $(G_n)$ with $G_{n_0} = P \neq 0$ and $G_n = 0$ if $n \neq n_0$. Indeed, since $P_n \to I$ strongly, one has $\|Q_n g\| < 1/2$ for $n$ large enough, whence $\|P_n K Q_n K^* P_n\| < 1/2$ for these $n$, and for the remaining (finitely many) $n$ one has $\|P_n K Q_n K^* P_n\| < 1$ as we have seen above. Since $\bK^r \in \cJ$ and $\cJ$ is closed, the sequence $(G_n)$ has the claimed properties. \hfill \qed \\[3mm]
The stability of a sequence in $\cS(\eT(C))$ is related with its coset modulo $\cG := \cG^\cP$. So let us see what Theorem \ref{t14.21} tells us about the quotient algebra $\cS(\eT(C))/\cG$. Since the ideal $\cG$ lies in the kernel of the homomorphisms $W$ and $\widetilde{W}$, the mapping
\begin{equation} \label{e14.27}
\mbox{smb} : \cS(\eT(C))/\cG \to L(l^2(\sZ^+)) \times L(l^2(\sZ^+)), \quad \bA + \cG \mapsto (W(\bA), \, \widetilde{W}(\bA))
\end{equation}
is a well defined homomorphism. From Theorem \ref{t14.21} and (\ref{e14.23}), (\ref{e14.24}) we derive that the intersection of the kernels of $W$ and $\widetilde{W}$ is just the ideal $\cG$, which implies the following.
\begin{theo} \label{t14.29}
The mapping $\mbox{\rm smb}$ is a $^*$-isomorphism from $\cS(\eT(C))/\cG$ onto the $C^*$-subalgebra of $L(l^2(\sZ^+)) \times L(l^2(\sZ^+))$ which consists of all pairs $(W(\bA), \, \widetilde{W}(\bA))$ with $\bA \in \cS(\eT(C))$.
\end{theo}
\begin{coro} \label{t14.28}
A sequence $\bA \in \cS(\eT(C))$ is stable if and only if $\mbox{\rm smb}(\bA + \cG)$ is invertible in $L(l^2(\sZ^+)) \times L(l^2(\sZ^+))$.
\end{coro}
Indeed, by the inverse closedness of $C^*$-algebras, the coset $\bA + \cG \in \cS(\eT(C))/\cG$ is invertible in $\cF/\cG$ if and only if it is invertible in $\cS(\eT(C))/\cG$. So we arrived at a classical result:
\begin{coro} \label{c160309.2}
Let $a \in C(\sT)$ and $K$ compact. The finite sections sequence $\bA = (P_n (T(a) + K) P_n)$ is stable if and only of the operator $T(a) + K$ is invertible.
\end{coro}
Indeed, using some special properties of Toeplitz operators it is easy to see that the invertibility of $W(\bA) = T(a) + K$ already implies the invertibility of $\widetilde{W} (\bA) = T(\tilde{a})$.
\section{Fractality} \label{s2}
Clearly, a subsequence of a stable sequence is stable again. Does, conversely, the stability of a certain (infinite) subsequence of a sequence $\bA$ imply the stability of the full sequence? In general certainly not; but this implication holds indeed if $\bA$ belongs to the algebra $\cS(\eT(C))$ of the finite sections method for Toeplitz operators. The argument is simple: The homomorphisms $W$ and $\widetilde{W}$ defined in the previous section are given by certain strong limits. Thus, the operators $W(\bA)$ and $\widetilde{W} (\bA)$ can be computed if only a subsequence of $\bA$ is known. Moreover, if this subsequence is stable, then the operators $W(\bA)$ and $\widetilde{W} (\bA)$ are already invertible. This implies the stability of the full sequence $\bA$ via Corollary \ref{t14.28}.

Employing Theorem \ref{t14.29} instead of Corollary \ref{t14.28} we can state this observation in a slightly different way: every sequence in $\cS(\eT(C))$ can be rediscovered from each of its (infinite) subsequences up to a sequence tending to zero in the norm. In that sense, the essential information on a sequence in $\cS(\eT(C))$ is stored in each of its subsequences. Subalgebras of $\cF$ with this property were called {\em fractal} in \cite{RoS5} in order to emphasize exactly this self-similarity aspect. We will see some of the remarkable properties of fractal algebras in the following sections. We start with the official definition of fractal algebras.
\subsection{Fractal algebras} \label{ss22}
Let $\eta : \sN \to \sN$ be a strictly increasing sequence. By $\cF_\eta$ we denote the set of all subsequences $(A_{\eta(n)})$ of sequences $(A_n)$ in $\cF$. As in Section \ref{ss11}, one can make $\cF_\eta$ to a $C^*$-algebra in a natural way. The $^*$-homomorphism
\[
R_\eta : \cF \to \cF_\eta, \quad (A_n) \mapsto (A_{\eta(n)})
\]
is called the restriction of $\cF$ onto $\cF_\eta$. It maps the ideal $\cG$ of $\cF$ onto a closed ideal $\cG_\eta$ of $\cF_\eta$. For every subset $\cS$ of $\cF$, we abbreviate $R_\eta \cS$ by $\cS_\eta$.
\begin{defi}
Let $\cA$ be a $C^*$-subalgebra of $\cF$. A $^*$-homomorphism $W$ from $\cA$ into a $C^*$-algebra $\cB$ is called {\em fractal} \index{homomorphism!fractal} if, for every strictly increasing sequence $\eta : \sN \to \sN$, there is a mapping $W_\eta : \cA_\eta \to \cB$ such that $W = W_\eta R_\eta|_\cA$.
\end{defi}
Thus, the image of a sequence in $\cA$ under a fractal homomorphism can be reconstructed from each of its (infinite) subsequences. It is easy to see that $W_\eta$ is necessarily a $^*$-homomorphism again.
\begin{lemma} \label{lf1.4}
Let $\cA$ be a $C^*$-subalgebra of $\cF$, $\cB$ a $C^*$-algebra, and $W : \cA \to \cB$ a $^*$-homomorphism. The following assertions are equivalent: \\[1mm]
$(a)$ the homomorphism $W$ is fractal; \\[1mm]
$(b)$ every sequence $(A_n) \in \cA$ with $\liminf \|A_n\| = 0$ belongs to $\ker W$; \\[1mm]
$(c)$ $\ker (R_\eta|_\cA) \subseteq \ker W$ for every strictly increasing sequence $\eta : \sN \to \sN$. \\[1mm]
In particular, $\cA \cap \cG \subseteq \ker W$ for every fractal homomorphism $W$ on $\cA$.
\end{lemma}
{\bf Proof.} For the implication $(a) \Rightarrow (b)$, let $\bA = (A_n)$ be such that $\liminf \|A_n\| = 0$. Then $\|A_{\mu(n)}\| \to 0$ for a certain strictly increasing sequence $\mu$. Given $\varepsilon > 0$, choose $n_0$ such
that $\|A_{\eta(n)}\| \le \varepsilon$ for $n \ge n_0$. Put $\eta(n) := \mu(n + n_0)$. Then $\eta$ is an increasing sequence, so the fractality of $W$ implies that
\[
\|W(\bA)\| = \|(W_\mu R_\eta)(\bA)\| \le \|W_\eta\| \, \|R_\eta(\bG)\| \le \varepsilon.
\]
Since $\varepsilon$ is arbitrary, $W(\bA) = 0$. The implication $(b) \Rightarrow (c)$ is evident. For the implication $(c) \Rightarrow (a)$, let $\bA_1$ and $\bA_2$ be sequences in $\cA$ such that $R_\eta \bA_1 = R_\eta \bA_2$. Then
$\bA_1 - \bA_2 \in \ker (R_\eta|_\cA)$, which implies $W(\bA_1) = W(\bA_2)$ by condition $(c)$. Thus, the mapping $W_\eta : \cA_\eta \to \cB$, $R_\eta \bA \mapsto W(\bA)$ is correctly defined, and $W = W_\eta R_\eta|_\cA$. Hence, $W$ is fractal. \hfill \qed \\[3mm]
The following is the archetypal example of a fractal homomorphism. Let $\cP = (P_n)$ be a filtration on a Hilbert space $H$ and consider the $C^*$-algebra $\cF^\cP$ with consistency map $W : \cF^\cP \to L(H)$, $(A_n) \mapsto \mbox{s-lim}_{n \to \infty} A_n P_n$. The homomorphism $W$ is fractal, since the strong limit of a sequence is determined by each of its subsequences. Formally, given a strictly increasing sequence $\eta$, define
\[
W_\eta : \cF^\cP_\eta \to L(H), \quad (A_{\eta(n)}) \mapsto \mbox{s-lim}_{n \to \infty} A_{\eta(n)} P_{\eta(n)}.
\]
Then $W_\eta$ is a homomorphism and $W = W_\eta R_\eta$. \hfill \qed \\[3mm]
The fractal subalgebras of $\cF$ are distinguished by their property that every sequence in the algebra can be rediscovered from each of its (infinite) subsequences up to a sequence tending to zero. Here is the formal definition.
\begin{defi} \label{df1.5}
$(a)$ A $C^*$-subalgebra $\cA$ of $\cF$ is called {\em fractal} \index{algebra!fractal} if the canonical homomorphism $\pi : \cA \to \cA/(\cA \cap \cG), \; \bA \mapsto \bA + (\cA \cap \cG)$ is fractal. \\[1mm]
$(b)$ A sequence $\bA \in \cF$  is called {\em fractal} \index{sequence!fractal} if the smallest $C^*$-subalgebra of $\cF$ which contains the sequence $\bA$ and the identity sequence is fractal.
\end{defi}
Note that, by this definition, a fractal sequence always lies in a {\em unital} fractal algebra, whereas a fractal algebra needs not to be unital. The following results provide equivalent characterizations of fractal algebras.
\begin{theo} \label{tf1.8}
$(a)$ A $C^*$-subalgebra $\cA$ of $\cF$ is fractal if and only if the implication
\begin{equation} \label{ef1.9}
R_\eta (\bA) \in \cG_\eta \; \Rightarrow \; \bA \in \cG
\end{equation}
holds for every sequence $\bA \in \cA$ and every strictly increasing sequence $\eta$. \\[1mm]
$(b)$ If $\cA$ is a fractal $C^*$-subalgebra of $\cF$, then $\cA_\eta \cap \cG_\eta = (\cA \cap \cG)_\eta$ for every strictly increasing sequence $\eta$. \\[1mm]
$(c)$ If $\cA$ is a fractal $C^*$-subalgebra of $\cF$, then the algebra $\cA + \cG$ is fractal.
\end{theo}
{\bf Proof.} Assertion $(a)$ is an immediate consequence of Lemma \ref{lf1.4},
since the kernel of the canonical homomorphism $\cA \to \cA/(\cA \cap \cG)$ is
$\cA \cap \cG$. The inclusion $\supseteq$ in assertion $(b)$ is evident, and the
reverse inclusion follows from assertion $(a)$. Finally, for assertion $(c)$,
let $\bB \in \cA + \cG$ and $R_\eta (\bB) \in \cG_\eta$ for a certain strictly
increasing sequence $\eta$. Write $\bB$ as $\bA + \bG$ with $\bA \in \cA$ and
$\bG \in \cG$. Then $R_\eta (\bA) + R_\eta (\bG) \in \cG_\eta$, whence $R_\eta
(\bA) \in \cG_\eta$. Since $\cA$ is fractal, this implies $\bA \in \cG$ via
assertion $(a)$. Consequently, $\bB = \bA + \bG \in \cG$. So $\cA + \cG$ is
fractal, again by assertion $(a)$. \hfill \qed \\[3mm]
Thus, when working with fractal algebras $\cA$ one may always assume that $\cG \subset \cA$. Of course, the converse of assertion $(c)$ of the previous theorem is also true, i.e. if $\cA + \cG$ is fractal, then $\cA$ is fractal. This is a special case of assertion $(a)$ of the next result.
\begin{theo} \label{tf1.11}
$(a)$ $C^*$-subalgebras of fractal algebras are fractal. \\[1mm]
$(b)$ A $C^*$-subalgebra of $\cF$ is fractal if and only if each of its singly generated $C^*$-subalgebras is fractal. \\[1mm]
$(c)$ A $C^*$-subalgebra of $\cF$ which contains the identity sequence is fractal if and only if each of its elements is fractal. \\[1mm]
$(d)$ Restrictions of fractal algebras are fractal.
\end{theo}
{\bf Proof.} $(a)$ Let $\cA$ be a fractal $C^*$-subalgebra of $\cF$ and $\cB$ a $C^*$-subalgebra of $\cA$.  Let $\bB \in \cB$ and $R_\eta \bB \in \cG_\eta$ for a strictly increasing sequence $\eta$. Since $\cA$ is fractal, $\bB \in \cA \cap \cG$ by Theorem \ref{tf1.8} $(a)$. Hence $\bB \in \cB \cap \cG$, whence the fractality of $\cB$ by the same theorem. \\[1mm]
$(b)$ If $\cA$ is a fractal $C^*$-subalgebra of $\cF$, then each of its (singly generated or not) subalgebras is fractal by $(a)$. If $\cA$ is not fractal then, by Theorem \ref{tf1.8} $(a)$, there are a sequence $\bA \in \cA$ and a strictly increasing sequence $\eta$ such that $R_\eta (\bA) \in \cG_\eta$, but $\bA \not\in \cG$. By Theorem \ref{tf1.8} again, the subalgebra of $\cA$ singly generated by $\bA$ is not fractal. \\[1mm]
Assertion $(c)$ follows in the same way as $(b)$, and for $(d)$ note that each restriction of a restriction of a fractal algebra $\cA$ can be viewed as a restriction of $\cA$.  \hfill \qed \\[3mm]
The following theorem will offer a simple way to verify the fractality of many specific algebras of approximation methods, where the homomorphism $W$ appearing in the theorem is typically a direct product of fractal homomorphisms.
\begin{theo}\label{tf1.12}
A unital $C^*$-subalgebra $\cA$ of $\cF$ is fractal if and only if there is a unital and fractal $^*$-homomorphism $W$ from $\cA$ into a unital $C^*$-algebra $\cB$ such that, for every sequence $\bA \in \cA$, the coset $\bA + \cA \cap \cG$ is invertible in $\cA/(\cA \cap \cG)$ if and only if $W(\bA)$ is invertible in $\cB$.
\end{theo}
{\bf Proof.} If $\cA$ is fractal then the family $\{ \pi \}$, the only element
of which is the canonical homomorphism $\pi : \cA \to \cA/(\cA \cap \cG)$, has
the desired properties. Conversely, let $W$ be a $^*$-homomorphism which is
subject to the conditions of the theorem. Then $\cA \cap \cG \subseteq \ker W$
by Lemma \ref{lf1.4}, and the quotient map
\[
W^\pi : \cA/(\cA \cap \cG) \to \cB, \quad \bA + \cA \cap \cG \mapsto W(\bA)
\]
is correctly defined. Since $W^\pi$ preserves spectra, $W^\pi$ is a $^*$-isomorphism between $\cA/(\cA \cap \cG)$ and $W(\cA)$, and $(W^\pi)^{-1} W$ is equal to the canonical homomorphism $\pi : \cA  \to \cA/(\cA \cap \cG)$. Let now $\eta$ be a strictly increasing sequence. Since $W$ is fractal, $W = W_\eta R_\eta$ with a certain mapping $W_\eta$. Hence, $\pi = (W^\pi)^{-1} W_\eta R_\eta$ is fractal. \hfill \qed \\[3mm]
It is now easy to see that the algebra $\cS(\eT(C))$ is indeed a fractal subalgebra of $\cF$ in the formal sense of the definition, as expected. The fractality of the homomorphisms $W$ and $\widetilde{W}$ (both acting via strong limits) implies that
\[
\mbox{smb}^{(o)} : \cS(\eT(C)) \to L(l^2(\sZ^+)) \times L(l^2(\sZ^+)), \quad \bA \mapsto (W(\bA), \, \widetilde{W}(\bA))
\]
is a fractal mapping. Hence, for each strictly increasing sequence $\eta : \sN \to \sN$, there is a mapping $\mbox{smb}^{(o)}_\eta$ such that $\mbox{smb}^{(o)} = \mbox{smb}^{(o)}_\eta \circ R_\eta$. Further, from Theorem \ref{t14.29} we know that the mapping
\[
\mbox{smb} : \cS(\eT(C))/\cG \to L(l^2(\sZ^+)) \times L(l^2(\sZ^+)), \quad \bA + \cG \mapsto (W(\bA), \, \widetilde{W}(\bA))
\]
is an isomorphism. Hence, $\mbox{smb}^{-1} \circ \mbox{smb}^{(o)}_\eta \circ R_\eta$ is the canonical homomorphism from $\cS(\eT(C))$ onto $\cS(\eT(C))/\cG$.
\subsection{Consequences of fractality} \label{ss23}
The results in this section give a first impression of the power of fractality.
\begin{prop} \label{pf2.1}
Let $\cA$ be a unital fractal $C^*$-subalgebra of $\cF$. Then a sequence in $\cA$ is stable if and only if it possesses a stable subsequence.
\end{prop}
{\bf Proof.} Let $\bA = (A_n) \in \cA$, and let $\eta : \sN \to \sN$ be a strictly increasing sequence such that the sequence $R_\eta (\bA) = (A_{\eta(n)})$ is stable. One can assume without loss that $A_{\eta(n)}$ is invertible for every $n \in \sN$ (otherwise take a subsequence of $\eta$). Due to the inverse closedness of $\cA_\eta$ in $\cF_\eta$, there is a sequence $\bB \in \cA$ such that
\begin{equation} \label{ef2.2}
R_\eta (\bA) \, R_\eta (\bB) = R_\eta (\bB) \, R_\eta (\bA) = R_\eta (\bI)
\end{equation}
with $\bI$ the identity sequence. By hypothesis, the canonical homomorphism $\pi : \cA \to \cA/(\cA \cap \cG)$ factors into $\pi = \pi_\eta R_\eta$. Applying the homomorphism $\pi_\eta$ to (\ref{ef2.2}), we thus get the invertibility of $\pi (\bA) = \bA + (\cA \cap \cG)$ in $\cA/(\cA \cap \cG)$. Hence, $\bA$ is a stable sequence. The reverse implication is obvious. \hfill \qed
\begin{prop} \label{pf2.3}
Let $\cA \subset \cF$ be a fractal $C^*$-algebra and $\bA = (A_n) \in \cA$. Then, \\[1mm]
$(a)$ for every strictly increasing sequence $\eta : \sN \to \sN$,
\[
\|\bA + \cG\|_{\cF/\cG} = \|R_\eta (\bA) + \cG_\eta\|_{\cF_\eta/\cG_\eta}.
\]
$(b)$ the limit $\lim_{n \to \infty} \|A_n\|$ exists and is equal to $\|\bA + \cG\|$.
\end{prop}
{\bf Proof.} $(a)$ By the third isomorphy theorem,
\begin{equation} \label{ef2.4}
\|\bA + \cG\|_{\cF/\cG} = \|\bA + \cG\|_{(\cA + \cG)/\cG} = \|\bA + (\cA \cap \cG)\|_{\cA/(\cA \cap \cG)}
\end{equation}
for each sequence $\bA$ in a (not necessarily fractal) $C^*$-subalgebra $\cA$ of $\cF$. If $\cA$ is fractal and $\bG \in \cA \cap \cG$, then
\[
\begin{array}{rclll}
\|\bA + \cG\|_{\cF/\cG} & = & \|\pi (\bA + \bG)\|_{\cA/(\cA \cap \cG)} & &(\mbox{by} \, (\ref{ef2.4})) \\[1mm]
& = & \|\pi_\eta R_\eta (\bA + \bG)\|_{\cA/(\cA \cap \cG)} & & (\mbox{fractality of} \; \pi) \\[1mm]
& \le & \|R_\eta (\bA + \bG)\|_{\cA_\eta}. & &
\end{array}
\]
Taking the infimum over all sequences $\bG \in \cA \cap \cG$ and applying Theorem \ref{tf1.8} $(b)$, we obtain
\begin{eqnarray*}
\|\bA + \cG\|_{\cF/\cG} & \le & \|R_\eta (\bA) + (\cA \cap \cG)_\eta\|_{\cA_\eta/(\cA \cap \cG)_\eta} \\
& = & \|R_\eta (\bA) + (\cA_\eta \cap \cG_\eta)\|_{\cA_\eta/(\cA_\eta \cap \cG_\eta)} \\
& = & \|R_\eta (\bA) + \cG_\eta\|_{\cF_\eta/\cG_\eta}
\end{eqnarray*}
where we used (\ref{ef2.4}) again. The reverse estimate is a consequence of the $\limsup$-formula (\ref{e170309.3}):
\[
\|(A_{\eta(n)}) + \cG_\eta\|_{\cF_\eta/\cG_\eta} = \limsup \|A_{\eta(n)}\| \le \limsup \|A_n\| = \|(A_n) + \cG\|_{\cF/\cG}.
\]
$(b)$ Choose a strictly increasing sequence $\eta : \sN \to \sN$ such that $\lim \|A_{\eta(n)}\| = \liminf \|A_n\|$.
By part $(a)$ of this proposition and by (\ref{e170309.3}),
\begin{eqnarray*}
\lefteqn{\limsup \|A_n\| = \|(A_n) + \cG\|_{\cF/\cG} = \|(A_{\eta(n)}) + \cG_\eta\|_{\cF_\eta/\cG_\eta} } \\
& & = \limsup \|A_{\eta(n)}\| = \lim \|A_{\eta(n)}\| = \liminf \|A_n\|
\end{eqnarray*}
which gives the assertion.  \hfill \qed \\[3mm]
Our next goal is convergence properties of the spectra $\sigma (A_n)$ for fractal normal sequences $(A_n)$ (there is no hope to say something substantial in case the sequence $(A_n)$ is not normal). We will need the following notions.

Let $(M_n)_{n \in \sN}$ be a sequence in the set $\sC_{comp}$ of the non-empty compact subsets of the complex plane. The {\em limes superior} $\limsup M_n$ \index{limes superior} $($also called the {\em partial limiting set}$)$ resp. the {\em limes inferior} $\liminf M_n$ \index{limes inferior} $($or the {\em uniform limiting set}$)$ of the sequence $(M_n)$ consists of all points $x \in \sC$ which are a partial limit resp. the limit of a sequence $(m_n)$ of points $m_n \in M_n$. Observe that both $\limsup M_n$ and $\liminf M_n$ are closed sets. The partial limiting set $\limsup M_n$ is never empty if $\cup_n M_n$ is bounded, whereas the uniform limiting set of a bounded set sequence can be empty.

The following result is the analog of the limsup formula (\ref{e170309.3}) for norms.
\begin{prop} \label{t16.4}
Let $(A_n) \in \cF$ be a normal sequence. Then
\[
\limsup \sigma (A_n) = \sigma_{\cF/\cG} ((A_n) + \cG).
\]
\end{prop}
{\bf Proof.} Let $\lambda \in \sigma ((A_n) + \cG)$. Then $(A_n - \lambda I_n)$ is a stable sequence. Since the norm $\|M\|$ is not less than the spectral radius $\rho(M)$ of a matrix $M$, there is an $n_0 \in \sN$ such that
\[
\sup_{n \ge n_0} \rho ((A_n - \lambda I_n)^{-1}) =: m < \infty.
\]
Then, for all $n \ge n_0$,
\[
m \ge \sup \, \{ |t| : t \in \sigma((A_n - \lambda I_n)^{-1}) \} = \sup \, \{ |t|^{-1} : t \in \sigma(A_n - \lambda I_n) \}
\]
whence
\[
1/m \le \inf \, \{ |t| : t \in \sigma(A_n) - \lambda \} = \inf \, \{ |t - \lambda| : t \in \sigma(A_n) \}
\]
for all $n \ge n_0$. Hence, $\lambda$ cannot belong to $\limsup \sigma(A_n)$.

For the reverse inclusion assume the sequence $(A_n - \lambda I_n)$ fails to be stable. Then either there is an infinite subsequence $(A_{n_k} - \lambda I_{n_k})$ which consists of non-invertible matrices only, or all matrices $A_n - \lambda I_n$ with sufficiently large $n$ are invertible, but
\[
\|(A_{n_k} - \lambda I_{n_k})^{-1}\| = \rho ((A_{n_k} - \lambda I_{n_k})^{-1}) \to \infty \quad \mbox{as} \; k \to \infty
\]
for a certain subsequence (note that $\|M\| = \rho(M)$ for every normal matrix $M$). In the first case one has $\lambda \in \sigma(A_{n_k})$ for every $k$, whence $\lambda \in \limsup \sigma(A_n)$. In the second case one finds numbers $t_{n_k} \in \sigma(A_{n_k})$ such that $|t_{n_k} - \lambda|^{-1} \to \infty$ resp. $|t_{n_k} - \lambda| \to 0$ as $k \to \infty$ which implies $\lambda \in \limsup \sigma(A_n)$ also in this case.  \hfill \qed
\\[3mm]
In the context of fractal algebras one can say again more. Note that the equality $\liminf M_n = \limsup M_n$ for a sequence of non-empty compact subsets of $\sC$ is equivalent to the convergence of the sequence $(M_n)$ with respect to the Hausdorff metric
\[
h(L, \, M) := \max \{ \max_{l \in L} \dist (l, \, M), \, \max_{m \in M} \dist (m, \, L) \}.
\]
Recall in this connection that $(\sC_{comp}, \, h)$ is a complete metric space and that every bounded sequence in $(\sC_{comp}, \, h)$ possesses a convergent subsequence. Thus, the relatively compact subsets of the metric space $(\sC_{comp}, \, h)$ are precisely its bounded subsets.
\begin{prop} \label{pf2.7}
Let $\cA$ be a fractal unital $C^*$-subalgebra of $\cF$. If $(A_n) \in \cA$ is normal, then
\begin{equation} \label{ef2.8}
\limsup \sigma (A_n) = \liminf \sigma (A_n) \quad (= \sigma_{\cF/\cG} ((A_n) + \cG)).
\end{equation}
\end{prop}
{\bf Proof.} Let $\lambda \in \sC \setminus \liminf \sigma(A_n)$. Then there are a $\delta > 0$ and a strictly
increasing sequence $\eta : \sN \to \sN$ such that $\dist (\lambda, \, \sigma(A_{\eta(n)})) \ge \delta$ for all $n$. Thus, and since the $A_n$ are normal,
\[
\sup_n \|(A_{\eta(n)} - \lambda I_{\eta(n)})^{-1}\| = \sup_n \rho ((A_{\eta(n)} - \lambda I_{\eta(n)})^{-1})
< 1/\delta.
\]
This shows that the sequence $(A_{\eta(n)} - \lambda I_{\eta(n)})$ is stable. Then, by Proposition \ref{pf2.1}, the sequence $(A_n - \lambda I_n)$ itself is stable. Hence, $\lambda \not\in \sigma ((A_n) + \cG) = \limsup \sigma(A_n)$ by Proposition \ref{t16.4}, which gives $\limsup \sigma(A_n) \subseteq \liminf \sigma(A_n)$. The reverse inclusion is evident. \hfill \qed \\[3mm]
Results as in Propositions \ref{pf2.3} and \ref{pf2.7} can be derived also for other spectral quantities, for example for the sequences of the condition numbers, the sets of the singular values, the $\epsilon$-pseudospectra, and the numerical ranges of the $A_n$. For details see Chapter 3 in \cite{HRS2}.
\subsection{Fractal restrictions of separable algebras} \label{ss24}
The results in the previous section indicate that, given an (in general non-fractal) approximation sequence, it is of vital importance to single out (one of) its fractal subsequences. In this moment it is not clear whether such subsequences exist at all, and how one can find them. Whereas the existence of a fractal subsequence of an approximation sequence will be established in this section in full generality, we are able to construct such sequences only in specific situations.

Our starting point is the simple but useful observation that the converse of Proposition \ref{pf2.3} $(b)$ is also true: If, for every sequence $(A_n)$ in a $C^*$-subalgebra $\cA$ of $\cF$, the sequence of the norms $\|A_n\|$ converges, then the algebra $\cA$ is fractal. Moreover, it is sufficient to have the convergence of the norms for a dense subset of $\cA$.
\begin{prop} \label{pf3.0a}
Let $\cA$ be a $C^*$-subalgebra of $\cF$ and $\cL$ a dense subset of $\cA$. If the sequence of the norms $\|A_n\|$ converges for each sequence $(A_n) \in \cL$, then the algebra $\cA$ is fractal.
\end{prop}
Note that already the convergence of the norms $\|A_n\|$ for each sequence $(A_n)$ in a dense subset of $\cA_{sa}$, \index{$\cA_{sa}$} the set of the self-adjoint elements of $\cA$, is sufficient for the fractality of $\cA$, which follows immediately from the $C^*$-axiom. By the same reason, the convergence of the norms for each sequence in a dense subset of the set of all positive elements of $\cA$ is sufficient for fractality. \\[3mm]
{\bf Proof.} First we show that if the sequence of the norms converges for each sequence in $\cL$ then it converges for each sequence in $\cA$. Let $(A_n) \in \cA$ and $\varepsilon > 0$. Choose $(L_n) \in \cL$ such that $\|(A_n) - (L_n)\|_\cF = \sup \|A_n - L_n\| < \varepsilon/3$, and let $n_0 \in \sN$ be such that $\left| \|L_n\| - \|L_m\| \right| < \varepsilon/3$ for all $m, \, n \ge n_0$. Then, for $m, \, n \ge n_0$,
\begin{eqnarray*}
\left| \|A_n\| - \|A_m\| \right| & \le &  \left| \|A_n\| - \|L_n\| \right| + \left| \|L_n\| - \|L_m\| \right| + \left| \|L_m\| - \|A_m\| \right| \\
& \le & \|A_n - L_n\| + \left| \|L_n\| - \|L_m\| \right| + \|L_m - A_m\| \; \le \; \varepsilon.
\end{eqnarray*}
Thus, $(\|A_n\|)$ is a Cauchy sequence, hence convergent. But the convergence of the norms for each sequence in $\cA$ implies the fractality of $\cA$ by Theorem \ref{tf1.8}. Indeed, if a subsequence of a sequence $(A_n) \in \cA$ tends to zero, then $0 = \liminf \|A_n\| = \lim \|A_n\|$, whence $(A_n) \in \cG$. \hfill \qed \\[3mm]
The following {\em fractal restriction theorem} was first proved in \cite{Roc0}. The proof given there was based on the converse of Proposition \ref{pf2.7} and rather involved. The surprisingly simple proof given below is based on Proposition \ref{pf3.0a} instead.
\begin{theo} \label{tf3.20}
Let $\cA$ be a separable $C^*$-subalgebra of $\cF$. Then there exists a strictly increasing sequence $\eta : \sN \to \sN$ such that the restricted algebra $\cA_\eta = R_\eta \cA$ is a fractal subalgebra of $\cF_\eta$.
\end{theo}
Since finitely generated $C^*$-algebras are separable, this result immediately implies:
\begin{coro} \label{cf3.21}
Every sequence in $\cF$ possesses a fractal subsequence.
\end{coro}
One cannot expect that Theorem \ref{tf3.20} holds for arbitrary $C^*$-subalgebras of $\cF$; for example it is certainly not valid for the algebra $l^\infty$. On the other hand, there {\em are} of course non-separable but fractal algebras; the algebra $\cS(\eT(PC))$ of the finite sections method for Toeplitz operators with piecewise continuous generating function can serve as an example. \\[3mm]
{\bf Proof of Theorem \ref{tf3.20}.} Let $\{\bA^m\}_{m \in \sN}$ be a dense countable subset of $\cA$ which consists of sequences $\bA^m = (A_n^m)_{n \in \sN}$. Let $\eta_1 : \sN \to \sN$ be a strictly increasing sequence such that the sequence of the norms $\|A_{\eta_1(n)}^1\|$ converges. Next let $\eta_2$ be a strictly increasing subsequence of $\eta_1$ such that the sequence $(\|A_{\eta_2(n)}^2\|)_{n \in \sN}$ converges. We proceed in this way and find, for each $k \ge 2$, a strictly increasing subsequence $\eta_k$ of $\eta_{k-1}$ such that the sequence $(\|A_{\eta_k(n)}^k\|)_{n \in \sN}$ converges. Set $\eta(n) := \eta_n(n)$ for $n \in \sN$. Then $\eta$ is a strictly increasing sequence, and the sequence $(\|A_{\eta(n)}^k\|)_{n \in \sN}$ converges for every $k \in \sN$.

Since the sequences $R_\eta (\bA^m)$ with $k \in \sN$ form a dense subset of the restricted algebra $\cA_\eta$ and since each sequence $R_\eta (\bA^m) = (A_{\eta(n)}^k)_{n \in \sN}$ has the property that the sequence of the norms $\|A_{\eta(n)}^k \|$ converges, the assertion follows from Proposition \ref{pf3.0a}. \hfill \qed
\subsection{Spatial discretization of Cuntz algebras} \label{s6}
We will digress for a moment in order to illustrate the usefulness of the fractal restriction theorem in the analysis of discretizations of concrete operator algebras. Our running example, the Toeplitz algebra $\eT(C)$, is (isomorphic to) the universal algebra generated by one isometry. Now we go one step further and consider the discretization of algebras which are generated by a finite number of non-commuting non-unitary isometries, namely the Cuntz algebras.

Let $N \ge 2$. The {\em Cuntz algebra} $\cO_N$ is the universal $C^*$-algebra generated by $N$ isometries $s_0, \, \ldots, \, s_{N-1}$ with the property that
\begin{equation} \label{e91.2}
s_0 s_0^* + \ldots + s_{N-1} s_{N-1}^* = I.
\end{equation}
For basic facts about Cuntz algebras see Cuntz' pioneering paper \cite{Cun1}. A nice introduction is also in \cite{Dav1}. The importance of Cuntz algebras in theory and applications cannot be overestimated. Let me only mention Kirchberg's deep result that a separable $C^*$-algebra is exact if and only if it embeds in the Cuntz algebra $\cO_2$, and the role that representations of Cuntz algebras play in wavelet theory and signal processing (see \cite{Bla2,BrJ1} and the references therein).

To discretize the Cuntz algebra $\cO_N$ by the finite sections method, we
represent this algebra as an algebra of operators on  $l^2(\sZ^+)$. Since Cuntz
algebras are simple, every $C^*$-algebra which is generated by $N$ isometries
$S_0, \, \ldots, \, S_{N-1}$ which fulfill (\ref{e91.2}) in place of the $s_i$
is $^*$-isomorphic to $\cO_N$. Thus, $\cO_N$ is $^*$-isomorphic to the smallest
$C^*$-subalgebra of $L(l^2(\sZ^+))$ which contains the operators
\begin{equation} \label{e91.3}
S_i : (x_k)_{k \ge 0} \mapsto (y_k)_{k \ge 0}
\quad \mbox{with} \quad
y_k := \left\{
\begin{array}{ll}
x_r & \mbox{if} \; k = rN + i \\
0   & \mbox{else}
\end{array}
\right.
\end{equation}
for $i = 0, \, \ldots, \, N-1$. We denote the (concrete) Cuntz algebra generated by the operators $S_i$ in  (\ref{e91.3}) by $\eO_N$. We use the same the filtration $\cP = (P_n)$ as for the Toeplitz algebra $\eT(C)$ and consider the smallest closed subalgebra $\cS(\eO_N)$ of $\cF = \cF^\cP$ which contains all finite sections sequences $(P_nAP_n)$ with $A \in \eO_N$. Since $(P_n A P_n)^* = (P_n A^* P_n)$, $\cS(\eO_N)$ is a $C^*$-algebra.

One should mention that the abstract Cuntz algebra $\cO_N$ has an uncountable set of equivalence classes of irreducible representations. Representations of $\cO_N$ different from (\ref{e91.3}) will certainly lead to sequence algebras different from $\cS(\eO_N)$. The relation between these algebras is not yet understood. The chosen representation of $\cO_N$ is distinguished by the fact that it is both irreducible and {\em permutative} in the sense that every isometry $S_i$ maps elements of the standard basis to elements of the standard basis.

The algebra $\cS(\eT(C))$ of the finite sections method for Toeplitz operators is the smallest closed $C^*$-subalgebra of $\cF$ which contains the sequence $(P_n V_1 P_n)$. A similar description holds for the algebra $\cS(\eO_N)$. Set  $\Omega := \{0, \, 1, \, \ldots, \, N-1\}$.
\begin{lemma} \label{l92.1}
$\cS(\eO_N)$ is the smallest $C^*$-subalgebra of $\cF$ which contains all sequences $(P_n S_j P_n)$ with $j \in \Omega$.
\end{lemma}
{\bf Proof.} Let $\cS^\prime$ denote the smallest $C^*$-subalgebra of $\cF$ which contains all sequences $(P_n S_j P_n)$ with $j \in \Omega$. Evidently, $\cS^\prime \subseteq \cS(\eO_N)$. For the reverse inclusion, note that
\begin{equation} \label{e92.2}
S_i^* S_j = 0 \quad \mbox{whenever} \; i \neq j.
\end{equation}
Indeed, this follows by straightforward calculation, but it also a consequence of the Cuntz axiom (\ref{e91.2}): Multiply (\ref{e91.2}) from the left by $S_i^*$ and from the right by $S_i$ and take into account that a sum of positive elements in a $C^*$-algebra is zero if and only if each of the elements is zero. From (\ref{e92.2}) we conclude that every finite word with letters in the alphabet $\{ S_1, \, \ldots, \, S_N, \, S_1^*, \, \ldots,
S_N^* \}$ is of the form
\begin{equation} \label{e92.3}
S_{i_1} S_{i_2} \ldots S_{i_k} S^*_{j_1} S^*_{j_2} \ldots S^*_{j_l}
\quad \mbox{with} \quad i_s, \, j_t \in \Omega
\end{equation}
(Lemma 1.3 in \cite{Cun1}). Further one easily checks that
\begin{equation} \label{e92.4}
P_n S_j = P_n S_j P_n \quad \mbox{and} \quad S_j^* P_n = P_n S_j^* P_n
\end{equation}
for every $j \in \Omega$ and every $n \in \sN$. Thus, if $A$ is
any word of the form (\ref{e92.3}), then
\[
P_n A P_n = P_n S_{i_1} P_n \cdot P_n S_{i_2} P_n \ldots P_n S_{i_k} P_n \cdot P_n S^*_{j_1} P_n \cdot P_n S^*_{j_2} P_n \ldots P_n S^*_{j_l} P_n \in \cS^\prime.
\]
Since the set of all linear combinations of the words (\ref{e92.3}) is dense in $\eO_N$, it follows that $\cS(\eO_N)
\subseteq \cS^\prime$. \hfill \qed \\[3mm]
For a closer look at the generators of $\cS(\eO_N)$, recall that an element $S$ of a $C^*$-algebra is called a partial isometry if $SS^*S = S$. If $S$ is a partial isometry, then $SS^*$ and $S^*S$ are projections (i.e., self-adjoint idempotents), called the {\em range projection} and the {\em initial projection} of $S$, respectively. Conversely, if $S^*S$ (or $SS^*$) is a projection for an element $S$, then $S$ is a partial isometry. Recall also that projections $P$ and $Q$ are called orthogonal if $PQ = 0$.
\begin{lemma} \label{l92.5}
Every sequence $(P_n S_i P_n)$, $i \in \Omega$, is a partial isometry in $\cF$, and the corresponding range projections are orthogonal if $i \neq j$. Moreover,
\begin{equation} \label{e92.6}
P_n S_i^* P_n S_j P_n = 0 \quad \mbox{if} \; \; i \neq j,
\end{equation}
and
\begin{equation} \label{e92.7}
P_n S_0 P_n S_0^* P_n + \ldots + P_n S_{N-1} P_n S_{N-1}^* P_n = P_n.
\end{equation}
\end{lemma}
{\bf Proof.}  The identities (\ref{e92.4}) imply that $P_n S_i S_i^* P_n = P_n S_i P_n S_i^* P_n$ for every $i \in \Omega$ and every $n \in \sZ^+$. The operators $S_i S_i^*$ are projections, and their matrices with respect to the standard basis of $l^2(\sZ^+)$ are of diagonal form. Hence, the projections $S_i S_i^*$ and $P_n$ commute, which implies that $P_n S_i S_i^* P_n$ is a projection. Hence, $(P_n S_i P_n)$ is a partial isometry in $\cF$, and $(P_n S_i S_i^* P_n)$ is the associated range projection.

Let $i \neq j$ be in $\Omega$. The fact that $P_n$ and $S_i S_i^*$ commute further implies together with (\ref{e92.2}) that
\[
(P_n S_i S_i^* P_n) (P_n S_j S_j^* P_n) = P_n S_i S_i^* S_j S_j^* P_n = 0.
\]
Multiplying $P_n S_i S_i^* P_n S_j S_j^* P_n = 0$ from the left by $P_n S_i^* P_n$ and from the right by $P_n S_j P_n$ yields (\ref{e92.6}). Finally, (\ref{e92.7}) follows by summing up the equalities (\ref{e92.6}) over $i \in \Omega$ and from axiom (\ref{e91.2}). \hfill \qed \\[3mm]
Thus, the generating sequences $(P_n S_i P_N)$, $i \in \Omega$, are still subject of the Cuntz axiom (\ref{e91.2}), but note they are partial isometries only and no longer isometries. Next we look at products of these generating sequences. For $i = (i_1, \, i_2, \, \ldots, \, i_k) \in \Omega^k$, abbreviate $S_i := S_{i_1} S_{i_2} \ldots S_{i_k}$. Further, for every real number $x$, let $\{x\}$ denote the smallest integer which is greater than or equal to $x$. The first assertion of the following proposition follows as in Lemma \ref{l92.5}, the second one by straightforward calculation.
\begin{prop} \label{p92.10}
Let $i = (i_1, \, i_2, \, \ldots, \, i_k) \in \Omega^k$. Each sequence
\[
(P_n S_{i_1} P_n S_{i_2} P_n \ldots P_n S_{i_k} P_n)
\]
is a partial isometry in $\cF$. The corresponding range projection is given by
\begin{equation} \label{e92.11}
P_n S_i^* P_n S_i P_n = P_{\{ (n - v_{i,k})/N^k \} },
\end{equation}
where $v_{i, k} := i_1 + i_2 N + \ldots + i_k N^{k-1}$.
\end{prop}
We specialize the result of Proposition \ref{p92.10} to the case $k = 1$. If $n
= jN$ is a multiple of $N$, then $\{ (n-i)/N \} = \{ (j N-i)/N \} = \{ j - i/N
\} = j$, whence
\begin{equation} \label{e92.12}
P_{jN} S_i^* P_{jN} S_i P_{jN} = P_j \quad \mbox{for all} \; i \in \Omega.
\end{equation}
On the other hand, one has
\begin{equation} \label{e92.13}
P_n S_0^* P_n S_0 P_n - P_n S_1^* P_n S_1 P_n = \left\{
\begin{array}{ll}
P_{j+1} - P_j & \mbox{if} \; n = jN+1, \\
0             & \mbox{else}.
\end{array}
\right.
\end{equation}
Thus, the sequence
\begin{equation} \label{e92.14}
(P_n S_0^* P_n S_0 P_n - P_n S_1^* P_n S_1 P_n)_{n \ge 1}
\end{equation}
possesses both a subsequence consisting of zeros only (take $\eta(n) := nN$) and a subsequence consisting of non-zero
projections (if $\eta(n) := nN + 1$). This shows that the algebra $\cS(\eO_N)$ cannot be fractal!

This is the place where the idea of forcing fractality by restriction comes into play. For simplicity (a sequence of zeros looks simpler than a sequence of non-zero projections) one would like to choose $\eta(n) := nN$ and consider the restricted algebra $\cS_\eta(\eO_N)$. Unfortunately, a similar argument shows that also this restricted algebra $\cS_\eta(\eO_N)$ cannot be fractal. This time one is attempted to choose $\eta(n) := nN^2$, which again leads to similar problems. Further manipulations with the generating sequences convinced me that one can avoid the above mentioned problems if one chooses
\begin{equation} \label{e92.16}
\eta(n) := N^n.
\end{equation}
It turned out that this is indeed the right choice, i.e., the algebra $\cS_\eta(\eO_N)$ with $\eta$ specified by (\ref{e92.16}) is fractal. The full approach is in \cite{Roc11}. Here I will add only a few comments.

Coburn's already mentioned result suggests to consider the Toeplitz algebra $\eT(C)$ as the Cuntz algebra $\eO_1$. But one should have in mind that the main properties of $\eO_1$ and of $\eO_N$ for $N > 1$ are quite different from each other. For example, the compact operators $K(l^2(\sZ^+))$ form a closed ideal of $\eO_1$, and the quotient $\eO_1/K(l^2(\sZ^+))$ is isomorphic to $C(\sT)$, whereas $\eO_N$ is simple if $N \ge 2$. These differences continue to the corresponding sequence algebras $\cS(\eO_1)$ and $\cS(\eO_N)$ for $N > 1$. A main point is that $\cS(\eO_1)/\cG$ contains two ideals, each of which isomorphic to the ideal $K(l^2(\sZ^+))$, and that the irreducible representations $W$ and $\widetilde{W}$ of $\cS(\eO_1)$ coming from these ideals own the property that a sequence $\bA = (A_n)$ in $\cS(\eO_1)$ is stable if and only if $W(\bA)$ and $\widetilde{W}(\bA)$ are invertible. We have seen that this fact implies an effective criterion to check the stability of a sequence in $\cS(\eO_1)$.

In contrast to these facts, if $N > 1$, then $\cS_\eta(\eO_N)/\cG$ has only one non-trivial ideal. There is an injective representation of this ideal, which extends to a representation, $W^\prime$ say, of $\cS_\eta(\eO_N)/\cG$ which is injective on all of $\cS(\eO_N)/\cG$. Thus, roughly speaking, the stability result of \cite{Roc11} says that a sequence $\bA$ in $\cS(\eO_N)$ is stable if and only if the operator $W^\prime(\bA)$ is invertible. At the first glance, this result might seem to be quite useless since the stability of $\bA$ is not easier to check than the invertibility of $W^\prime(\bA)$. So why this effort, if many canonical homomorphisms on $\cS_\eta(\eO_N)/\cG$ own the same property as $W^\prime$: the identical mapping and the faithful representation via the GNS-construction, for example. What is important is the concrete form of the mapping $W^\prime$: it can be defined by means of strong limits of operator sequences, and this special form implies the fractality of $\cS_\eta(\eO_N)$ immediately.
\subsection{Fractality of self-adjoint sequences} \label{ssf3.1}
Let again $\cF$ be an algebra of matrix sequences (or, more generally, the product of a sequence $(\cC_n)_{n \in \sN}$ of unital $C^*$-algebras) with associated ideal $\cG$ of zero sequences. If $(A_n) \in \cF$ is a normal fractal sequence then
\begin{equation} \label{ef3.1}
\limsup \sigma(A_n) = \liminf \sigma(A_n)
\end{equation}
by Proposition \ref{pf2.7}. We will see now that (\ref{ef3.1}) is the {\em only} obstruction for a normal sequence to be fractal.
\begin{theo} \label{tf3.2}
A normal sequence $(A_n) \in \cF$ is fractal if and only if $(\ref{ef3.1})$ holds.
\end{theo}
{\bf Proof.} The 'only if'-part is Proposition \ref{pf2.7}. For the reverse conclusion suppose that (\ref{ef3.1}) holds. Let $\cA$ denote the smallest closed subalgebra of $\cF$ which contains the sequences $(A_n)$ and $(A_n^*)$ and the identity sequence $(I_n)$. Further, let $\eta : \sN \to \sN$ be strictly increasing. We abbreviate $R_\eta \cA$ to  $\cA_\eta$ and claim that the mapping
\begin{equation} \label{ef3.4x}
(B_n) + (\cA \cap \cG) \mapsto (B_{\eta(n)}) + (\cA_\eta \cap \cG_\eta).
\end{equation}
establishes a $^*$-isomorphism
\begin{equation} \label{ef3.3}
\cA/(\cA \cap \cG) \cong \cA_\eta/(\cA_\eta \cap \cG_\eta).
\end{equation}
By the third isomorphy theorem, the algebras $\cA/(\cA \cap \cG)$ and $\cA_\eta/(\cA_\eta \cap \cG_\eta)$ are $^*$-isomorphic to $(\cA + \cG)/\cG$ and $(\cA_\eta + \cG_\eta)/\cG_\eta$, respectively. The latter are (as unital $C^*$-algebras) singly generated by their elements $(A_n) + \cG$ and $(A_{\eta(n)}) + \cG_\eta$, and the spectra of these cosets are $\limsup \sigma(A_n)$ and $\limsup \sigma(A_{\eta(n)})$, respectively, by Proposition \ref{t16.4}. The assumption (\ref{ef3.1}) guarantees that these spectra coincide. Hence, (\ref{ef3.3}) is a consequence of the Gelfand-Naimark theorem for singly generated $C^*$-algebras.

Let now $\pi$ stand for the canonical homomorphism $\cA \to \cA/(\cA \cap \cG)$. Then $\pi = \varphi_\eta \psi_\eta R_\eta$ where $\psi_\eta$ is the canonical homomorphism from $\cA_\eta$ onto $\cA_\eta/ (\cA_\eta \cap \cG_\eta)$ and $\varphi_\eta$ if the inverse of the isomorphism (\ref{ef3.4x}). Hence, $\pi$ is fractal. \hfill \qed \\[3mm]
The previous theorem offers an alternate way to verify Corollary \ref{cf3.21} for {\em normal} sequences.
\begin{coro} \label{tf3.19}
Every normal sequence in $\cF$ possesses a fractal subsequence.
\end{coro}
{\bf Proof.} Let $(A_n) \in \cF$ be normal and abbreviate $M_n := \sigma (A_n)$. Since the sequence $(M_n)$ is bounded in $\sC$, it has a subsequence $(M_{\eta(n)})_{n \in \sN}$ which converges with respect to the Hausdorff metric. Hence, $\limsup (M_{\eta(n)}) = \liminf (M_{\eta(n)})$, and Theorem \ref{tf3.2} implies the fractality of the sequence
$(A_{\eta(n)})_{n \ge 1}$. \hfill \qed \\[3mm]
As another application of Theorem \ref{tf3.2}, we derive a fractality criterion for the sequence of the finite sections of a self-adjoint operator. Note in this connection that the spectrum of a self-adjoint operator $A$ for which the finite sections method $(P_n A P_n)$ is fractal can be as complicated as possible: Given a non-empty compact subset $K$ of the real line, choose a sequence $(k_n)_{n \in \sN}$ which is dense in $K$ and consider the diagonal operator $A := \diag (k_1, \, k_2, \, k_3, \, \dots )$ on $l^2(\sN)$. Then
\[
\limsup \sigma(P_n A P_n) = \liminf \sigma(P_n A P_n) = K.
\]
Hence, the sequence $(P_n A P_n)$ is fractal by Theorem \ref{tf3.2}.
\begin{theo} \label{tf3.5}
Let $A \in L(H)$ be a self-adjoint operator with connected spectrum and let $\cP = (P_n)$ be a filtration on $H$. Then the sequence $(P_n A P_n)$ is fractal.
\end{theo}
We consider the sequence $(P_n A P_n)$ as an element of the algebra $\cF^\cP$. The proof of Theorem \ref{tf3.5} is based on the following result by Arveson.
\begin{theo} \label{tf3.6}
Let $(A_n)$ be a normal sequence in $\cF^\cP$ with strong limit $A$. Then
\begin{equation} \label{ef3.7}
\sigma (A) \subseteq \liminf \sigma(A_n).
\end{equation}
\end{theo}
In particular, $\liminf \sigma(A_n)$ is not empty for a $^*$-strongly convergent normal sequence $(A_n)$. \\[3mm]
{\bf Proof.} Let $\lambda \in \sC \setminus \liminf \sigma(A_n)$. We have to show that $A - \lambda I$ is invertible. Since $\lambda$ is not in the lower limit of the spectra, there are an $\varepsilon > 0$ and an infinite subset $\sM$ of $\sN$ such that the distance of $\sigma(A_n)$ to $\lambda$ is at least $\varepsilon$ for each $n \in \sM$. Since the $A_n$ are normal, this implies that the operators $A_n - \lambda I_n|_{\im P_n}$ are invertible and that their inverses are uniformly bounded,
\begin{equation} \label{ef3.8}
\sup_{n \in \sM} \|(A_n - \lambda I_n|_{\im P_n})^{-1}\| \le 1/\varepsilon.
\end{equation}
The strong convergence of $(A_n - \lambda I_n) P_n$ to $A - \lambda I$ together with the uniform boundedness (\ref{ef3.8}) imply the strong convergence of the sequence of the inverses $(A_n - \lambda I_n|_{\im P_n})^{-1} P_n$. Write $B$ for the strong limit of this sequence. Letting $n$ go to infinity in
\[
(A_n - \lambda I_n|_{\im P_n})^{-1} P_n \, (P_n A P_n - \lambda I) P_n = P_n \quad \mbox{for} \; n \in \sM
\]
we obtain $B (A - \lambda I) = I$, thus, $A - \lambda I$ is invertible from the left-hand side. The invertibility from the right-hand side follows analogously.  \hfill \qed \\[3mm]
{\bf Proof of Theorem \ref{tf3.5}.} Let $A - \lambda I$ be invertible for some $\lambda \in \sR$. Since $A - \lambda I$ has a connected spectrum, this operator is either positively or negatively definite. But any kind of definiteness implies the stability of the sequence $(P_n(A - \lambda I)P_n)$. Conversely, if the sequence $(P_n(A - \lambda I)P_n)$ is stable, then the operator $A - \lambda I$ is invertible. Thus,
\[
\sigma(A) = \sigma_{\cF^\cP/\cG^\cP} ((P_n A P_n) + \cG^\cP).
\]
Further we know from Proposition \ref{t16.4} that
\[
\sigma_{\cF^\cP/\cG^\cP} ((P_n A P_n) + \cG^\cP) = \limsup \sigma(P_n A P_n),
\]
whereas we infer from Theorem \ref{tf3.6} that
\[
\sigma(A)  \subseteq  \liminf \sigma(P_nAP_n).
\]
These inclusions yield $\limsup \sigma(P_nAP_n) = \liminf \sigma(P_nAP_n)$. Hence, $(P_nAP_n)$ is a fractal sequence by Theorem \ref{tf3.2}. \hfill \qed
\subsection{Minimal stability spectra} \label{ssf3.4}
The stability spectrum $\sigma_{stab} (\bA) = \sigma_{\cF/\cG} (\bA + \cG)$ of a sequence $\bA \in \cF$ cannot increase when passing from $\bA$ to one of its subsequences $\bA_\eta := R_\eta(\bA)$ for a strictly increasing sequence $\eta$,
\[
\sigma_{\cF_\eta/\cG_\eta} (\bA_\eta + \cG_\eta) \subseteq \sigma_{\cF/\cG} (\bA + \cG).
\]
It is natural to ask how far the stability spectrum can decrease by passing to subsequences. We will say that a  sequence $\bA \in \cF$ has {\em minimal stability spectrum} \index{stability spectrum!minimal} if no subsequence of $\bA$ has a stability spectrum which is strictly less than that of $\bA$. The following is an immediate consequence of Proposition \ref{pf2.1} and Corollary \ref{cf3.21}.
\begin{prop} \label{pf3.31}
$(a)$ Every fractal sequence has a minimal stability spectrum. \\[1mm]
$(b)$ Every sequence possesses a subsequence with minimal stability spectrum.
\end{prop}
There are several questions coming up naturally for a sequence $\bA \in \cF$.
\begin{itemize}
\item
What is $\sigma_{inf} (\bA) := \cap_\eta \sigma_{stab} (R_\eta (\bA))$, \index{$\sigma_{inf} (\bA)$} the intersection taken over all strictly increasing sequences $\eta : \sN \to \sN$?
\item
Is there a sequence $\eta^*$ with $\sigma_{stab} (R_{\eta^*} (\bA)) = \sigma_{inf} (\bA)$, i.e., a sequence for which this intersection is attained?
\end{itemize}
Here is a partial answer to the first question.
\begin{prop} \label{pf3.32}
For every sequence $\bA = (A_n) \in \cF$,
\begin{equation} \label{ef3.33}
\liminf \sigma(A_n) \subseteq \sigma_{inf} (\bA).
\end{equation}
If $\bA$ is a normal sequence, then equality holds in $(\ref{ef3.33})$.
\end{prop}
{\bf Proof.} In the first part of the proof of Proposition \ref{t16.4} we have seen that
\begin{equation} \label{ef3.34}
\limsup \sigma(A_n) \subseteq \sigma_{stab} (\bA)
\end{equation}
for every sequence $\bA \in \cF$. Applying this inclusion to each subsequence of $\bA$ gives
\[
\cap_\eta \limsup \sigma(A_{\eta(n)}) \subseteq \cap_\eta \sigma_{stab} (R_\eta(\bA)) = \sigma_{inf} (\bA).
\]
The intersection on the left-hand side coincides with $\liminf \sigma(A_n)$, whence (\ref{ef3.33}). If the sequence $\bA$ is normal, then equality holds in (\ref{ef3.34}) by Proposition \ref{t16.4}. \hfill \qed \\[3mm]
We will show next that the answer to the second question is negative even if we restrict our attention to sequences of finite sections of self-adjoint operators. To give an appropriate example, we need some more facts about self-adjoint operators and their finite sections.

Let $H$ be a separable Hilbert space with orthonormal basis $(e_i)_{i \in \sN}$, and let $P_n$ be the orthogonal projection from $H$ onto the linear span of the first $n$ elements of the basis. By the Riesz-Fischer theorem, we can assume without loss that $H = l^2(\sN)$, provided with its standard basis. Every operator $P_n A P_n$ on $\im P_n$ will be identified via its matrix representation with respect to the standard basis with an $n \times n$-matrix. We denote the eigenvalues of a self-adjoint $n \times n$-matrix $B$ by
\[
\sigma_1(B) \le \sigma_2(B) \le \ldots \le \sigma_n(B)
\]
and call $\zs (B) := (\sigma_1(B), \, \sigma_2(B), \, \ldots, \, \sigma_n(B))$ the {\em ordered tuple of eigenvalues of} $B$. Further, a sequence $(\za^{(n)})_{n \in \sN}$ of ordered $n$-tuples $\za^{(n)} = (\alpha_1^{(n)}, \ldots, \, \alpha_k^{(n)})$ of real numbers is called {\em interlacing} \index{sequence!interlacing} if
\[
\alpha_1^{(n+1)} \le \alpha_1^{(n)} \le \alpha_2^{(n+1)} \le \alpha_2^{(n)} \le \ldots \le \alpha_n^{(n+1)} \le \alpha_n^{(n)} \le \alpha_{n+1}^{(n+1)}
\]
for every $n \in \sN$. The following is known as {\em Cauchy`s interlacing theorem}. \index{interlacing theorem} A proof is in \cite{Bha1}, Corollary III.1.5.
\begin{theo} \label{tf3.34a}
If $A \in L(H)$ is self-adjoint, then the sequence $(\zs (P_n A P_n))_{n \in \sN}$ of ordered tuples of eigenvalues is interlacing.
\end{theo}
Moreover, the sequence $(\zs (P_n A P_n))$ is bounded and
\[
|\sigma_k(P_n A P_n)| \le \|P_n A P_n\| \le \|A\| \quad \mbox{for every} \; n \in \sN.
\]
It turns out that the interlacing property and the boundedness are the only obstructions for the eigenvalues of the finite sections of a self-adjoint operator.
\begin{theo} \label{tf3.35}
Let $(\za^{(n)})_{n \in \sN}$ be a bounded and interlacing sequence of ordered $k$-tuples of real numbers. Then there is a self-adjoint operator $A \in L(l^2(\sN))$ with
\[
\zs (P_n A P_n) = \za^{(n)} \quad \mbox{for every} \; n \in \sN.
\]
\end{theo}
The proof of Theorem \ref{tf3.35} will follow by repeated application of the following {\em converse interlacing theorem}, which is Theorem III.1.9 in \cite{Bha1}. \index{interlacing theorem!converse}
\begin{theo} \label{tf3.36}
Let $\za = (\alpha_1, \, \ldots, \, \alpha_n)$ and $\zb = (\beta_1, \, \ldots, \, \beta_{n+1})$ be ordered tuples of real numbers with
\[
\beta_1 \le \alpha_1 \le \beta_2 \le \alpha_2 \le \ldots \le \beta_n \le \alpha_n \le \beta_{n+1},
\]
and let $A \in L(\im P_n)$ be a self-adjoint operator with $\zs(A) = \za$. Then there is a self-adjoint operator  $B \in L(\im P_{n+1})$ such that
\[
P_n B P_n|_{\im P_n} = A \quad \mbox{and} \quad \zs(B) = \zb.
\]
\end{theo}
{\bf Proof.} We identify $\im P_{n+1}$ with the orthogonal sum $\im P_n \oplus \sC e_{n+1}$. Accordingly, each operator $B \in L(\im P_{n+1})$ can be written as a block matrix
\[
C = \pmatrix{C_{11} & C_{12} \cr C_{21} & C_{22}} \quad \mbox{with} \; C_{11} \in L(\im P_n) \; \mbox{and} \; C_{22} \in \sC.
\]
Set $\widetilde{A} := \diag (\alpha_1, \, \ldots, \, \alpha_n)$, and let $U \in L(\im P_n)$ be a unitary operator such that $U^* A U = \widetilde{A}$. We will show that there is a matrix $C \in L(\im P_{n+1})$ of the form
\begin{equation} \label{ef3.37}
C = \pmatrix{\widetilde{A} & Z^* \cr Z & z_{n+1}} \quad \mbox{with} \; Z = (z_1, \, \ldots, \, z_n) \; \mbox{and} \; z_{n+1} \in \sR
\end{equation}
with $\zs(C) = \zb$. Then the matrix
\[
B := \pmatrix{U & 0 \cr 0 & 1} C  \pmatrix{U^* & 0 \cr 0 & 1}
\]
has the desired properties, since the unitary transformation does not effect the eigenvalues. So we are left with showing that there is a row matrix $Z$ and a real number $z_{n+1}$ such that the matrix $C$ in (\ref{ef3.37}) has $\zb$ as its ordered tuple of eigenvalues. Let
\[
P(\lambda) := \det \pmatrix{
\alpha_1 - \lambda & 0                  & \ldots & 0                  & \overline{z_1} \cr
0                  & \alpha_2 - \lambda & \ldots & 0                  & \overline{z_2} \cr
\vdots             & \vdots             & \ddots & \vdots             & \vdots         \cr
0                  & 0                  & \ldots & \alpha_n - \lambda & \overline{z_n} \cr
z_1                & z_2                & \ldots & z_n                & z_{n+1}}
\]
be the characteristic polynomial of $C$. Evaluating the determinant with respect to the last row we get
\[
P(\lambda) = (z_{n+1} - \lambda) \prod_{i=1}^n (\alpha_i - \lambda) - \sum_{i=1}^n |z_i|^2 \prod_{j=1, \, j\neq i}^n (\alpha_i - \lambda).
\]
On the other hand, in order to get $\zs(C) = \zb$ we must have $P(\lambda) = \prod_{j=1}^n (\beta_j - \lambda)$. Thus, we have to determine $z_1, \, \ldots, \, z_n \in \sC$ and $z_{n+1} \in \sR$ such that
\begin{equation} \label{ef3.37a}
(z_{n+1} - \lambda) \prod_{i=1}^n (\alpha_i - \lambda) - \sum_{i=1}^n |z_i|^2 \prod_{j=1, \, j\neq i}^n (\alpha_i - \lambda) = \prod_{j=1}^n (\beta_j - \lambda)
\end{equation}
for all $\lambda \in \sC$. We will take into account the multiplicities of the $\alpha_j$. Write
\[
(\alpha_1 - \lambda) \ldots (\alpha_n - \lambda) = (\gamma_1 - \lambda)^{k_1} \ldots (\gamma_r - \lambda)^{k_r}
\]
with pairwise different numbers $\gamma_i$ and positive integers $k_i$ such that $k_1 + \ldots k_r = n$. If $k_i > 1$ for some $i$, then the interlacing property implies that at least $k_i - 1$ of the $\beta_i$ are equal to $\gamma_i$. Thus, at least
\[
(k_1-1) + \ldots + (k_r-1) = k_1 + \ldots k_r - r = n-r
\]
of the $\beta_i$ appear among the $\gamma_j$. Denote the remaining $\beta_i$ by $\delta_1, \, \ldots, \, \delta_{r+1}$. Inserting the new notation into (\ref{ef3.37a}) and canceling the common factor $\prod_{j=1}^r (\gamma_j - \lambda)^{k_j - 1}$ we get
\begin{equation} \label{ef3.37b}
(z_{n+1} - \lambda) \prod_{i=1}^r (\gamma_i - \lambda) - \sum_{i=1}^n |z_i|^2 \frac{\prod_{j=1}^r (\gamma_i - \lambda)}{\alpha_i - \lambda} = \prod_{i=1}^{r+1} (\delta_i - \lambda).
\end{equation}
Evaluating this equation at $\lambda = \gamma_{i_0}$ with $1 \le i_0 \le r$ we find
\[
- \sum_{i: \alpha_i = \gamma_{i_0}} |z_i|^2 \prod_{j=1, j \neq i_0}^r (\gamma_i - \gamma_{i_0}) = \prod_{i=1}^{r+1} (\delta_i - \gamma_{i_0}).
\]
To get a unique solution, we seek only non-negative $z_i$ and assume moreover that $z_i = z_j$ if $\alpha_i = \alpha_j$. Let $z_j =: x_i$ if $\alpha_j = \gamma_i$. Then $x_{i_0}$ satisfies
\[
x_{i_0}^2 \cdot k_{i_0}  \prod_{j=1, j \neq i_0}^r (\gamma_i - \gamma_{i_0}) = \prod_{i=1}^{r+1} (\delta_i - \gamma_{i_0})
\]
and is uniquely determined by this equation. It remains to compute $z_{n+1}$. For this goal, we insert the already determined numbers $z_1, \, \ldots, \, z_n$ into (\ref{ef3.37b}) and evaluate this equation at an arbitrarily chosen real point $\lambda^\prime \neq \gamma_i$ for $1 \le i \le r$. This equation has a unique real solution $z_{n+1}$.
\hfill \qed \\[3mm]
{\bf Proof of Theorem \ref{tf3.35}.} We use the converse interlacing theorem to construct a sequence $(A_n)$ of $n \times n$-matrices $A_n$ as follows. Let $\za^{(1)} = (\alpha_1^{(1)})$ and set $A_1 := (\alpha_1^{(1)})$. Then $\zs(A_1) = \za^{(1)}$. Assume we have already found an $n \times n$-matrix $A_n$ with $\zs(A_n) = \za^{(n)}$. Then we construct $A_{n+1}$ by means of Theorem \ref{tf3.36} such that $\zs(A_{n+1}) = \za^{(n+1)}$.
We identify each matrix $A_n$ with an operator on the subspace $\im P_n$ of $H$ and denote this operator by $A_n$ again. The operators $A_n \in L(H)$ are uniformly bounded (the norm of the self-adjoint operator $A_n$ is equal to its spectral radius, and the spectral radii are uniformly bounded by hypothesis), and they converge strongly on the linear span of the basis $\{e_i\}$, which is a dense subset of $H$. By the Banach-Steinhaus theorem, the operators $A_n$ converge strongly to a bounded linear operator $A$ on $H$. This operator is self-adjoint, and $P_n A P_n = A_n$ by construction. \hfill \qed \\[3mm]
The following example shows that the intersection of the stability spectra of subsequences of a given sequence is not necessarily the stability spectrum of a subsequence again, thus answering the second of the above questions in the negative. Moreover, this example shows that the stability spectra of the subsequences of a given sequences are not necessarily linearly ordered.
\begin{example} \label{exf3.40}
Define a sequence $(\za^{(n)})_{n \in \sN}$ of ordered $k$-tuples as follows. For $n = 2k$ with $k \ge 2$, set
\[
\alpha_1^{(n)} = \ldots = \alpha_k^{(n)} = 0, \quad \alpha_{k+1}^{(n)} = 1, \quad \alpha_{k+2}^{(n)} = \ldots =
\alpha_n^{(n)} = 3,
\]
and for $n = 2k+1$ set
\[
\alpha_1^{(n)} = \ldots = \alpha_{k+1}^{(n)} = 0, \quad \alpha_{k+2}^{(n)} = 2, \quad \alpha_{k+3}^{(n)} = \ldots =
\alpha_n^{(n)} = 3.
\]
For $n = 1, \, 2, \, 3$ we choose $\alpha_i^{(n)} \in \{0, \, 1, \, 2, \, 3\}$ such that the interlacing property is
satisfied. Then $(\za^{(n)})_{n \in \sN}$ is a bounded sequence of ordered $k$-tuples with interlacing property. By  Theorem \ref{tf3.35}, there is a bounded self-adjoint operator $A$ on $l^2(\sN)$ such that $\zs(P_nAP_n) = \za^{(n)}$ for every $n \in \sN$. Thus, for $k \ge 2$,
\[
\sigma(P_{2k}AP_{2k}) = \{0, \, 1, \, 3\} \quad \mbox{and} \quad \sigma(P_{2k+1}AP_{2k+1}) = \{0, \, 2, \, 3\}.
\]
The sequences $(P_{2k}AP_{2k})$ and $(P_{2k+1}AP_{2k+1})$ are fractal by Theorem \ref{tf3.2} and have, thus, minimal stability spectra. The self-adjointness of these sequences further implies that
\[
\sigma_{stab} ((P_{2k}AP_{2k})_{k \in \sN}) = \{0, \, 1, \, 3\} \quad \mbox{and} \quad \sigma_{stab} ((P_{2k+1}AP_{2k+1})_{k \in \sN}) = \{0, \, 2, \, 3\}.
\]
Hence, $\sigma_{inf} ((P_nAP_n)_{n \in \sN}) = \{0, \, 3\}$, but there is no subsequence of $(P_nAP_n)$ which has $\{0, \, 3\}$ as its stability spectrum.  \hfill \qed
\end{example}
\section{Essential Fractality} \label{s3}
Recall that a $C^*$-subalgebra $\cA$ of $\cF$ is fractal if each sequence $(A_n) \in \cA$ can be rediscovered from each of its (infinite) subsequences modulo a sequence in the ideal $\cG$. There are plenty of subalgebras of $\cF$ which arise from concrete discretization methods and which are fractal (we have seen the finite sections algebra $\cS(\eT(C))$ for Toeplitz operators as one example). On the other hand, the algebra of the finite sections method for general band-dominated operators is an example of an algebra which fails to be fractal, as the finite sections of the operator
\[
A = \diag \left( \pmatrix{0 & 1 \cr 1 & 0}, \, \pmatrix{0 & 1 \cr 1 & 0}, \, \pmatrix{0 & 1 \cr 1 & 0}, \, \ldots \right)
\]
show. But this algebra enjoys a weaker form of fractality which we call essential fractality and which we are going to introduce now. First we have to introduce an distinguished ideal of $\cF$ which plays a similar role as the ideal $K(H)$ of the compact operators on an infinite-dimensional Hilbert space $H$.

To motivate the definition of this ideal, we come back to our running example, the algebra $\cS(\eT(C))$ of the finite sections method for Toeplitz operators. There is an ideal hidden in the algebra $\cS(\eT(C))$ which did not appear explicitly in the previous considerations but which always
acted as a player in the background and which will play (together
with its relatives) an outstanding role in what follows. Let
\begin{equation} \label{e14.35a}
\cJ = \{ (P_n K P_n + R_n L R_n + G_n) : K, \, L \; \mbox{compact}, \; (G_n) \in \cG \}.
\end{equation}
\begin{theo} \label{t14.36}
$(a)$ $\cJ$ is a closed ideal of $\cS(\eT(C))$. \\[1mm]
$(b)$ The quotient algebra $\cS(\eT(C))/\cJ$ is $^*$-isomorphic to $C(\sT)$, and the mapping $(P_n T(a) P_n) + \cJ \mapsto a$ is a $^*$-isomorphism between these algebras.
\end{theo}
{\bf Proof.} $(a)$ First note that $\cJ \subset \cS(\eT(C))$ by Theorem \ref{t14.21}. The closedness of $\cJ$ in $\cS(\eT(C))$ follows by standard arguments. To check that $\cJ$ is a left ideal, let $a \in C(\sT)$ and let $K$ and $L$ be compact. Then
\begin{eqnarray*}
\lefteqn{P_n T(a) P_n (P_n K P_n + R_n L R_n)} \\
&& = P_n T(a) P_n K P_n + R_n (R_n T(a) R_n) L R_n \\
&& = P_n T(a) P_n K P_n + R_n T(\tilde{a}) P_n L R_n \\
&& = P_n T(a)K P_n + R_n T(\tilde{a}) L R_n - P_n T(a) Q_n K P_n - R_n T(\tilde{a}) Q_n L R_n
\end{eqnarray*}
with $Q_n := I - P_n$. The operators $T(a)K$ and $T(\tilde{a}) L$ are compact. Since the operators $Q_n$ converge strongly to zero and $K$ and $L$ are compact, the last two operators converge to zero in the norm. Hence, $(P_n T(a) P_n) \, (P_n K P_n + R_n L R_n) \in \cJ$. Similarly one checks that $\cJ$ is a right ideal. \\[2mm]
$(b)$ Widom's identity (\ref{e170309.4}) together with the compactness of Hankel operators with continuous generating function imply that the mapping $a \mapsto (P_n T(a) P_n) + \cJ$ is a $^*$-homomorphism from $C(\sT)$ into $\cS(\eT(C))/\cJ$. This homomorphism is surjective by Theorem \ref{t14.21}. To get its injectivity, let $(P_n T(a) P_n) \in \cJ$ for a continuous function $a$. Then there are compact operators $K$ and $L$ and a zero sequence $(G_n)$ such that
\[
P_n T(a) P_n = P_n K P_n + R_n L R_n + G_n \quad \mbox{for all}
\, n \in \sN.
\]
Letting $n$ go to infinity yields the compactness of $T(a)$. But then $a$ is the zero function. \hfill \qed \\[3mm]
The importance of the ideal $\cJ$ results from several facts:
\begin{itemize}
\item
The algebra $\cS(\eT(C))/\cJ$ is commutative, hence subject to Gelfand-Naimark theory of commutative $C^*$-algebras. Similarly, factorization of a subalgebra $\cA$ of $\cF$ by $\cJ$ (or by an ideal with similar properties; see below) often yields quotient algebras $\cA/\cJ$ which can be effectively studied by tools like central localization or other non-commutative generalizations of Gelfand theory.
\item
The algebra $\cJ/\cG$ has exactly two non-equivalent irreducible representations which are given by the homomorphisms $W$ and $\widetilde{W}$. These representations extend to representations of $\cS(\eT(C))$ (of course, also given by $W$ and $\widetilde{W}$) with the property that a sequence $\bA$ in $\cS(\eT(C))$ is stable if and only if the operators $W(\bA)$ and $\widetilde{W}(\bA)$ are invertible. In this sense, the irreducible representations of $\cJ$ yield a sufficient family of irreducible representations of $\cS(\eT(C))$. Similar effects can be observed in numerous instances.
\item
Invertibility modulo $\cJ$ can be lifted in the following sense. Let $\cF^\cJ$ stand for the largest subalgebra of $\cF$ for which $\cJ$ is an ideal. Then the mappings $W$ and $\widetilde{W}$ extend to irreducible representations of $\cF^\cJ$, and a sequence $\bA \in \cF^\cJ$ is stable if and only if the operators $W(\bA)$ and $\widetilde{W}(\bA)$ are invertible and if the coset $\bA + \cJ$ is invertible in the quotient $\cF^\cJ/\cJ$. Again, such a lifting result holds in a much more general context.
\end{itemize}
The ideal $\cJ$ is clearly related with compact operators. We will now introduce a larger ideal $\cK$ of sequences of compact type. Throughout this section, $\cF$ will be an algebra of matrix sequences with dimension function $\delta$.
\subsection{Compact sequences} \label{ss32}
A sequence $(K_n)$ in the $C^*$-algebra $\cF$ is a {\em sequence of rank one matrices} \index{sequence!of rank one matrices} if every matrix $K_n$ has range dimension less than or equal to one. The smallest closed ideal of $\cF$ which contains all sequences of rank one matrices will be denoted by $\cK$. Thus, a sequence $(A_n) \in \cF$ belongs to $\cK$ if and only if, for every $\varepsilon > 0$, there is a sequence $(K_n) \in \cF$ such that
\begin{equation} \label{e41.1}%
\sup_n \|A_n - K_n\| < \varepsilon \quad \mbox{and} \quad \sup_n \, \rank K_n < \infty.
\end{equation}
We refer to the elements of $\cK$ as {\em compact sequences}. The role of the ideal $\cK$ in numerical analysis can be compared with the role of the ideal of the compact operators in operator theory.

Note that $\cG \subseteq \cK$. Indeed, given a sequence $(G_n) \in \cG$ and an $\varepsilon > 0$, set $K_n := G_n$ if $\|G_n\| \ge \varepsilon$ and $K_n := 0$ otherwise. Then $(G_n)$ satisfies (\ref{e41.1}) in place of $(A_n)$ since there are only finitely many operators $K_n$ which are not zero.

An appropriate notion of the rank of a sequence in $\cF$ can be introduced as
follows. We say that a sequence $\bA \in \cF$ has {\em finite essential rank}
\index{sequence!of finite essential rank} if it is the sum of a sequence $(G_n)$
in $\cG$ and of a sequence $(K_n)$ with $\sup_n \rank K_n < \infty$. If $\bA$ is
of finite essential rank, then there is a smallest integer $r \ge 0$ such that
$\bA$ can be written as $(G_n) + (K_n)$ with $(G_n) \in \cG$ and $\sup_n \rank
K_n = r$. We call this integer the {\em essential rank}
\index{rank!essential} of $\bA$ and write $\ess \rank \bA = r$. If $\bA$ is not
of finite essential rank, then we put $\ess \rank \bA = \infty$. Thus, the
sequences of essential rank 0 are just the sequences in $\cG$. Clearly, the
sequences of finite essential rank form an ideal of $\cF$ which is dense in
$\cK$, and
\[
\ess \rank (\bA + \bB) \le \ess \rank \bA + \ess \rank \bB,
\]
\[
\ess \rank (\bA \bB) \le \min \, \{ \ess \rank \bA, \, \ess \rank
\bB \}
\]
for arbitrary sequences $\bA, \bB \in \cF$.

Consider our running example. It is not hard to see that the intersection of the algebra $\cS(\eT(C))$ with the ideal $\cK$ is just the distinguished ideal $\cJ$ which we examined in the beginning of the preceding section. Moreover, the essential rank of the sequence $(P_nKP_n + R_nLR_n + G_n)$ turns out to be $\rank K + \rank L$.

There are several equivalent characterizations of compact sequences. Since the entries $A_n$ of the sequences are $n \times n$-matrices, a characterization of the compactness property and of the essential rank via the singular values of the $A_n$ will be particularly useful for our purposes. Recall from linear algebra that the singular values of an $n \times n$ matrix $A$ are the non-negative square roots of the eigenvalues of $A^*A$. We denote them by \index{$\Sigma_k(A)$}
\begin{equation} \label{e41.5}%
\|A\| = \Sigma_1(A) \ge \Sigma_2(A) \ge \ldots \ge  \Sigma_n (A) \ge 0
\end{equation}
if they are ordered decreasingly and by \index{$\sigma_k(A)$}
\begin{equation} \label{e41.6}%
0 \le \sigma_1(A) \le \sigma_2(A) \le \ldots \le \sigma_n(A) = \|A\|
\end{equation}
in case of increasing order. Thus, $\sigma_k(A) = \Sigma_{n-k+1}(A)$. Since the matrices $A^*A$ and $AA^*$ are unitarily equivalent, $\Sigma_k(A) = \Sigma_k (A^*)$ for every $k$. We will also need the fact that every $n \times n$ matrix $A$ has a {\em singular value decomposition}
\[
A = E^* \, \diag (\Sigma_1 (A), \, \ldots, \, \Sigma_n (A)) F
\]
with unitary matrices $E$ and $F$. Another simple fact from linear algebra which we will use several times is the following.
\begin{lemma} \label{l41.9}
Let $A$ be an $n \times n$ matrix with rank $r$. Then $\rank A^\prime \ge r$ for each $n \times n$ matrix $A^\prime$ with $\|A - A^\prime\| < \Sigma_r (A)$.
\end{lemma}
The announced characterization of compact sequences in terms of singular values reads as follows.
\begin{theo} \label{t41.12}
The following conditions are equivalent for a sequence $(K_n) \in \cF$: \\[1mm]
$(a)$ $\lim_{k \to \infty} \sup_{n \ge k} \Sigma_k (K_n) = 0$; \\[1mm]
$(b)$ $\lim_{k \to \infty} \limsup_{n \to \infty} \Sigma_k (K_n) = 0$; \\[1mm]
$(c)$ the sequence $(K_n)$ is compact.
\end{theo}
Since the sequence $k \mapsto \sup_{n \ge k} \Sigma_k (K_n)$ is decreasing, the limit in $(a)$ and $(b)$ can be replaced by an infimum. \\[3mm]
{\bf Proof.} The implication $(a) \Rightarrow (b)$ is evident. Let $(K_n) \in \cF$ be a sequence which satisfies condition $(b)$, and let
\[
K_n = E_n^* \, \diag (\Sigma_1 (K_n), \, \ldots, \, \Sigma_{\delta(n)} (K_n)) F_n
\]
be the singular value decomposition of $K_n$. For every $n \in \sN$ and $k \ge 1$, set
\[
K^{(k)}_n := \left\{
\begin{array}{lll}
E_n^* \, \diag (\Sigma_1(K_n) \, \ldots, \, \Sigma_{k-1} (K_n), 0, \ldots, \, 0) F_n & \mbox{if} & 1 < k \le n, \\
0   & \mbox{if} & 1 = k \le n, \\
K_n & \mbox{if} & n < k.
\end{array} \right.
\]
Then, for $n > k$,
\[
\| K_n - K^{(k)}_n \| = \| E_n^* \, \diag (0, \, \ldots, \, 0, \, \Sigma_k(K_n), \, \ldots, \, \Sigma_{\delta(n)} (K_n)) F_n \| = \Sigma_k (K_n),
\]
and the limsup formula (\ref{e170309.3}) for the norm of a coset in $\cF/\cG$ yields
\[
\|(K_n) - (K^{(k)}_n) + \cG \|_{\cF/\cG} = \limsup_{n \to \infty} \Sigma_k (K_n).
\]
Together with property $(b)$, this implies that
\[
\lim_{k \to \infty} \| (K_n) - (K^{(k)}_n) + \cG\|_{\cF/\cG} = \lim_{k \to \infty} \limsup_{n \to \infty} \Sigma_k (K_n) = 0.
\]
Thus, for each $k \in \sN$, there is a sequence $(C_n^{(k)})$ in $\cG$ such that
\[
\lim_{k \to \infty} \| (K_n) - (K^{(k)}_n) - (C_n^{(k)})\|_\cF = 0,
\]
i.e., the sequence $(K_n)$ is the limit as $k \to \infty$ of the sequences $(K^{(k)}_n + C_n^{(k)})_{n \in \sN}$. Since $\rank K^{(k)}_n \le k-1$ by definition, each of these sequences belongs to $\cK$. Hence, $(K_n)$ is a compact sequence.

For the implication $(c) \Rightarrow (a)$, take a compact sequence $(K_n)$. The sequence $(\sup_{n \ge k} \Sigma_k (K_n))_{k \ge 1}$ is monotonically decreasing and bounded below (by zero), hence, convergent. Assume that the limit of this sequence is positive. Then there is a $C > 0$ such that $\sup_{n \ge k} \, \Sigma_k (K_n) > C$ for all $k \ge 1$. Thus, there are numbers $n_k \ge k$ such that
\begin{equation} \label{e41.13}
\Sigma_k (K_{n_k}) > C \quad \mbox{for all} \; k \ge 1.
\end{equation}
On the other hand, since the sequence $(K_n)$ is compact, there is a sequence $(R_n) \in \cF$ with
\begin{equation} \label{e41.14}
\sup_n \, \rank R_n < \infty \quad \mbox{and} \quad \sup_n \|K_n - R_n \| < C.
\end{equation}
In particular, for each $k$ one has $\|K_{n_k} - R_{n_k}\| < C$, which implies via Lemma \ref{l41.9} and (\ref{e41.13}) that $\rank R_{n_k} \ge k$. Since $k$ can be chosen arbitrarily large, this contradicts the first condition in
(\ref{e41.14}). Hence, the sequence $(\sup_{n \ge k} \Sigma_k (K_n))_{k \ge 1}$ cannot have a positive limit, whence condition $(a)$. \hfill \qed \\[3mm]
In the same vein one can prove the following characterization of sequences of essential rank $r$.
\begin{coro} \label{c41.15}
A sequence $(K_n) \in \cF$ is of essential rank $r$ if and only if
\[
\limsup_{n \to \infty} \Sigma_r (K_n) > 0 \quad \mbox{and} \quad \lim_{n \to \infty} \Sigma_{r+1} (K_n) = 0.
\]
\end{coro}
One consequence is the {\em lower semi-continuity} of the essential rank function.
\begin{coro} \label{c41.16}
If $\ess \rank (K_n) = r$, then $\ess \rank (K_n^\prime) \ge r$ for all sequences $(K_n^\prime)$ which are sufficiently close to $(K_n)$.
\end{coro}
Another corollary concerns the behavior of the small singular values of $K_n$.
\begin{coro} \label{c41.17}
Let $(K_n) \in \cK$. Then the limit $\lim_{n \to \infty} \sigma_k (K_n)$ exists and is equal to $0$ for every $k$.
\end{coro}
{\bf Proof.} Let $\varepsilon > 0$. By Theorem \ref{t41.12}, there is a $k_0$ such that $\sup_{n \ge k_0} \Sigma_{k_0} (K_n) < \varepsilon$. Then, for all $n \ge n_0 := k_0 + k -1$,
\[
\sigma_k (K_n) = \Sigma_{\delta(n)-k+1} (K_n) \le \Sigma_{k_0} (K_n) \le \sup_{n \ge k_0} \Sigma_{k_0} (K_n) < \varepsilon,
\]
which gives the assertion. \hfill \qed \\[3mm]
We still mention two other characterizations of the ideal of the compact sequences in $\cF$.
\begin{theo} \label{t41.27}
$\cK$ is the smallest closed ideal of $\cF$ which contains the constant sequence $(P_1)$.
\end{theo}
\begin{theo} \label{t41.30}
$(a)$ A sequence $\bK \in \cF$ belongs to the ideal $\cK$ if and only if $W(\bK)$ is compact for every irreducible representation $W$ of $\cF$. \\[1mm]
$(b)$ A coset $\bK + \cG \in \cF/\cG$ belongs to the ideal $\cK/\cG$ if and only if $W(\bK + \cG)$ is compact for every irreducible representation $W$ of $\cF/\cG$.
\end{theo}
A crucial step in the proof is to show that $\rank W(\bK) \le 1$ for each sequence $\bK \in \cF$ of rank one matrices and each irreducible representation $W$ of $\cF$. For proofs of the preceding theorems and further characterizations of the ideal $\cK$, see \cite{Roc10}.
\subsection{Fredholm sequences} \label{ss33}
Corresponding to the ideal $\cK$ we introduce an appropriate class of Fredholm sequences by calling a sequence $(A_n) \in \cF$ {\em Fredholm} \index{sequence!Fredholm} if it is invertible modulo the ideal $\cK$ of the compact sequences. The following properties of Fredholm sequences are obvious. \\[1mm]
-- Stable sequences are Fredholm.\\
-- Adjoints of Fredholm sequences are Fredholm.\\
-- Products of Fredholm sequences are Fredholm.\\
-- The sum of a Fredholm and a compact sequence is Fredholm.\\
-- The set of all Fredholm sequences is open in $\cF$. \\[1mm]
For alternate characterizations of Fredholm sequences, let $\sigma_1 (A) \le \ldots \le \sigma_n (A)$ denote the singular values of an $n \times n$ matrix $A$.
\begin{theo} \label{t42.1}
The following conditions are equivalent for a sequence $(A_n)
\in \cF$: \\[1mm]
$(a)$ The sequence $(A_n)$ is Fredholm. \\[1mm]
$(b)$ There are sequences $(B_n) \in \cF$ and $(J_n) \in \cK$ with $\sup_n \mbox{\rm rank} \, J_n < \infty$ such that
\begin{equation} \label{e42.2}%
B_n A_n = I_n + J_n \quad \mbox{for all} \; n \in \sN.
\end{equation}
$(c)$ There is a $k \in \sN$ such that
\begin{equation} \label{e42.3}%
\liminf_{n \to \infty} \sigma_{k+1} (A_n) > 0.
\end{equation}
\end{theo}
{\bf Proof.} $(a) \Rightarrow (b)$: Let $(A_n) \in \cF$ be a Fredholm sequence. Then there are sequences $(C_n) \in \cF$ and $(K_n) \in \cK$ such that $(C_n) (A_n) = (I_n) + (K_n)$. Choose a sequence $(L_n) \in \cK$ with $\|(L_n) - (K_n)\|_\cF < 1/2$ and $\sup \rank L_n < \infty$. Then
\[
(C_n) (A_n) = (I_n) + (K_n - L_n) + (L_n).
\]
Since $(I_n) + (K_n - L_n)$ is invertible in $\cF$, we obtain (\ref{e42.2}) with
\[
B_n := (I_n + K_n - L_n)^{-1} C_n \quad \mbox{and} \quad J_n := (I_n + K_n - L_n)^{-1} L_n.
\]
$(b) \Rightarrow (c)$: Let the singular value decomposition of $A_n$ be given by
\[
A_n = E_n^* \Sigma_n F_n := E_n^* \, \diag (\sigma_1 (A_n), \, \ldots, \, \sigma_{\delta(n)} (A_n)) F_n.
\]
After multiplication by $F_n$ and $F_n^*$, the identity (\ref{e42.2}) becomes
\[
(F_n B_n E_n^*) (\Sigma_n) = (I_n) + (F_n J_n F_n^*).
\]
Abbreviating $C_n := F_n B_n E_n^*$ and $K_n := F_n J_n F_n^*$ we get
\begin{equation} \label{e42.5}
C_n \Sigma_n = C_n \diag (\sigma_1 (A_n), \, \ldots, \, \sigma_{\delta(n)} (A_n)) = I_n + K_n \quad \mbox{for all} \; n \in \sN
\end{equation}
where still $\sup_n \rank K_n < \infty$. Let $k := \limsup_{n \to \infty} \rank K_n$. We claim that $\liminf_{n \to \infty} \sigma_{k+1} (A_n) > 0$. Contrary to what we want to show, assume that there is an infinite subsequence $(n_l)_{l \ge 1}$ of $\sN$ with $\lim_{l \to \infty} \sigma_{k+1} (A_{n_l}) = 0$. Multiplying (\ref{e42.5}) from both sides by $P_{k+1}$, we get
\[
P_{k+1} C_{n_l} \Sigma_{n_l} P_{k+1} = P_{k+1} + P_{k+1} K_{n_l} P_{k+1}.
\]
Since
\[
\| \Sigma_{n_l} P_{k+1} \| = \| \diag (\sigma_1 (A_{n_l}), \, \ldots, \, \sigma_{k+1} (A_{n_l}), \, 0, \, \ldots, \, 0) \| = \sigma_{k+1} (A_{n_l}) \to 0,
\]
one has
\[
\lim_{l \to \infty} \|P_{k+1} + P_{k+1} K_{n_l} P_{k+1}\| = 0.
\]
Thus, the matrices $P_{k+1} K_{n_l} P_{k+1} \in \sC^{(k+1) \times (k+1)}$ are invertible for all sufficiently large $n_l$. But this is impossible since $P_{k+1}$ has rank $k+1$, whereas $\rank K_{n_l} \le k$. This proves the claim which, on its hand, implies assertion $(c)$. \\[2mm]
$(c) \Rightarrow (a)$: As in the previous part of the proof, let $A_n = E_n^* \Sigma_n F_n$ refer to the singular value decomposition of $A_n$, and let $k$ be a non-negative integer such that
\[
\liminf_{n \to \infty} \sigma_{k+1} (A_n) > 0.
\]
Then the sequence $(\Sigma_n + P_k)_{n \ge 1}$ (with $P_0 := 0$) is stable, and so is the sequence $(A_n + E_n^* P_k F_n)_{n \in \sN}$. Thus, there are sequences $(C_n) \in \cF$ and $(G_n), \, (H_n) \in \cG$ such that
\[
(C_n)(A_n + E_n^* P_k F_n) = (I_n) + (G_n) \quad \mbox{and} \quad (A_n + E_n^* P_k F_n)(C_n) = (I_n) + (H_n),
\]
whence
\[
(C_n)(A_n) = (I_n) + (G_n) - (C_n E_n^* P_k F_n)
\]
and
\[
(A_n)(C_n) = (I_n) + (H_n) - (E_n^* P_k F_n C_n).
\]
The sequences $(G_n) - (C_n E_n^* P_k F_n)$ and $(H_n) - (E_n^* P_k F_n C_n)$ are of finite essential rank. Hence, $(A_n)$ is invertible modulo $\cK$. \hfill \qed \\[3mm]
The preceding theorem suggests to introduce the {\em $\alpha$-number} $\alpha(\bA)$ \index{$\alpha$-number} of a Fredholm sequence $\bA = (A_n)$, which corresponds to the kernel dimension of a Fredholm operator. By definition, $\alpha (\bA)$ is the smallest non-negative integer $k$ for which (\ref{e42.3}) is true. Equivalently, $\alpha (\bA)$ is the smallest non-negative integer $k$ for which there exist a sequence $(B_n) \in \cF$ and a sequence $(J_n) \in \cK$ of essential rank $k$ such that $B_n A_n^* A_n = I_n + J_n$ for all $n \in \sN$. The latter fact follows easily from the proof of the preceding theorem.

The {\em index} \index{index} of a Fredholm sequence $\bA$ is the integer
\[
\ind (\bA) := \alpha (\bA) - \alpha (\bA^*).
\]
It turns out that, in the case at hand, the index of a Fredholm sequence always zero. This is a consequence of the fact that the operators $A_n$ act on finite dimensional spaces which implies that $A_n^* A_n$ and $A_n A_n^*$ have the same eigenvalues, even with respect to their multiplicity. So the more interesting quantity associated with a Fredholm sequence seems to be its $\alpha$-number. On the other hand, the vanishing of the index of a Fredholm sequence allows one to make use of the index as a conservation quantity.

If $A$ is a Fredholm operator with index 0, then there is an operator $K$ with finite rank such that $A+K$ is invertible. The analog for Fredholm sequences reads as follows. Notice that there is no index obstruction since the index of a Fredholm sequence is always 0.
\begin{theo} \label{t42.10}
If $\bA \in \cF$ is a Fredholm sequence, then there is a sequence $\bK \in \cK$ with $\ess \rank \bK \le \alpha (\bA)$ such that $\bA + \bK$ is a stable sequence.
\end{theo}
{\bf Proof.} Let $k$ denote the $\alpha$-number of $\bA =: (A_n)$, and let
\[
A_n = E_n^* \, \diag (\sigma_1 (A_n), \, \ldots, \, \sigma_{\delta(n)} (A_n)) F_n
\]
be the singular value decomposition of $A_n$. Set
\[
K_n :=  E_n^* \, \diag (\sigma_{k+1} (A_n) - \sigma_1 (A_n), \, \ldots, \, \sigma_{k+1} (A_n) - \sigma_k (A_n), \, 0, \, \ldots, \, 0) F_n
\]
and $\bK := (K_n)$. Then $\rank K_n \le k$ for each $n$, hence $\ess \rank \bK \le k$, and the sequence $\bA + \bK$ is stable. \hfill \qed
\subsection{Fractality of quotient maps} \label{ss43.1}
We still let $\cF$ be the algebra of matrix sequences with dimension function $\delta$ and $\cG$ the associated ideal of zero sequences, but in the present and the subsequent subsections we could also allow for $\cF$ to be the product of a family $(\cC_n)_{n \in \sN}$ of unital $C^*$-algebras. Given a strictly increasing sequence $\eta :\sN \to \sN$ and a $C^*$-subalgebra $\cA$ of $\cF$, we define the restriction mapping $R_\eta : \cF \to \cF_\eta$ and the image $\cA_\eta$ of $\cA$ under this mapping as in Section \ref{ss22}.
\begin{theo} \label{t43.1}
Let $\cA$ be a $C^*$-subalgebra of $\cF$ and $\cJ$ a closed ideal of $\cA$. The canonical homomorphism $\pi^\cJ : \cA \to \cA/\cJ$ is fractal if and only if the following implication holds for every sequence $\bA \in \cA$ and every strictly increasing sequence $\eta : \sN \to \sN$
\begin{equation} \label{e43.2}
R_\eta (\bA) \in \cJ_\eta \quad \Longrightarrow \quad \bA \in \cJ.
\end{equation}
\end{theo}
{\bf Proof.} Let $\pi^\cJ$ be fractal, i.e., for each $\eta$, there is a mapping $\pi^\cJ_\eta$ such that $\pi^\cJ = \pi^\cJ_\eta R_\eta|_\cA$. Let $R_\eta (\bA) \in \cJ_\eta$ for a sequence $\bA \in \cA$. We choose a sequence $\bJ \in \cJ$ such that $R_\eta (\bA) = R_\eta (\bJ)$. Applying the homomorphism $\pi^\cJ_\eta$ to both sides of this equality we obtain $\pi^\cJ (\bA) = \pi^\cJ (\bJ) = 0$, whence $\bA \in \cJ$.

For the reverse implication, let $\bA$ and $\bB$ be sequences in $\cA$ with $R_\eta (\bA) = R_\eta (\bB)$. Then $R_\eta (\bA - \bB) = 0 \in \cJ_\eta$, and (\ref{e43.2}) implies that $\bA - \bB \in \cJ$. Thus, the mapping
\[
\pi^\cJ_\eta : \cA_\eta \to \cA/\cJ, \quad R_\eta (\bA) \mapsto \bA + \cJ
\]
is correctly defined, and it satisfies $\pi^\cJ_\eta R_\eta|_\cA = \pi^\cJ$.  \hfill \qed \\[3mm]
Let now $\cJ$ be a closed ideal of $\cF$. Then $\cA \cap \cJ$ is a closed ideal of $\cA$, and the preceding theorem states that the canonical mapping $\pi^{\cA \cap \cJ} : \cA \to \cA/(\cA \cap \cJ)$ is fractal if and only if the implication
\begin{equation} \label{e43.2a}
R_\eta (\bA) \in (\cA \cap \cJ)_\eta \quad \Longrightarrow \quad \bA \in \cJ
\end{equation}
holds for every sequence $\bA \in \cA$ and every strictly increasing sequence $\eta$. It would be much easier to check this implication if one would have
\begin{equation} \label{e43.2b}
(\cA \cap \cJ)_\eta = \cA_\eta \cap \cJ_\eta
\end{equation}
foe every $\eta$, in which case the implication (\ref{e43.2a}) reduces to $R_\eta (\bA) \in \cJ_\eta \Rightarrow \bA \in \cJ$. Recall from Theorem \ref{tf1.8} $(b)$ that (\ref{e43.2b}) indeed holds if $\cJ = \cG$ and if the canonical homomorphism $\pi^{\cA \cap \cG} : \cA \to \cA/(\cA \cap \cG)$ is fractal. The following example shows that one cannot expect an analogous result for arbitrary closed ideals $\cJ$ of $\cF$.
\begin{example} \label{ex43.6}
Let $\cA := \cS(\eT(C))$, the algebra of the finite sections method for Toeplitz operators, and $\cK$ the ideal of the compact sequences in $\cF$. Then
\[
\cJ := \{ (K_n) \in \cK : \lim_{n \to \infty} \|K_{2n}\| = 0 \}
\]
is a closed ideal of $\cF$. We claim that
\begin{equation} \label{e43.7}
\cS(\eT(C)) \cap \cJ = \cG.
\end{equation}
The inclusion $\cG \subseteq \cS(\eT(C)) \cap \cJ$ is evident. For the reverse inclusion, let $(A_n) \in \cS(\eT(C)) \cap \cJ$. By Theorem \ref{t14.21} and its proof, there is a unique representation
\[
A_n = P_n T(a) P_n + P_n K P_n + R_n L R_n + G_n
\]
with a continuous function $a$, compact operators $K$ and $L$ and a sequence $(G_n) \in \cG$, and one has
\[
\mbox{s-lim} \, A_n P_n = T(a) + K \quad \mbox{and} \quad \mbox{s-lim} \, R_n A_n R_n = T(\tilde{a}) + L
\]
with $\tilde{a}(t) = a(1/t)$. Since $\|A_{2n}\| \to 0$, this implies $T(a) + K = T(\tilde{a}) + L = 0$, whence $a = 0$ and $K = L = 0$. Thus, $(A_n) \in \cG$.

As a consequence of (\ref{e43.7}), the canonical homomorphism $\pi^{\cS(\eT(C)) \cap \cJ}$ coincides with $\pi^\cG$and is, thus, fractal. But $\cG_\eta = (\cS(\eT(C)) \cap \cJ)_\eta$ is a proper subset of $\cS(\eT(C))_\eta \cap \cJ_\eta$ for the sequence $\eta(n) := 2n+1$, since the sequence $(P_{2n+1} K P_{2n+1})$ belongs to $\cS(\eT(C))_\eta \cap \cJ_\eta$ for each compact operator $K$. \hfill \qed
\end{example}
\subsection{$\cJ$-fractal algebras} \label{ss43.2}
The considerations in the previous subsection suggest the following definitions.
\begin{defi} \label{d43.8}
Let $\cA$ be a $C^*$-subalgebra of $\cF$. \\[1mm]
$(a)$ If $\cJ$ is a closed ideal of $\cA$ then $\cA$ is called {\em $\cJ$-fractal} \index{algebra!$\cJ$-fractal}
if the canonical homomorphism $\pi^\cJ : \cA \to \cA/\cJ$ is fractal. \\[1mm]
$(b)$ If $\cJ$ is a closed ideal of $\cF$ then $\cA$ is called {\em $\cJ$-fractal} if $\cA$ is $(\cA \cap \cJ)$-fractal and if $(\cA \cap \cJ)_\eta = \cA_\eta \cap \cJ_\eta$ for each strictly increasing sequence $\eta : \sN \to \sN$.
\end{defi}
Note that both definitions coincide if $\cJ$ is a closed ideal of $\cA$ {\em and} $\cF$. Fractality of an algebra $\cA$ in the sense of Definition \ref{df1.5} is just $\cA \cap \cG$-fractality of $\cA$ in sense of Definition \ref{d43.8}, and it coincides with $\cG$-fractality of $\cA$ by Theorem \ref{tf1.8} $(b)$.

The following results show that $\cJ$-fractality implies what one expects: A sequence in a $\cJ$-fractal algebra belongs to $\cJ$ or is invertible modulo $\cJ$ if and only if at least one of its subsequences has this property.
\begin{theo} \label{t43.9}
Let $\cJ$ be a closed ideal of $\cF$. A $C^*$-subalgebra $\cA$ of $\cF$ is $\cJ$-fractal if and only if the following implication holds for every sequence $\bA \in \cA$ and every strictly increasing sequence $\eta : \sN \to \sN$
\begin{equation} \label{e43.10}
R_\eta (\bA) \in \cJ_\eta \quad \Longrightarrow \quad \bA \in \cJ.
\end{equation}
\end{theo}
{\bf Proof.} Let $\cA$ be $\cJ$-fractal and $\bA \in \cA$ a sequence with $R_\eta (\bA) \in \cJ_\eta$. Then
$R_\eta (\bA) \in \cA_\eta \cap \cJ_\eta = (\cA \cap \cJ)_\eta$, and the $(\cA \cap \cJ)$-fractality of $\cA$ implies $\bA \in \cJ$ via Theorem \ref{t43.1}.

Conversely, let (\ref{e43.10}) hold for each strictly increasing sequence $\eta$. From Theorem \ref{t43.1} we conclude that $\cA$ is $(\cA \cap \cJ)$-fractal. Further, the inclusion $\subseteq$ in $(\cA \cap \cJ)_\eta = \cA_\eta \cap \cJ_\eta$ is obvious. For the reverse inclusion, let $\bA$ be a sequence in $\cF$ with $R_\eta (\bA) \in \cA_\eta \cap \cJ_\eta$. Then there are sequences $\bB \in \cA$ and $\bJ \in \cJ$ such that $R_\eta (\bA) = R_\eta (\bB) = R_\eta (\bJ)$. Since $R_\eta (\bB) \in \cJ_\eta$, the implication (\ref{e43.10}) gives $\bB \in \cJ$. Hence, $R_\eta (\bB) \in (\cA \cap \cJ)_\eta$, and since $R_\eta (\bB) = R_\eta (\bA)$, one also has $R_\eta (\bA) \in (\cA \cap \cJ)_\eta$. \hfill \qed
\begin{theo} \label{t43.11}
Let $\cJ$ be a closed ideal of $\cF$ and $\cA$ a $\cJ$-fractal and unital $C^*$-subalgebra of $\cF$. Then the following implication holds for every sequence $\bA \in \cA$ and every strictly increasing sequence $\eta : \sN \to \sN$
\begin{equation} \label{e43.12}
R_\eta (\bA) + \cJ_\eta \; \mbox{invertible in} \; \cF_\eta/\cJ_\eta \quad \Longrightarrow \quad \bA + \cJ \;
\mbox{invertible in} \; \cF/\cJ.
\end{equation}
\end{theo}
{\bf Proof.} Let $\bA \in \cA$ be such that $R_\eta (\bA)  + \cJ_\eta$ is invertible in $\cF_\eta/\cJ_\eta$. Since $C^*$-algebras are inverse closed, this coset is also invertible in $(\cA_\eta + \cJ_\eta)/\cJ_\eta$. The latter algebra is $^*$-isomorphic to $\cA_\eta/(\cA_\eta \cap \cJ_\eta)$ by the third isomorphy theorem, hence, to $\cA_\eta/(\cA \cap \cJ)_\eta$ by $\cJ$-fractality of $\cA$. Thus, the coset $R_\eta (\bA) + (\cA \cap \cJ)_\eta$ is invertible in $\cA_\eta/ (\cA \cap \cJ)_\eta$. Choose sequences $\bB \in \cA$ and $\bJ \in \cA \cap \cJ$ such that
\[
R_\eta (\bA) \, R_\eta(\bB) = R_\eta(\bI) + R_\eta(\bJ)
\]
where $\bI$ stands for the identity element of $\cF$. Applying the homomorphism $\pi^{\cA \cap \cJ}_\eta$ to both sides of this equality one gets
\[
\pi^{\cA \cap \cJ} (\bA) \, \pi^{\cA \cap \cJ} (\bB) = \pi^{\cA \cap \cJ} (\bI) + \pi^{\cA \cap \cJ} (\bJ)
\]
which shows that $\bA \bB - \bI \in \cJ$. Hence, $\bA$ is invertible modulo $\cJ$ from the right-hand side. The
invertibility from the left-hand side follows analogously. \hfill \qed
\begin{coro} \label{c43.13}
Let $\cJ$ be a closed ideal of $\cF$ and $\cA$ a $\cJ$-fractal and unital $C^*$-subalgebra of $\cF$. Then a sequence $\bA \in \cA$ \\[1mm]
$(a)$ belongs to $\cJ$ if and only if there is a strictly increasing sequence $\eta$ such that $\bA_\eta$ belongs to $\cJ_\eta$. \\[1mm]
$(b)$ is invertible modulo $\cJ$ if and only there is a strictly increasing sequence $\eta$ such that $\bA_\eta$ is invertible modulo $\cJ_\eta$.
\end{coro}
Another feature of $\cJ$-fractal algebras is that quotient norms are preserved when passing to subsequences.
\begin{theo} \label{t43.14}
Let $\cA$ be a $C^*$-subalgebra of $\cF$ and $\cJ$ a closed ideal of $\cA$. If $\cA$ is $\cJ$-fractal, then
\[
\|R_\eta (\bA) + \cJ_\eta\| = \|\bA + \cJ\|
\]
for every sequence $\cA \in \cA$ and every strictly increasing sequence $\eta$.
\end{theo}
{\bf Proof.} For arbitrary sequences $\bA \in \cA$ and $\bJ \in \cJ$,
\[
\|R_\eta (\bA) + \cJ_\eta\| \le \|R_\eta (\bA) + R_\eta (\bJ)\| \le \|\bA + \bJ\|.
\]
Taking the infimum over all sequences $\bJ \in \cJ$ yields $\|R_\eta (\bA) + \cJ_\eta\| \le \|\bA + \cJ\|$. The reverse estimate follows from \begin{eqnarray*}
\|\bA + \cJ\| & = & \|\bA + \bJ + \cJ\| = \|\pi^\cJ (\bA + \bJ)\| \\
& = & \| \pi^\cJ_\eta R_\eta (\bA + \bJ)\| \le \|R_\eta (\bA) + R_\eta (\bJ)\|,
\end{eqnarray*}
which holds for arbitrary sequences $\bJ \in \cJ$, by taking again the infimum over all sequences $\bJ \in \cJ$ on the right-hand side. \hfill \qed
\begin{prop} \label{p43.15}
Let $\cJ$ be a closed ideal of $\cF$ and $\cA$ a $\cJ$-fractal $C^*$-subalgebra of $\cF$. Then \\[1mm]
$(a)$ every $C^*$-subalgebra of $\cA$ is $\cJ$-fractal. \\[1mm]
$(b)$ if $\cI$ is an ideal of $\cF$ with $\cJ \subseteq \cI$ and $(\cA \cap \cI)_\eta = \cA_\eta \cap \cI_\eta$ for every strictly increasing sequence $\eta$, then $\cA$ is $\cI$-fractal.
\end{prop}
{\bf Proof.} $(a)$ Let $\cB$ be a $C^*$-subalgebra of $\cA$, and let $\bB$ be a sequence in $\cB$ with $R_\eta(\bB) \in \cJ_\eta$ for a certain strictly increasing sequence $\eta$. Then $R_\eta(\bB) \in \cA_\eta \cap \cJ_\eta$. Since $\cA$ is $\cJ$-fractal, Theorem \ref{t43.9} implies that $\bB \in \cJ$. Hence $\cB$ is $\cJ$-fractal, again by Theorem \ref{t43.9}. \\[1mm]
$(b)$ Let $R_\eta (\bA) \in \cI_\eta$ for a sequence $\bA \in \cA$ and a strictly increasing sequence $\eta$. By
hypothesis, $R_\eta (\bA) \in (\cA \cap \cI)_\eta$. Choose a sequence $\bJ \in \cA \cap \cI$ with $R_\eta (\bA) = R_\eta (\bJ)$. The $\cJ$-fractality of $\cA$ implies that $\bA - \bJ \in \cJ$, whence $\bA \in \bJ + \cJ \subseteq \cI$. By Theorem \ref{t43.9}, $\cA$ is $\cI$-fractal. \hfill \qed
\subsection{Essential fractality and Fredholm property} \label{ss43.3}
Let $\cF$ be the algebra of matrix sequences with dimension function $\delta$ and $\cK$ the associated ideal of compact sequences. We will call the $\cK$-fractal $C^*$-subalgebras of $\cF$ {\em essentially fractal}. \index{algebra!essentially fractal} Note that every restriction $\cF_\eta$ of $\cF$ is again an algebra of matrix sequences (with dimension function $\delta \circ \eta$); hence, the restriction $\cK_\eta$ of $\cK$ is just the ideal of the compact sequences related with $\cF_\eta$. If we speak on {\em compact subsequences} and {\em Fredholm subsequences} in what follows, we thus mean sequences $R_\eta \bA \in \cK_\eta$ and sequences $R_\eta \bA$ which are invertible modulo $\cK_\eta$, respectively. In these terms, Corollary \ref{c43.13} reads as follows.
\begin{coro} \label{c43.16}
Let $\cA$ be an essentially fractal and unital $C^*$-subalgebra of $\cF$. Then a sequence $\bA \in \cA$ is compact (Fredholm) and only if one of its subsequences is compact (Fredholm), respectively.
\end{coro}
The following is a consequence of Proposition \ref{p43.15}.
\begin{coro} \label{c43.16a}
Let $\cA$ be a fractal $C^*$-subalgebra of $\cF$. If $(\cA \cap \cK)_\eta = \cA_\eta \cap \cK_\eta$ for every strictly increasing sequence $\eta$, then $\cA$ is essentially fractal.
\end{coro}
The finite sections algebra for band-dominated operators is a real life example for an algebra which is essentially fractal but not fractal (Theorem \ref{t73.27} below). The following example shows an algebra which is fractal, but not essentially fractal.
\begin{example} \label{ex43.16a}
Define $(A_n) \in \cF$ by
\[
A_n := \left\{
\begin{array}{ll}
\mbox{diag} \, (0, \, 0, \, \ldots 0, \, 1) & \mbox{if $n$ is even} \\
\mbox{diag} \, (0, \, 1, \, \ldots 1, \, 1) & \mbox{if $n$ is odd.}
\end{array}
\right.
\]
The sequence $(A_n)$ is fractal by Theorem \ref{tf3.2}, but it is not essential fractal: its subsequence $(A_{2n})$ is compact whereas its subsequence $(A_{2n+1})$ is Fredholm. \hfill \qed
\end{example}
Essential fractality has striking consequences for the behavior of the smallest singular values. Let again $\sigma_1(A) \le \ldots \le \sigma_n(A)$ denote the increasingly ordered singular values of an $n \times n$-matrix $A$.
\begin{theo} \label{t43.17}
Let $\cA$ be an essentially fractal and unital $C^*$-subalgebra of $\cF$. A sequence $(A_n) \in \cA$ is Fredholm if and only if there is a $k \in \sN$ such that
\begin{equation} \label{e43.18}
\limsup_{n \to \infty} \sigma_k (A_n) > 0.
\end{equation}
\end{theo}
{\bf Proof.} If $(A_n)$ is Fredholm then, by Theorem \ref{t42.1} $(c)$, $\liminf_{n \to \infty} \sigma_k \, (A_n) > 0$ for some $k \in \sN$, whence (\ref{e43.18}). Conversely, let (\ref{e43.18}) hold for some $k$. We choose a strictly increasing sequence $\eta$ such that $\lim_{n \to \infty} \, \sigma_k (A_{\eta(n)}) > 0$. Thus, the restricted sequence $(A_{\eta(n)})_{n \ge 1}$ is Fredholm by Theorem \ref{t42.1}. Since $\cA$ is essentially fractal, Corollary \ref{c43.16} $(b)$ implies the Fredholm property of the sequence $(A_n)$ itself. \hfill \qed \\[3mm]
Consequently, if a sequence $(A_n)$ in an essentially fractal and unital $C^*$-subalgebra of $\cF$ is {\em not} Fredholm, then
\begin{equation} \label{e43.19}
\lim_{n \to \infty} \, \sigma_k(A_n) = 0 \qquad \mbox{for each} \; k \in \sN.
\end{equation}
In analogy with operator theory, we call a sequence $(A_n)$ with property (\ref{e43.19}) {\em not normally solvable}. \index{sequence!not normally solvable}
\begin{coro} \label{c43.20}
Let $\cA$ be an essentially fractal and unital $C^*$-subalgebra of $\cF$. Then a sequence in $\cA$ is either Fredholm or not normally solvable.
\end{coro}
\begin{example} \label{ex43.21}
Consider the finite sections algebra $\cS(\eT(C))$ for Toeplitz operators. Using the description of $\cS(\eT(C))$ given in Theorem \ref{t14.21}, one easily checks that $(\cS(\eT(C)) \cap \cK)_\eta = \cS(\eT(C))_\eta \cap \cK_\eta$ for each strictly increasing sequence $\eta$. Since $\cS(\eT(C))$ is fractal and $\cG \subset \cK$, the algebra $\cS(\eT(C))$ is essentially fractal by Corollary \ref{c43.16a}. \hfill \qed
\end{example}
\subsection{Essential fractality of $\cS(\eB \eD \eO(\sN))$}
A bounded linear operator $A$ on $l^2(\sN)$ is said to be a {\em band operator} if its matrix representation $(a_{ij})$ with respect to the standard basis of $l^2(\sN)$ is a band matrix, i.e., if there is a constant $M$ such that $a_{ij} = 0$ whenever $|i-j| > M$. The norm closure of the set of all band operators on $l^2(\sN)$ is a $C^*$-subalgebra of $L(l^2(\sN))$ which we denote by $\eB \eD \eO(\sN)$ and the elements of which we call {\em band-dominated operators}.

To discretize the algebra $\eB \eD \eO(\sN)$ we choose the same filtration $\cP = (P_n)$ as for the Toeplitz algebra and consider the smallest closed subalgebra $\cS(\eB \eD \eO(\sN))$ of $\cF^\cP$ which contains all sequences $(P_n A P_n)$ with $A$ band-dominated. It is not difficult to derive a stability criterion for sequences in $\cS(\eB \eD \eO(\sN))$. It rests on the observation that a sequence $\cA = (A_n) \in \cF^\cP$ is stable if and only if the associated operator
\[
\op (\bA) := \diag (A_1, \, A_2, \, A_3, \, \ldots) \in L(l^2(\sN))
\]
is a Fredholm operator. In general, this observation is of less use, but if $\bA
\in \cS(\eB \eD \eO(\sN))$, then $\op (\bA)$ is a band-dominated operator. Thus,
the known limit operators criterion for the Fredholm property of band-dominated
operators implies a criterion for the stability of sequences in $\cS(\eB \eD
\eO(\sN))$. For details see Chapter 6 in \cite{RRS4} and \cite{Lin2,RRS1,Roc10}.
One should mention that Roe \cite{Roe2} generalized the limit operators Fredholm
criterion to band-dominated operators on general exact discrete groups, which
was used in \cite{RaR9} to derive a stability criterion for certain finite
sections methods for these operators.

Now we turn to fractality properties of the algebra $\cS(\eB \eD \eO(\sN))$. We have already seen an example which shows that this algebra fails to be fractal.
\begin{theo} \label{t73.27}
The algebra $\cS(\eB \eD \eO(\sN))$ is essentially fractal.
\end{theo}
{\bf Sketch of the proof.} Let $\bA = (A_n)$ be a sequence in $\cS(\eB \eD \eO(\sN))$ with $R_\eta (\bA) = (A_{\eta(n)}) \in \cK_\eta$ for some strictly increasing sequence $\eta$. Let $A$ denote the strong limit of the sequence $(A_n)$, and write $\bA$ as $(P_n A P_n) + (K_n)$. Then $A$ is a band-dominated operator, and the sequence $(K_n)$ lies in the kernel of the consistency map. This kernel is a closed ideal in $\cS(\eB \eD \eO(\sN))$ which is generated (as an ideal) by sequences of the form $(P_n B P_n C P_n) - (P_n BC P_n) =: (P_n B Q_n C P_n)$ with band operators $B, \, C$ by Theorem 2.15 in \cite{Roc10}. It is an easy exercise to check that all sequences of this form are compact. Thus, letting $n$ go to infinity in
\[
A_{\eta(n)} = P_{\eta(n)} A P_{\eta(n)} + K_{\eta(n)}
\]
yields the compactness of $A$ by Proposition 4.2 $(a)$ in \cite{Roc10}. Consequently, $(P_n A P_n)$ is a compact sequence. Hence, the sequence $\bA$ itself is compact. By Theorem \ref{t43.9}, the algebra $\cS(\eB \eD \eO(\sN))$ is essentially fractal. \hfill \qed
\subsection{Essential fractal restriction} \label{ss43.4}
Our next goal is an analogue of Theorem \ref{tf3.20} for essential fractality. Recall that we based the proof of Theorem \ref{tf3.20} on the fact that there is a sequence $\eta$ such that the norms $\|A_{\eta(n)}\|$ converge for each sequence $(A_n)$. We start with showing that $\eta$ can be even chosen such that not only the sequences
$(\|A_{\eta(n)}\|) = (\Sigma_1(A_{\eta(n)}))$ converge, but {\em every} sequence $(\Sigma_k(A_{\eta(n)}))$ with $k \in \sN$. Again, $\Sigma_1(A) \ge \ldots \ge \Sigma_n(A)$ denote the decreasingly ordered singular values of an $n \times n$-matrix $A$.
\begin{prop} \label{p43.22}
Let $\cA$ be a separable $C^*$-subalgebra of $\cF$. Then there is a strictly increasing sequence $\eta : \sN \to \sN$ such that the sequence $(\Sigma_k(A_{\eta(n)}))_{n \ge 1}$ converges for every sequence $(A_n)_{n \ge 1} \in \cA$ and every $k \in \sN$.
\end{prop}
{\bf Proof.} First consider a single sequence $(A_n) \in \cA$. We choose a strictly increasing sequence $\eta_1 : \sN \to \sN$ such that the sequence $(\Sigma_1(A_{\eta_1(n)}))_{n \ge 1}$ converges, then a subsequence $\eta_2$ of $\eta_1$ such that the sequence $(\Sigma_2(A_{\eta_2(n)}))_{n \ge 1}$ converges, and so on. The sequence $\eta(n) := \eta_n(n)$ has the property that the sequence $(\Sigma_k(A_{\eta(n)}))_{n \ge 1}$ converges for every $k \in \sN$.

Now let $(\bA^m)_{m \ge 1}$ be a countable dense subset of $\cA$, consisting of sequences $\bA^m = (A_n^m)_{n \ge 1}$.
We use the result of the previous step to find a strictly increasing sequence $\eta_1 : \sN \to \sN$ such that the sequences $(\Sigma_k(A^1_{\eta_1(n)}))_{n \ge 1}$ converge for every $k \in \sN$, then a subsequence $\eta_2$ of $\eta_1$ such that the sequences $(\Sigma_k(A^2_{\eta_2(n)}))_{n \ge 1}$ converge for every $k$, and so on. Then the sequence $\eta(n) := \eta_n(n)$ has the property that the sequences $(\Sigma_k(A^m_{\eta(n)}))_{n \ge 1}$ converge for every pair $k, \, m \in \sN$.

Let $\eta$ be as in the previous step, i.e., the sequences $(\Sigma_k(A^m_{\eta(n)}))_{n \ge 1}$ converge for every $k \in \sN$ and for every sequence $\bA^m = (A_n^m)_{n \ge 1}$ in a countable dense subset of $\cA$. We show that then the sequences $(\Sigma_k(A_{\eta(n)}))_{n \ge 1}$ converge for every $k \in \sN$ and every sequence $\bA = (A_n)$ in $\cA$. Fix $k \in \sN$ and let $\varepsilon > 0$. Using the well known inequality $|\Sigma_k(A) - \Sigma_k(B)| \le \|A-B\|$ we obtain
\begin{eqnarray*}
\lefteqn{|\Sigma_k(A_{\eta(n)}) - \Sigma_k(A_{\eta(l)})|} \\
&& \le |\Sigma_k(A_{\eta(n)}) - \Sigma_k(A^m_{\eta(n)})| + |\Sigma_k(A^m_{\eta(n)}) - \Sigma_k(A^m_{\eta(l)})| \\
&& \quad + \; |\Sigma_k(A^m_{\eta(l)}) - \Sigma_k(A_{\eta(l)})| \\
&& \le \|A_{\eta(n)} - A^m_{\eta(n)}\| + |\Sigma_k(A^m_{\eta(n)}) - \Sigma_k(A^m_{\eta(l)})| + \|A^m_{\eta(l)} - A_{\eta(l)}\| \\
&& \le 2 \, \|\bA - \bA^m\|_\cF + |\Sigma_k(A^m_{\eta(n)}) - \Sigma_k(A^m_{\eta(l)})|.
\end{eqnarray*}
Now choose $m \in \sN$ such that $\|\bA - \bA^m\|_\cF < \varepsilon/3$ and then $N \in \sN$ such that $|\Sigma_k(A^m_{\eta(n)}) - \Sigma_k(A^m_{\eta(l)})| < \varepsilon/3$ for all $n, \, l \ge N$. Then $|\Sigma_k(A_{\eta(n)}) - \Sigma_k(A_{\eta(l)})| < \varepsilon$ for all $n, \, l \ge N$. Thus, $(\Sigma_k(A_{\eta(n)}))_{n \ge 1}$ is a Cauchy sequence, hence convergent. \hfill \qed
\begin{prop} \label{p43.23}
Let $\cA$ be a $C^*$-subalgebra of $\cF$ with the property that the sequences $(\Sigma_k(A_n))_{n \ge 1}$ converge for every sequence $(A_n) \in \cA$ and every $k \in \sN$. Then $\cA$ is essentially fractal.
\end{prop}
{\bf Proof}. Let $\bK = (K_n) \in \cA$ and let $\eta: \sN \to \sN$ be a strictly increasing sequence such that $\bK_\eta \in \cK_\eta$. Then, by assertion $(b)$ of Theorem \ref{t41.12},
\begin{equation} \label{e43.24}
\lim_{k \to \infty} \limsup_{n \to \infty} \Sigma_k (K_{\eta(n)}) = 0.
\end{equation}
By hypothesis, $\limsup_{n \to \infty} \Sigma_k (K_{\eta(n)}) = \lim _{n \to \infty} \Sigma_k (K_n)$. Hence, (\ref{e43.24}) implies $\lim_{k \to \infty} \lim_{n \to \infty} \Sigma_k (K_n) = 0$, whence $\bK \in \cK$ by assertion $(b)$ of Theorem \ref{t41.12}. Thus, every sequence in $\cA$ which has a compact subsequence is compact itself. Thus $\cA$ is essentially fractal by Theorem \ref{t43.9}. \hfill \qed
\begin{theo} \label{t43.25}
Let $\cA$ be a separable $C^*$-subalgebra of $\cF$. Then there is a strictly increasing sequence $\eta : \sN \to \sN$
such that the restricted algebra $\cA_\eta = R_\eta \cA$ is essentially fractal.
\end{theo}
Indeed, if $\eta$ is as in Proposition \ref{p43.22}, then the restriction $\cA_\eta$ is essentially fractal by Proposition \ref{p43.23}. \hfill \qed \\[3mm]
We know from Theorems \ref{tf3.20} and \ref{t43.25} that every separable $C^*$-subalgebra of $\cF$ has both a fractal and an essentially fractal restriction. It is an open question if this fact generalizes to arbitrary closed ideals $\cJ$ of $\cF$ in place of $\cG$ or $\cK$, i.e., if can one always force $\cJ$-fractality by a suitable restriction?
\subsection{Essential spectra of self-adjoint sequences} \label{ss44.1}
This section is addressed to the following question: Can one discover the {\em essential spectrum} \index{spectrum!essential} of a sequence $\bA = (A_n) \in \cF$ (i.e., the spectrum of the coset $\bA + \cK$, considered as an element of the quotient algebra $\cF/\cK$) from the behavior of the eigenvalues of the matrices $A_n$?

Given a self-adjoint matrix $A$ and a subset $M$ of $\sR$, let $N(A, \, M)$ \index{$N(A, \, M)$} denote the number of eigenvalues of $A$ which lie in $M$, counted with respect to their multiplicity. If $M = \{\lambda\}$ is a singleton, we write $N(A, \, \lambda)$ in place of $N(A, \, \{\lambda\})$. Thus, if $\lambda$ is an eigenvalue of $A$, then $N(A, \, \lambda)$ is its multiplicity.

Let $\bA = (A_n) \in \cF$ be a self-adjoint sequence. Following Arveson \cite{Arv7,Arv3,Arv4}, a point $\lambda \in \sR$ is called {\em essential} \index{point!essential} for this sequence if, for every open interval $U$ containing $\lambda$,
\[
\lim_{n \to \infty} N(A_n, \, U) = \infty,
\]
and $\lambda \in \sR$ is called {\em transient} \index{point!transient} for $\bA$ if there is an open interval $U$ which contains $\lambda$ such that
\[
\sup_{n \in \sN} \, N(A_n , \, U) < \infty.
\]
Thus, $\lambda \in \sR$ is {\em not} essential for $\bA$ if and only if $\lambda$ is transient for a subsequence of $\bA$, and $\lambda$ is {\em not} transient for $\bA$ if and only if $\lambda$ is essential for a subsequence of $(A_n)$. Moreover, if a point $\lambda$ is transient (resp. essential) for $\bA$, then is is also transient (resp. essential) for every subsequence of $\bA$.
\begin{theo} \label{t44.1}
Let $\bA \in \cF$ be a self-adjoint sequence. A point $\lambda \in \sR$ belongs to the essential spectrum of $\bA$ if and only if it is not transient for the sequence $\bA$.
\end{theo}
{\bf Proof.} Let $\bA = (A_n)$ be a bounded sequence of self-adjoint matrices. First let $\lambda \in \sR \setminus \sigma (\bA + \cK)$. We set $B_n := A_n - \lambda I_n$ and have to show that 0 is transient for the sequence $(B_n)$. Since $\lambda \in \sR \setminus \sigma (\bA + \cK)$, the sequence $(B_n)$ is Fredholm. Let $k$ denote its  $\alpha$-number. By Theorem \ref{t42.1} $(c)$ and the definition of the $\alpha$-number,
\[
\liminf_{n \to \infty} \sigma_{k+1} (B_n) =: C > 0 \quad \mbox{and} \quad \liminf_{n \to \infty} \sigma_k (B_n) = 0.
\]
Let $U := (-C/2, \, C/2)$. Since the singular values of a self-adjoint matrix are just the absolute values of the
eigenvalues of that matrix, we conclude that $N(B_n, \, U) \le k$ for all sufficiently large $n$. Thus, 0 is transient.

Conversely, let $\lambda \in \sR$ be transient for $(A_n)$. We claim that $(A_n - \lambda I_n)$ is a Fredholm
sequence. By transiency, there is an interval $U = (\lambda - \varepsilon, \, \lambda + \varepsilon)$ with $\varepsilon > 0$ such that $\sup_{n \in \sN} N(A_n, \, U) =: k < \infty$. Let $T_n$ denote the orthogonal projection from $\sC^{\delta(n)}$ onto the $U$-spectral subspace of $A_n$. Then $\rank T_n$ is not greater than $k$. It is moreover obvious that the matrices $B_n := (A_n - \lambda P_n)(I - T_n) + T_n$ are invertible for all $n \in \sN$ and that their inverses are uniformly bounded by the maximum of $1/\varepsilon$ and $1$. Hence, $(B_n^{-1}) \in \cF$ and
\begin{equation} \label{e44.1a}
(A_n - \lambda P_n)(I - T_n) B_n^{-1} = I - T_n B_n^{-1}.
\end{equation}
Since $(T_n)$ is a compact sequence (of essential rank not greater than $k$), this identity shows that the coset $(A_n - \lambda I_n) + \cK$ is invertible from the right-hand side. Since this coset is self-adjoint, it is then invertible from both sides. Thus, $(A_n - \lambda I_n)$ is a Fredholm sequence. \hfill \qed
\begin{prop} \label{p44.1a1}
The set of the non-transient points and the set of the essential points of a self-adjoint sequence $\bA \in \cF$ are compact.
\end{prop}
{\bf Proof.} The first assertion is an immediate consequence of Theorem \ref{t44.1}. The second assertion will follow once we have shown that the set of the essential points of $\bA$ is closed.

Let $(\lambda_k)$ be a sequence of essential points for $\bA = (A_n)$ with limit $\lambda$. Assume that $\lambda$ is not essential for $\bA$. Then there is a strictly increasing sequence $\eta : \sN \to \sN$ such that $\lambda$ is transient for $\bA_\eta$. Let $U$ be an open neighborhood of $\lambda$ with $\sup_{n \in \sN} \, N(A_{\eta(n)}, \, U) =: c < \infty$. Since $\lambda_k \to \lambda$ and $U$ is open, there are a $k \in \sN$ and an open neighborhood $U_k$ of $\lambda_k$ with $U_k \subseteq U$. Clearly, $N(A_{\eta(n)}, \, U_k) \le N(A_{\eta(n)}, \, U) \le c$. On the other hand, since $\lambda_k$ is also essential for the restricted sequence $\bA_\eta$, one has $N(A_{\eta(n)}, \, U_k) \to \infty$ as $n \to \infty$, a contradiction. \hfill \qed \\[3mm]
Note that the set of the non-transient points of a self-adjoint sequence is non-empty by Theorem \ref{t44.1}, whereas it is easy to construct self-adjoint sequences without any essential point: take a sequence which alternates between the zero and the identity matrix. In contrast to this observation, the following result shows that sequences which arise by discretization of a self-adjoint operator, always possess essential points. Let $H$ be an infinite dimensional separable Hilbert space with filtration $\cP := (P_n)$.
\begin{theo} \label{t44.2}
Let $\bA := (A_n) \in \cF^\cP$ be a self-adjoint sequence with strong limit $A$. Then every point in the essential spectrum of $A$ is an essential point for $\bA$.
\end{theo}
{\bf Proof.} We show that $A - \lambda I$ is a Fredholm operator if $\lambda \in \sR$ is not essential for $\bA$. Then $\lambda$ is transient for a subsequence of $\bA$, i.e., there are an infinite subset $\sM$ of $\sN$ and an interval $U = (\lambda - \varepsilon, \, \lambda + \varepsilon)$ with $\varepsilon > 0$ such that
\begin{equation} \label{e44.2a}
\sup_{n \in \sM} N(A_n, \, U) =: k < \infty.
\end{equation}
Let $T_n$ denote the orthogonal projection from $H$ onto the $U$-spectral subspace of $A_n P_n$. By (\ref{e44.2a}), the rank of the projection $T_n$ is not greater than $k$ if $n \in \sM$. So we conclude as in the proof of Theorem \ref{t44.1} that the operators
\[
B_n := (A_n - \lambda P_n)(I - T_n) + T_n
\]
are invertible for all $n \in \sM$ and that their inverses are uniformly bounded by the maximum of $1/\varepsilon$ and $1$. Hence,
\begin{equation} \label{e44.3}
(A_n - \lambda P_n)(I - T_n) B_n^{-1} = I - T_n B_n^{-1}
\end{equation}
for all $n \in \sM$. By the weak sequential compactness of the unit ball of $L(H)$, one finds weakly convergent subsequences $((I - T_{n_r}) B_{n_r}^{-1})_{r \ge 1}$ of $((I - T_n) B_n^{-1}) _{n \in \sM}$ and $(T_{n_r} B_{n_r}^{-1})_{r \ge 1}$ of $(T_n B_n^{-1})_{n \in \sM}$ with limits $B$ and $T$, respectively. The product of a weakly convergent sequence with limit $C$ and a $^*$-strongly convergent sequence with limit $D$ is weakly convergent with limit $CD$. Thus, passing to subsequences and taking the weak limit in (\ref{e44.3}) yields $(A - \lambda I) B = I - T$. Further, the rank of $T$ is not greater than $k$ by Proposition 4.2 $(a)$ in \cite{Roc10}. Thus, $(A - \lambda I) B - I$ is a compact operator. The compactness of $B (A - \lambda I) - I$ follows similarly. Hence, $A$ is a Fredholm operator. \hfill \qed \\[3mm]
Arveson gave a first example where the inclusion in Theorem \ref{t44.2} is proper. Specifically, he constructed a self-adjoint unitary operator $A \in L(l^2(\sN))$ with
\begin{equation}\label{e44.4}
\sigma(A) = \sigma_{ess} (A) = \{-1, \, 1\}
\end{equation}
such that 0 is an essential point of the sequence $(P_nAP_n)$. For completeness, Arveson's example is restated below.
\begin{example} \label{ex44.5}
Let $E_1 := 4 \sN$ and $E_2 := 2 \sN \setminus 4 \sN$, and define a function $f$ on $E_1$ by $f(k) = k^2+1$. Set $O_1 := f(E_1)$ and $O_2:=(2 \sN-1) \setminus O_1$. Clearly, both $E_2$ and $O_2$ are infinite subsets of $\sN$. Let $f$ be any bijection from $E_2$ onto $O_2$. This construction implies a permutation $\pi$ of $\sN$ with $\pi^2$ being the identity via
\[
\pi(k) := \left\{
\begin{array}{ll}
f(k)      & {\rm if} \; k \; {\rm even} \\
f^{-1}(k) & {\rm if} \; k \; {\rm odd}.
\end{array}  \right.
\]
We claim that the operator $A : l^2(\sN) \to l^2(\sN), \; (a_n)_{n \in \sN} \mapsto (a_{\pi(n)})_{n \in \sN}$ has the announced properties. Evidently, $A = A^* = A^{-1} \neq \pm I$, whence (\ref{e44.4}) follows. In order to see that 0 is an essential point of $(P_nAP_n)$, let $\sN_n$ denote the set $\{1, \, 2, \, \dots, \, n\}$ and write $\sharp S$ for the number of the elements of a set $S$. It is elementary to check that $\sharp (f(E_1 \cap \sN_n) \setminus \sN_n)$ tends to infinity as $n \to \infty$, which implies that $\lim_{n \to \infty} \, \sharp (f(E \cap \sN_n) \setminus \sN_n) = \infty$ and, consequently,
\begin{equation} \label{e44.6}%
\lim_{n \to \infty} \, \sharp (\pi(\sN_n) \setminus \sN_n) = \infty.
\end{equation}
Since $A$ maps the basis element $e_k := (0, \, \dots, \, 0, \, 1, \, 0, \, \dots)$ of $l^2(\sN)$ (with the $1$ standing at the $k$th place) to $e_{\pi(k)}$, it follows that $P_n A P_n e_k = 0$ for every $k$ belonging to the set $S_n := \{k \in \sN_n : \pi(k) \not\in \sN_n \}$. But the cardinality of $S_n$ tends to infinity as $n \to \infty$ due to (\ref{e44.6}). Hence, $0$ is an eigenvalue of $P_n A P_n$ for all sufficiently large $n$, and the multiplicity
of this eigenvalue tends to infinity. Consequently, $0$ is an essential point for $(P_nAP_n)$.
\hfill \qed
\end{example}
\subsection{Arveson dichotomy and essential fractality} \label{ss44.1a}
We say that a self-adjoint sequence $\bA \in \cF$ enjoys {\em Arveson's dichotomy} \index{Arveson's dichotomy} if every real number is either essential or transient for this sequence. Note that Arveson dichotomy is preserved when passing to subsequences. Arveson introduced and studied this property in a series of papers \cite{Arv7,Arv3,Arv4}. He proved the dichotomy of the finite sections sequence $(P_n A P_n)$ when $A$ is a self-adjoint band operator, and he extended this result to a class of self-adjoint band-dominated operators which satisfy a Wiener and a Besov space condition. It was shown in Theorem 7.6 in \cite{Roc10} that this result holds for arbitrary self-adjoint band-dominated operators. We will get this fact here by combining Theorems \ref{t73.27} and \ref{t44.1b} below.
\begin{theo} \label{t44.1a2}
The set of all self-adjoint sequences with Arveson dichotomy is closed in $\cF$.
\end{theo}
{\bf Proof.} Let $(\bA_n)_{n \in \sN}$ be a sequence of self-adjoint sequences in $\cF$ with Arveson dichotomy which converges to a (necessarily self-adjoint) sequence $\bA$ in the norm of $\cF$. Then $\bA_n + \cK \to \bA + \cK$ in the norm of $\cF/\cK$. Since $\bA_n + \cK$ and $\bA + \cK$ are self-adjoint elements of $\cF/\cK$, this implies that the spectra of $\bA_n + \cK$ converge to the spectrum of $\bA + \cK$ in the Hausdorff metric. Thus, by Theorem \ref{t44.1}, the sets of the non-transient points of $\bA_n$ converge to the set of the non-transient points of $\bA$. Since the $\bA_n$ have Arveson dichotomy by hypothesis, this finally implies that the sets of the essential points of $\bA_n$ converge to the set of the non-transient points of $\bA$ in the Hausdorff metric.

Let now $\lambda$ be a non-transient point for $\bA$ and assume that $\lambda$ is not essential for $\bA$. Then there is a strictly increasing sequence $\eta : \sN \to \sN$ such that $\lambda$ is transient for the restricted sequence $\bA_\eta$. As we have seen above, there is a sequence $(\lambda_n)$, where $\lambda_n$ is an essential point for $\bA_n$, with $\lambda_n \to \lambda$. Since the property of being an essential is preserved under passage to a subsequence, $\lambda_n$ is also essential for the restricted sequence $(\bA_n)_\eta$.

Since the sequences $(\bA_n)_\eta$ also have Arveson dichotomy and since $(\bA_n)_\eta \to \bA_\eta$ in the norm of $\cF_\eta$, we can repeat the above arguments to conclude that the sets $M_n$ of the essential points for $(\bA_n)_\eta$ converge to the set $M$ of the non-transient points for $\bA_\eta$ in the Hausdorff metric. Since $\lambda_n \in M_n$ by construction, this implies that $\lambda \in M$. This means that $\lambda$ in not transient for $\bA_\eta$, a contradiction. \hfill \qed \\[3mm]
We continue with a result which relates Arveson dichotomy with essential fractality.
\begin{theo} \label{t44.1b}
Let $\cA$ be a unital $C^*$-subalgebra of $\cF$. Then $\cA$ is essentially fractal if and only if every self-adjoint sequence in $\cA$ has Arveson dichotomy.
\end{theo}
{\bf Proof.} First let $\cA$ be essentially fractal. Let $\bA$ be a self-adjoint sequence in $\cA$ and $\lambda \in \sR$ a point which is not essential for $\bA$. Then $\lambda$ is transient for a subsequence of $\bA$, thus, $0$ is transient for a subsequence of $\bA - \lambda \bI$. From Theorem \ref{t44.1} we conclude that this subsequence has the Fredholm property. Then, by Corollary \ref{c43.16} $(b)$ and since $\cA$ is essentially fractal, the sequence $\bA - \lambda \bI$ itself is a Fredholm sequence. Thus, $0$ is transient for $\bA - \lambda \bI$ by Theorem \ref{t44.1} again, whence finally follows that $\lambda$ is transient for $\bA$. Hence, $\bA$ has Arveson dichotomy.

Now assume that $\cA$ is not essentially fractal. Then, by Theorem \ref{t43.9}, there are a sequence $\bA = (A_n) \in \cA$ and a strictly increasing sequence $\eta : \sN \to \sN$ such that the restricted sequence $\bA_\eta$ belongs to $\cK_\eta$ but $\bA \not\in \cK$. The self-adjoint sequence $\bA^* \bA$ has the same properties, i.e., $(\bA^* \bA)_\eta = \bA^*_\eta \bA_\eta \in \cK_\eta$, but $\bA^* \bA \not\in \cK$.

Since $\bA^*_\eta \bA_\eta \in \cK_\eta$, the essential spectrum of $\bA^*_\eta \bA_\eta$ (i.e., the spectrum of the coset $\bA^*_\eta \bA_\eta + \cK_\eta$ in $\cF_\eta/\cK_\eta$) consists of the point $0$ only. Thus, by Theorem \ref{t44.1}, $0$ is the only non-transient point for the restricted sequence $\bA^*_\eta \bA_\eta$.

Since $\bA^* \bA \not\in \cK$, there is a strictly increasing sequence $\mu : \sN \to \sN$ such that $\mu(\sN) \cap \eta(\sN) = \emptyset$ and $\bA^*_\mu \bA_\mu \not\in \cK_\mu$. Hence, the essential spectrum of $\bA^*_\mu \bA_\mu$ contains at least one point $\lambda \neq 0$, and this point is non-transient for $\bA^*_\mu \bA_\mu$ by Theorem \ref{t44.1} again. Hence, there is a subsequence $\nu$ of $\mu$ such that $\lambda$ is essential for $\bA^*_\nu \bA_\nu$, but $\lambda \neq 0$ is transient for $\bA^*_\eta \bA_\eta$ as we have seen above. Thus, $\lambda$ is neither transient nor essential for $\bA^* \bA$. Hence, the sequence $\bA^* \bA$ does not have Arveson dichotomy. \hfill \qed
\begin{coro} \label{c44.1a3}
Every self-adjoint sequence in $\cF$ possesses a subsequence with Arveson dichotomy.
\end{coro}
{\bf Proof.} Let $\bA$ be a self-adjoint sequence in $\cF$. The smallest closed subalgebra $\cA$ of $\cF$ which contains $\bA$ is separable. By Theorem \ref{t43.25}, there is an essentially fractal restriction $\cA_\eta$ of $\cA$.  Then $\bA_\eta$ is a subsequence of $\bA$ with Arveson dichotomy by the previous theorem. \hfill \qed
\section{Fractal algebras of compact sequences} \label{s4}
In this section we consider compact and Fredholm sequences in fractal algebras. The property of fractality has some striking consequences. For example, fractal ideals in $\cK$ are constituted of blocks which are isomorphic to the ideal of the compact operators on a Hilbert space. There will be also a nice formula for the alpha-number of a Fredholm sequence. We will be very brief in this section and omit many details and almost all proofs.
\subsection{Fractality and large singular values} \label{ss41}
First we will see that the singular values of fractal compact sequences behave as the singular values of compact operators on Hilbert space, i.e., the set of the singular values is countable and has 0 as its only possible accumulation point. Note that this does not hold for general compact sequences. For example, take an enumeration $(a_n)$ of the rational numbers in $[0, \, 1]$ and set
\[
K_n := a_n P_n P_1 P_n = \diag (a_n, \, 0, \, \ldots, \, 0).
\]
Then the sequence $(K_n)$ is compact (it consists of rank one matrices), but the spectrum of its coset $(K_n) + \cG$ is the closed interval $[0, \, 1]$.

Given an $n \times n$-matrix $A$, let again $\Sigma_1 (A) \ge \ldots \ge \Sigma_n (A) \ge 0$ denote the singular values of $A$, and write $\sigma_{sing} (A)$ for the set of the singular values of $A$. Since the singular values of $A$ are the eigenvalues of self-adjoint matrix $(A^*A)^{1/2}$, it is an immediate consequence of Proposition \ref{pf2.7} that, for each sequence $(A_n)$ in a fractal algebra, the sets $\sigma_{sing} (A_n)$ converge with respect to the Hausdorff metric. In particular,
\begin{equation} \label{e170309.1}
\limsup \sigma_{sing} (A_n) = \liminf \sigma_{sing} (A_n) = \sigma_{sing} \left( (A_n) + \cG \right) .
\end{equation}
If $(A_n) \in \cF$ is a fractal sequence, then the sequence $(\Sigma_1 (A_n))$ of the largest singular values of $A_n$ converges. This fact follows immediately from Proposition \ref{pf2.3} $(b)$ and the identity $\Sigma_1 (A_n) = \|A_n\|$. One cannot expect that the sequence of the second singular values $\Sigma_2 (A_n)$ converges, too. Indeed, the sequence defined by
\[
A_n := \left\{ \begin{array}{ll}%
\diag (1, \, 0, \, 0, \, \ldots, \, 0) & \mbox{if} \; n \; \mbox{is odd} \\
\diag (1, \, 1, \, 0, \, \ldots, \, 0) & \mbox{if} \; n \; \mbox{is even}
\end{array} \right.
\]
is fractal by Theorem \ref{tf3.2}, but the sequence of its second singular values alternates between 0 and 1 and has, thus, two accumulation points. In fact, one can show that the sequence $(\Sigma_2 (A_n))$ can possess {\em at most two} limiting points, at most one of which is different from $\lim \Sigma_1 (A_n)$. This fact holds more general.
\begin{prop} \label{p51.2}
If the sequence $(A_n) \in \cF$ is fractal, then the set
\[
\limsup_{n \to \infty} \, \{ \Sigma_1 (A_n), \, \ldots, \, \Sigma_k (A_n) \}
\]
contains at most $k$ elements.
\end{prop}
{\bf Proof.} Write $\Pi_j$ for the set of all partial limits of the sequence $(\Sigma_j (A_n))_{n \in \sN}$. We first verify that
\begin{equation} \label{e51.3}
\Pi_1 \cup \ldots \cup \Pi_k = \limsup_{n \to \infty} \,
\{\Sigma_1 (A_n), \, \ldots, \, \Sigma_k(A_n) \} \quad \mbox{for
every} \; k \in \sN.
\end{equation}
The inclusion $\subseteq$ is evident. Conversely, if $\lambda$ belongs to the right-hand side of (\ref{e51.3}), then there are a strictly increasing sequence $\eta : \sN \to \sN$ and numbers $k_n$ in $\{1, \, \ldots, \, k \}$ such that $\lambda = \lim_{n \to \infty} \Sigma_{k_n} (A_{\eta(n)})$. Since $k_n$ can take only finitely many values, there is a $k_0$ between 1 and $k$ and a subsequence $\mu$ of $\eta$ such that $\lambda = \lim_{n \to \infty} \Sigma_{k_0} (A_{\mu(n)})$. Hence, $\lambda \in \Pi_{k_0}$, what verifies (\ref{e51.3}).

We have already mentioned that $\Pi_1$ is a singleton. Next we show that for each $j \ge 1$ the difference $\Pi_{j+1} \setminus (\Pi_1 \cup \ldots \cup \Pi_j)$ contains at most one element. Assume there are points $\alpha$ and $\beta$ in $\Pi_{j+1} \setminus (\Pi_1 \cup \ldots \cup \Pi_j)$ with $\alpha > \beta$. Choose a subsequence $(\Sigma_{j+1} (A_{\eta(n)}))$ of $(\Sigma_{j+1} (A_n))$ which converges to $\beta$ as $n \to \infty$. Then $\alpha$ cannot belong to the partial limiting set $\limsup \sigma_{sing} (A_{\eta(n)})$. Indeed, if $\alpha \in \limsup \sigma_{sing} (A_{\eta(n)})$ then
\[
\alpha \in \limsup _{n \to \infty} \, \{\Sigma_1 (A_{\eta(n)}), \, \ldots, \, \Sigma_j (A_{\eta(n)}) \}
\]
due to monotonicity reasons. Then
\[
\alpha \in \limsup_{n \to \infty} \, \{\Sigma_1 (A_n) , \, \ldots, \, \Sigma_j (A_n) \} = \Pi_1 \cup \ldots \cup \Pi_j
\]
which was excluded. Hence, $\alpha \in \limsup \sigma_{sing} (A_n) \setminus \liminf \sigma_{sing} (A_n)$, which contradicts (\ref{e170309.1}). \hfill \qed \\[3mm]
Here is the announced result on singular values of fractal compact sequences.
\begin{theo} \label{t51.5}
Let $(K_n) \in \cK$ be a fractal sequence. Then the set
\[
\mbox{\rm h-lim} \, \sigma_{sing} (K_n) = \sigma_{sing} ((K_n) + \cG)
\]
is at most countable, it contains the point $0$, and $0$ is the only accumulation point of this set.
\end{theo}
{\bf Proof.} Let $\varepsilon > 0$. By Theorem \ref{t41.12} $(a)$, $\lim_{k \to \infty} \sup_{n \ge k} \Sigma_k (K_n) = 0$. Thus, there is a $k_0$ such that $\sup_{n \ge k} \Sigma_{k} (K_n) \le \varepsilon$ for each $k \ge k_0$, whence
\[
\limsup_{n \to \infty} \, \{ \Sigma_{k_0} (K_n), \, \ldots, \,
\Sigma_n (K_n) \} \subseteq [0, \, \varepsilon].
\]
Hence, every point in $\mbox{\rm h-lim} \, \sigma_{sing} (K_n) \setminus [0, \, \varepsilon]$ must lie in
\[
\limsup_{n \to \infty} \, \{ \Sigma_1 (K_n), \, \ldots, \,
\Sigma_{k_0-1} (K_n) \}
\]
which is a finite set by Proposition \ref{p51.2}. Consequently,
$\mbox{\rm h-lim} \, \sigma_{sing} (K_n)$ is at most countable and has 0 as only possible accumulation point. That $0$ indeed belongs to this set is a consequence of Corollary \ref{c41.17}. \hfill \qed
\subsection{Compact elements in $C^*$-algebras} \label{ss42}
There is a general notion of a compact element in a $C^*$-algebra $\cA$. A non-zero element $k$ of $\cA$ is said to be {\em of rank one} if, for each $a \in \cA$, there is a complex number $\mu$ such that $kak = \mu k$. We let $\cC(\cA)$ stand for the smallest closed subalgebra of $\cA$ which contains all elements of rank one. If such elements do not exist, we set $\cC(\cA) = \{0\}$. The elements of $\cC(\cA)$ are called compact.
Since the product of a rank one element with an arbitrary element of $\cA$ is zero or rank one again, $\cC(\cA)$ is a closed ideal of $\cA$. There are several equivalent descriptions of the ideal $\cC(\cA)$. To state the descriptions which are important in what follows, we need some more notation.

A $C^*$-algebra is called {\em elementary} if it is $^*$-isomorphic to the ideal $K(H)$ of the compact operators on some Hilbert space $H$. A $C^*$-algebra $\cJ$ is called {\em dual} \index{algebra!dual} if it is $^*$-isomorphic to a direct sum of elementary algebras. Thus, there is an index set $T$, for each $t \in T$ there is an elementary algebra $\cJ_t$, and $\cJ$ is $^*$-isomorphic to the $C^*$-algebra of all bounded functions $a$ which are defined on $T$, take a value $a(t)$ in $\cJ_t$ at $t \in T$, and which own the property that for each $\varepsilon > 0$, there are only finitely many $t \in T$ with $\|a(t)\| > \varepsilon$. An alternate way to think of dual algebras is the following. Let $\{\cJ_t\}_{t \in T}$ be a family of elementary ideals of a $C^*$-algebra $\cA$ with the property that $\cJ_s \cJ_t$ is the zero ideal whenever $s \neq t$. Then the smallest closed subalgebra of $\cA$ which contains all algebras $\cJ_t$ is a dual algebra, and each dual algebra is of this form.
\begin{theo} \label{t180309.1}
Let $\cA$ be a unital $C^*$-algebra and $\cJ$ a closed ideal of $\cA$. The following assertions are equivalent: \\[1mm]
$(a)$ $\cJ = \cC(\cJ)$. \\
$(b)$ $\cJ$ is a dual algebra. \\
$(c)$ The spectrum of every self-adjoint element of $\cJ$ is at most countable and has 0 as only possible accumulation point.
\end{theo}
For every dual ideal of a $C^*$-algebra there is a lifting theorem as follows. For a proof, see \cite{HRS2}. The first version of this theorem appeared in \cite{Sil1} and celebrates its 30th birthday in this year.
\begin{theo}[Lifting theorem for dual ideals] \label{tb1.6}
Let $\cA$ be a unital $C^*$-alge\-bra. For every element $t$ of a set $T$, let $\cJ_t$ be an elementary ideal of $\cA$ such that $\cJ_s \cJ_t = \{0\}$ whenever $s \neq t$, and let $W_t : \cA \to L(H_t)$ denote the irreducible representation of $\cA$ which extends the (unique up to unitary equivalence) irreducible representation of $\cJ_t$. Let further $\cJ$ stand for the smallest closed ideal of $\cA$ which contains all ideals $\cJ_t$. \\[1mm]
$(a)$ An element $a \in \cA$ is invertible if and only if the coset $a + \cJ$ is invertible in $\cA/\cJ$ and if all elements $W_t(a)$ are invertible in $\cB_t$. \\[1mm]
$(b)$ The separation property holds, i.e. $W_s (\cJ_t) = \{0\}$ whenever $s \neq t$. \\[1mm]
$(c)$ If $j \in \cJ$, then $W_t(j)$ is compact for every $t \in T$. \\[1mm]
$(d)$ If the coset $a + \cJ$ is invertible, then all operators
$W_t(a) \in L(H_t)$ are Fredholm, and there are at most finitely many of these operators which are not invertible.
\end{theo}
The following result is an immediate consequence of Theorems \ref{t51.5} and \ref{t180309.1}. This corollary implies that every unital and fractal $C^*$-subalgebra of $\cF$ which contains non-trivial compact sequences is a subject to the lifting theorem.
\begin{coro} \label{c190309.1}
Let $\cA$ be a unital and fractal $C^*$-subalgebra of $\cF$ which contains the ideal $\cG$. Then the ideal $(\cA \cap \cK)/\cG$ of $\cA/\cG$ is a dual
algebra.
\end{coro}
\subsection{Weights of elementary algebras of sequences} \label{ss43}
A projection in a $C^*$-algebra is a self-adjoint element $p$ with $p^2 = p$. We say that a closed ideal $\cJ$ of a $C^*$-algebra $\cA$ lifts projections, if every coset which is a projection in $\cA/\cJ$ contains a representative which is a projection in $\cA$. Note that, in general, closed ideals of $C^*$-algebras do not lift projections. For a simple example, take $\cA = C([0, \, 1])$ and $\cJ = \{f \in \cA: f(0) = f(1) = 0\}$. The following proposition shows that elementary ideals of $\cF/\cG$ lift projections. More general, every dual ideal has the projection lifting property.
\begin{prop} \label{p50.1}
Let $\cJ$ be an elementary $C^*$-subalgebra of $\cF/\cG$. \\[1mm]
$(a)$ Every projection $p \in \cJ$ lifts to a sequence $(\Pi_n) \in \cF$ of orthogonal projections, i.e., $(\Pi_n) + \cG = p$. \\[1mm]
$(b)$ If $p$ and $q$ are rank one projections in $\cJ$ which lift to sequences of projections $(\Pi^p_n)$ and $(\Pi^q_n)$, respectively, then
\[
\dim \im \Pi^p_n = \dim \im \Pi^q_n
\]
for all sufficiently large $n$.
\end{prop}
Thus, {\em for large} $n$, the entries of the sequence $(\dim \im \Pi^p_n)_{n \ge 1}$ are uniquely determined by the algebra $\cJ$; they do neither depend on the choice of the rank one projection $p$ nor on its lifting.

For a precise formulation, we define an equivalence relation $\sim$ in the set of all sequences of non-negative integers by calling two sequences $(\alpha_n), \, (\beta_n)$ {\em equivalent} if $\alpha_n = \beta_n$ for all sufficiently large $n$. Then Proposition \ref{p50.1} states that the equivalence class which contains the sequence $(\dim \im \Pi^p_n)_{n \ge 1}$ is uniquely determined by the algebra $\cJ$. We denote this equivalence class by $\alpha^\cJ$ and call it the {\em weight} of the elementary algebra $\cJ$. The algebra $\cJ$ is said to be an {\em algebra of positive weight} if the equivalence class $\alpha^\cJ$ contains a sequence consisting of positive numbers only, and $\cJ$ is an {\em algebra of weight one} if the equivalence class $\alpha^\cJ$ contains the constant sequence $(1, \, 1, \, \ldots)$. Note that the weight is {\em bounded} if $\cJ$ is in $\cK/\cG$, since then $(\Pi^p_n)$ is a compact sequence and has finite essential rank.
\subsection{Silbermann pairs and $\cJ$-Fredholm sequences} \label{ss44}
Next we are going to examine the consequences of the Lifting theorem \ref{tb1.6}
for sequence algebras. We will do this in the slightly more general context of
Silbermann pairs. A {\em Silbermann pair} $(\cA, \, \cJ)$ consists of a unital
$C^*$-subalgebra $\cA$ of $\cF$ and of a closed ideal $\cJ$ of $\cA$ which
contains $\cG$ properly and which consists of compact sequences only, and for
which $\cJ/\cG$ is a dual subalgebra of $\cK/\cG$. This property ensures that
the lifting theorem applies to Silbermann pairs. Every sequence in $\cA$ which
is invertible modulo $\cJ$ is called an {\em $\cJ$-Fredholm} sequence. Note that
each $\cJ$-Fredholm sequence is Fredholm in sense of Section \ref{ss33} (but, of
course, a Fredholm sequence in $\cA$ is not necessarily $\cJ$-Fredholm). Under
the conditions of Corollary \ref{c190309.1}, $(\cA, \cA \cap \cK)$ is a
Silbermann pair, and a sequence in $\cA$ is $(\cA \cap \cK)$-Fredholm if and
only if it is Fredholm. The study of Silbermann pairs (in the special case
when $\cJ/\cG$ is an elementary subalgebra of $\cK/\cG$) was initiated in
\cite{Sil8}.

Let $(\cA, \, \cJ)$ be a Silbermann pair. Then the algebra $\cJ/\cG$ is dual; hence, it is the direct sum of a family $(I_t)_{t \in T}$ of elementary algebras with associated bijective representations $W_t : I_t \to K(H_t)$. These representations extent to irreducible representations of $\cA$ into $L(H_t)$ which we denote by $W_t$ again. In this context, the Lifting theorem \ref{tb1.6} specifies as follows.
\begin{theo} \label{t50.5}
Let $(\cA, \, \cJ)$ form a Silbermann pair. \\[1mm]
$(a)$ A sequence $\bA \in \cA$ is stable if and only if it is $\cJ$-Fredholm and if the operators $W_t(\bA)$ are invertible for each $t \in T$. \\[1mm]
$(b)$ The separation property holds, i.e., $W_s(I_t) = \{0\}$ whenever $s \neq t$. \\[1mm]
$(c)$ If $\bJ \in \cJ$, then $W_t(\bJ)$ is a compact operator for every $t \in T$. \\[1mm]
$(d)$ If the sequence $\bA \in \cA$ is $\cJ$-Fredholm, then all operators $W_t(\bA)$ are Fredholm, and there are at most finitely many of these operators which are not invertible.
\end{theo}
For each $t \in T$, we choose and fix a representative $(\alpha^t_n)$ of the weight $\alpha^{I_t}$ of the elementary ideal $I_t$. Let the sequence $\bA := (A_n) \in \cA$ be $\cJ$-Fredholm. Assertion $(d)$ of the Lifting theorem \ref{t50.5} implies that the sum
\begin{equation} \label{e50.6}
\alpha_n(\bA) := \sum_{t \in T} \alpha_n^t \, \dim \ker W_t(\bA)
\end{equation}
is finite. Evidently, this definition depends on the choice of the representatives of the weight functions. But since only a finite number of items in the sum (\ref{e50.6}) is not zero, the equivalence class of the sequence $(\alpha_n(\bA))$ modulo $\sim$ is uniquely determined. Thus, the entries of that sequence are uniquely determined for sufficiently large $n$.

The main result of the present section is the following splitting property of the singular values of a $\cJ$-Fredholm sequence. The numbers $\sigma_k (A_n)$ with $1 \le k \le n$ denote again the increasingly ordered singular values of $A_n$.
\begin{theo} \label{t50.7}
Let $(\cA, \, \cJ)$ be a Silbermann pair, and let the sequence $\bA = (A_n)$ be $\cJ$-Fredholm. Then $\bA$ is a Fredholm sequence, and
\begin{equation} \label{e50.8}
\lim_{n \to \infty} \sigma_{\alpha_n(\bA)} (A_n) = 0 \quad \mbox{whereas} \quad \liminf_{n \to \infty} \sigma_{\alpha_n(\bA) +1} (A_n) > 0.
\end{equation}
\end{theo}
The proof makes use of results on lifting of families of mutually orthogonal projections and on generalized (or Moore-Penrose) invertibility. For details see \cite{Roc9}.

Theorem \ref{t50.7} has some remarkable consequences. First note that the number
\begin{equation} \label{e50.16}
\alpha(\bA) := \limsup_{n \to \infty} \alpha_n (\bA)
\end{equation}
is well defined and finite for every $\cJ$-Fredholm sequence $\bA \in \cA$. Since $(\alpha_n (\bA))$ is a sequence of non-negative integers, it possesses a constant subsequence the entries of which are equal to $\alpha(\bA)$ given by (\ref{e50.16}). Together with (\ref{e50.8}), this shows that
\begin{equation} \label{e50.17}
\liminf_{n \to \infty} \sigma_{\alpha(\bA)} (A_n) = 0 \quad \mbox{and} \quad \liminf_{n \to \infty} \sigma_{\alpha(\bA) + 1} (A_n) > 0.
\end{equation}
\begin{coro} \label{c50.18}
Let $(\cA, \, \cJ)$ be a Silbermann pair and $\bA \in \cA$ a $\cJ$-Fredholm sequence. Then the $\alpha$-number of the Fredholm sequence $\bA$ is given by $(\ref{e50.16})$.
\end{coro}
Let again $\bA = (A_n)$ be $\cJ$-Fredholm. Evidently, for large $n$, the singular values of $A_n$ are located in the union $[0, \, \varepsilon_n] \cup [d, \, \infty)$ where
\[
\varepsilon_n := \sigma_{\alpha_n(\bA)} (A_n) \quad \mbox{and} \quad d := \liminf_{n \to \infty} \sigma_{\alpha_n(\bA) +1} (A_n)/2.
\]
From (\ref{e50.8}) one concludes that $\varepsilon_n \to 0$ as $n \to \infty$ and $d > 0$. Thus, the singular values of the entries of $\cJ$-Fredholm sequence own the splitting property.

Note that the number of the singular values of $A_n$ which lie in $[0, \, \varepsilon_n]$ depends on $n$ in general (it is just given by the quantity $\alpha_n(\bA)$ in (\ref{e50.6})). A concrete instance where this dependence on $n$ can be observed occurs will be examined in Example \ref{ex50.23} below. The idea used there allows one to construct Silbermann pairs with arbitrarily prescribed weight sequences $(\alpha_n^t)$. On the other hand, many of the approximation methods used in practice have the property that every rank one projection in $\cJ/\cG$ lifts to a sequence of projections of rank one. Thus, in this case, the numbers $\alpha_n^t$ are independent on $n$ and can be chosen to be 1 for all $n$. For Silbermann pairs with this property, Theorem \ref{t50.7} and its Corollary \ref{c50.18} specify as follows.
\begin{coro} \label{c50.20}
Let $(\cA, \, \cJ)$ be a Silbermann pair where all weight sequences $(\alpha_n^t)$ are identically equal to one, and let $\bA \in \cA$ be a $\cJ$-Fredholm sequence. Then
\begin{equation} \label{e50.21}
\alpha(\bA) = \sum_{t \in T} \dim \ker W_t(\bA),
\end{equation}
and the sequence $\bA$ has the $\alpha(\bA)$-splitting property, i.e., the number of the singular values of $A_n$ which tend to zero is $\alpha(\bA)$.
\end{coro}
We are going to consider a few examples.
\begin{example} \label{ex50.22}
The simplest Silbermann pairs $(\cA, \, \cJ)$ arise when $\cJ/\cG$ is an elementary algebra. For a concrete model, let $\cP = (P_n)$ be a sequence of orthogonal projections of finite rank on a Hilbert space $H$ which converge strongly to the identity operator. Let $\cF^\cP$ denote the $C^*$-algebra of all sequences $\bA = (A_n)$ of operators $A_n : \im P_n \to \im P_n$ which converge $^*$-strongly to an operator $W(\bA)$. The set
\[
\cJ^\cP := \{ (P_n K P_n + G_n) : K \in K(H), \, (G_n) \in \cG \}
\]
forms a closed ideal of the algebra $\cF^\cP$, and $(\cF^\cP, \,
\cJ^\cP)$ is a Silbermann pair for which $\cJ^\cP/\cG^\cP$ is
$^*$-isomorphic to $K(H)$. Moreover, $\cJ^\cP$ is an algebra of weight one. In this setting, Theorems \ref{t50.5} and \ref{t50.7} and Corollary \ref{c50.20} specify as follows.
\begin{coro} \label{c50.22a}
Every $\cJ^\cP$-Fredholm sequence $\bA \in \cF^\cP$ owns the finite splitting property, and its splitting number $\alpha (\bA)$ is equal to $\dim \ker W(\bA)$. \hfill \qed
\end{coro}
\end{example}
\begin{example} \label{ex50.23}
Define sequences $\cP = (P_n)$ and $(R_n)$ as in Section \ref{ss13}. Consider the set $\cA$ of all sequences $\bA = (A_n)$ in $\cF^\cP$ for which the strong limits
\[
\mbox{s-lim} \, A_nP_n, \quad \mbox{s-lim} \, A_n^*P_n, \quad
\mbox{s-lim} \, R_n A_n R_n, \quad \mbox{s-lim} \, R_n A_n^* R_n
\]
exist. We denote the first and third of these strong limits by $W(\bA)$ and $\widetilde{W} (\bA)$, respectively. One can straightforwardly check that $\cA$ is a $C^*$-subalgebra of $\cF^\cP$, that $W$ and $\widetilde{W}$ are $^*$-homomorphisms on $\cA$, and that
\[
\cJ := \{ (P_nKP_n + R_nLR_n + G_n) : K, \, L \in K(l^2(\sZ^+)), \, (G_n) \in \cG^\cP \}
\]
is a closed ideal of $\cA$ for which $(\cA, \, \cJ)$ is a Silbermann pair. Moreover, the two involved weight sequences can be chosen to be identically one. The algebra $\cA$ contains all sequences $(P_n T(a) P_n)$ of the finite sections of Toeplitz operators $T(a)$ with generating function $a \in L^\infty(\sT)$. Thus, Corollary \ref{c50.20} implies the following result which holds for arbitrary bounded Toeplitz operators. Again, we set $\tilde{a}(t) = a(t^{-1})$.
\begin{coro} \label{c190309.2}
Let $a \in L^\infty(\sT)$. If the sequence $\bA := (P_n T(a) P_n)$ is invertible modulo the ideal $\cJ$, then $\bA$ is a Fredholm sequence with $\alpha$-number
\[
\alpha (\bA) = \dim \ker T(a) + \dim \ker T(\tilde{a}).
\]
\end{coro}
Note that neither a criterion for the Fredholm property of the finite sections sequence $(P_n T(a) P_n)$ with general $a \in L^\infty(\sT)$ nor an explicit formula for their $\alpha$-number is known. We also do not know anything on the fractality of such sequences. Recall in this connection that Treil constructed an invertible Toeplitz
operator for which the finite sections sequence $(P_n T(a) P_n)$ fails to be stable. It is not known if Treil's sequence is Fredholm. \hfill \qed
\end{example}
\begin{example} \label{ex50.24}
Here we present an example with non-constant weight. Let the operators $P_n$ and $R_n$ be as in Example \ref{ex50.23},
and set $\cP = (P_n)_{n \ge 1}$. Consider the smallest closed subalgebra $\cA$ of $\cF^\cP$ which contains the identity sequence $(I_n)$, the ideal $\cG$ of the zero sequences and all sequences $(K_n)$ of the form
\[
K_n := \left\{ \begin{array}{ll}
P_n K P_n & \mbox{if} \; n \; \mbox{is odd} \\
P_n K P_n + R_n K R_n & \mbox{if} \; n \; \mbox{is even}
\end{array} \right.
\]
where $K \in K(l^2(\sZ^+))$. One easily checks that if $(A_n)$ is a sequence in $\cA$ then every entry $A_n$ is of the form
\begin{equation} \label{e50.25}
A_n := \left\{ \begin{array}{ll}
\gamma I_n + P_n K P_n + G_n & \mbox{if} \; n \; \mbox{is odd} \\
\gamma I_n + P_n K P_n + R_n K R_n + G_n & \mbox{if} \; n \;
\mbox{is even}
\end{array} \right.
\end{equation}
where $\gamma \in \sC$, $K$ is compact, and $(G_n) \in \cG$.
Clearly, $\cA$ is a unital $C^*$-subalgebra of $\cF^\cP$, and the mapping
\[
W : \cA \to L(l^2(\sZ^+)), \quad (A_n) \mapsto \mbox{s-lim} \, A_n P_n
\]
is a representation of $\cA$ which maps the sequence $(A_n)$ given by (\ref{e50.25}) to the operator $\gamma I + K$. We show that the invertibility of the operator $\gamma I + K$ implies the stability of the sequence $\bA$ defined by $(\ref{e50.25})$. Indeed, consider the strictly increasing sequences $\mu, \, \eta: \sN \to \sN$ given by $\mu(n) = 2n$ and $\eta(n) = 2n+1$. By Theorem \ref{t14.21}, the restricted sequences $\bA_\mu$ and $\bA_\eta$ belong to the corresponding restricted algebras $\cS(\eT(C))_\mu$ and $\cS(\eT(C))_\eta$ of the finite sections method for Toeplitz operators, respectively. Since
\[
R_{2n} A_{2n} R_{2n} \to \gamma I + K = W(\bA) \quad \mbox{strongly as} \; n \to \infty,
\]
we conclude from Corollary \ref{c160309.2} that the invertibility of $W(\bA) = \gamma I + K$ implies the stability of both $\bA_\mu$ and $\bA_\eta$ and, hence of the sequence $\bA$.

It is further easy to check that the set $\cJ$ of all sequences of the form (\ref{e50.25}) with $\gamma = 0$ forms a closed ideal of
$\cA$ and that the quotient algebra $\cJ/\cG$ is $^*$-isometric
(via $W$) to $K(l^2(\sZ^+))$. Set $\Pi_n := P_n P_1 P_n$ if $n$ is odd and $\Pi_n := P_n P_1 P_n + R_n P_1 R_n$ if $n$ is even. Then
the sequence $(\Pi_n)$ belongs to $\cJ$, the coset $p := (\Pi_n) + \cG$ is a non-trivial minimal projection in $\cJ/\cG$, and one has
\[
\dim \im \Pi_n = \left\{ \begin{array}{ll}
1 & \mbox{if} \; n \; \mbox{is odd} \\
2 & \mbox{if} \; n \; \mbox{is even}.
\end{array} \right.
\]
Thus, the alternating sequence $(1, \, 2, \, 1, \, 2, \, \ldots)$ is a representative of the (only) weight related with $\cJ$, and the identity (\ref{e50.6}) specifies to
\[
\alpha_n (\bA) = \left\{ \begin{array}{ll}
\dim \ker W_t(\bA) & \mbox{if} \; n \; \mbox{is odd} \\
2 \, \dim \ker W_t(\bA) & \mbox{if} \; n \; \mbox{is even}
\end{array} \right.
\]
(which also could have been verified directly without effort). \hfill \qed
\end{example}
\begin{example}
The smallest closed subalgebra of $\cF$ which contains all sequences $(P_n T(a) P_n)$ where $a \in C(\sT)$ and $a = \tilde{a}$ is of constant weight two. It is left to the reader to work out the details. \hfill \qed
\end{example}
\subsection{Complete Silbermann pairs} \label{ss45}
Let $(\cA, \cJ)$ be a Silbermann pair. We call this pair {\em complete} \index{Silbermann pair!complete} if the ideal $\cG$ is properly contained in $\cJ$ and if the family $\{W_t\}_{t \in T}$ of the lifting homomorphisms of $(\cA, \cJ)$ is sufficient for stability in the sense that a sequence $\bA \in \cA$ is stable if and only if the operators $W_t(\bA)$ are invertible for every $t \in T$. We call the pair $(\cA, \cJ)$ {\em weakly complete} if a sequence $\bA \in \cA$ is stable if and only if the operators $W_t(\bA)$ are invertible for every $t \in T$ and if the norms of their inverses are uniformly bounded. In this case we call the family of the $W_t$ {\em weakly sufficient} for $\cA$. Finally, we call a unital fractal $C^*$-subalgebra $\cA$ of $\cF$ a {\em Silbermann algebra} \index{algebra!Silbermann} if $(\cA, \, \cA \cap \cK)$ is a weakly complete Silbermann pair. The latter "weak" notions will be needed in Section \ref{s55}.
\begin{theo} \label{t50.25a}
Let $(\cA, \, \cJ)$ be a complete Silbermann pair and let $\bA \in \cA$. Then \\[1mm]
$(a)$ $\bA$ is stable if and only if all operators $W_t(\bA)$ are invertible; \\[1mm]
$(b)$ $\|\bA + \cG\|_{\cF/\cG} = \max_{t \in T} \|W_t (\bA)\|$. \\[1mm]
$(c)$ $\bA$ is $\cJ$-Fredholm if and only if all operators $W_t(\bA)$ are Fredholm and if there are only finitely many of them which are not invertible; \\[1mm]
$(d)$ $\bA \in \cJ$ if and only if all operators $W_t(\bA)$ are compact and if, for each $\varepsilon > 0$, there are only finitely many of them with $\|W_t(\bA)\| > \varepsilon$.
\end{theo}
{\bf Proof.} Assertion $(a)$ is a re-formulation of the sufficiency condition. Assertion $(b)$ is a consequence of $(a)$, since every $^*$-homomorphism between $C^*$-algebras which preserves spectra is an isometry. \\[1mm]
$(c)$ The 'only if' part of assertion $(c)$ follows from the Lifting theorem \ref{t50.5} $(d)$. Conversely,
let $\bA \in \cA$ be a sequence for which all operators $W_t(\bA)$ are Fredholm and for which there is a finite subset $T_0$ of $T$ which consists of all $t$ such that $W_t(\bA)$ is not invertible. Then all operators $W_t(\bA^* \bA)$ are Fredholm, and they are invertible if $t \notin T_0$. Let $t \in T_0$. Then $W_t(\bA^*\bA)$ is a Fredholm operator of index $0$. Hence, there is a compact operator $K_t$ such that $W_t(\bA^* \bA) + K_t$ is invertible. Choose a sequence $\bK_t \in \cJ$ with $W_t (\bK_t) = K_t$ and $W_s(\bK_t) = 0$ for $s \neq t$ (which is possible by the separation property in Theorem \ref{t50.5}), and set
\[
\bK := \sum_{t \in T_0} \bK_t.
\]
Then $\bK$ belongs to the ideal $\cJ$, and all operators $W_t(\bA^* \bA + \bK)$ are invertible. By assertion $(a)$, the sequence $\bA^* \bA + \bK$ is stable. Similarly, one finds a sequence $\bL \in \cJ$ such that $\bA \bA^* + \bL$ is a stable sequence. Consequently, the sequence $\bA$ is invertible modulo $\cJ$, whence the $\cJ$-Fredholm property of that sequence. \\[1mm]
$(d)$ Since $\cJ/\cG$ is a dual algebra, the 'only if' part follows again from the Lifting theorem \ref{t50.5} $(c)$. For the 'if' part, let $\bK \in \cA$ be a sequence such that, for every $\varepsilon > 0$, there are only finitely many $t \in T$ with $\|W_t(\bA)\| > \varepsilon$. For $n \in \sN$, let $T_n$ stand for the (finite) subset of $T$ which collects all $t$ with $\|W_t (\bK)\| > 1/n$. For each $t \in T_n$, choose a sequence $\bK^t \in \cJ$ with $W_t (\bK^t) = W_t(\bK)$ and $W_s(\bK^t) = 0$ for $s \neq t$ (which can be done by the separation property in Theorem \ref{t50.5} again), and set $\bK_n := \sum_{t \in T_n} \bK^t$. Then $W_t(\bK - \bK_n) = 0$ for $t \in T_n$ and $W_t(\bK_n) = 0$ for $t \notin T_n$. Hence, $\sup_{t \in T} \|W_t(\bK - \bK_n)\| \le 1/n$ for every $n \in \sN$. By Theorem \ref{t50.25a} $(b)$, the left-hand side coincides with $\|\bK - \bK_n + \cG\|_{\cF/\cG}$. Being the norm limit of a sequence in $\cJ$, the sequence $\bK$ belongs to $\cJ$ itself. \hfill \qed \\[3mm]
An example for a complete Silbermann pair is $(\cS(\eT(C)), \, \cJ)$ consisting of the algebra of the finite sections method for Toeplitz operators and its distinguished ideal (\ref{e14.35a}). A sequence in the algebra $\cS(\eT(C))$ is Fredholm if and only its strong limit is a Fredholm operator (note that $T(a)$ and $T(\tilde{a})$ are Fredholm only simultaneously). Equivalently, the sequence $\bA := (P_n T(a) P_n + P_n K P_n + R_n L R_n + G_n)$ with $a \in C(\sT)$, $K, \, L$ compact and $(G_n) \in \cG$ is Fredholm if and only if $T(a)$ is a Fredholm operator. In this case,
\begin{equation} \label{e42.5a}
\alpha (\bA) = \dim \ker (T(a) + K) + \dim \ker (T(\tilde{a}) + L).
\end{equation}
In particular, if $K = L = 0$, then
\begin{eqnarray*}
\alpha (\bA) & = & \dim \ker T(a) + \dim \ker T(\tilde{a}) \\
& = & \max \{ \dim \ker T(a), \, \dim \ker T(\tilde{a}) \}
\end{eqnarray*}
where the second equality holds by a theorem of Coburn which states that one of the quantities $\dim \ker T(a)$ and $\dim \ker T(\tilde{a})$ for each non-zero Toeplitz operator.
\subsection{The extension-restriction theorem} \label{s55}
Let $\cF^\delta$ be the algebra of matrix sequences with dimension function $\delta$ and $\cG^\delta$ the associated  ideal of zero sequences. We say that a $C^*$-subalgebra $\cA^{ext}$ of $\cF^\delta$ is an extension of a $C^*$-subalgebra $\cA$ of $\cF^\delta$ by compact sequences if there is a subset $\cK^\prime$ of the ideal $\cK^\delta$ of the compact sequences in $\cF^\delta$ such that $\cA^{ext}$ is the smallest $C^*$-subalgebra of $\cF^\delta$ which contains $\cA$ and $\cK^\prime$.
\begin{theo} \label{t55.1}
Let $\cA$ be a unital separable $C^*$-subalgebra of an algebra $\cF^\delta$. Then there are an extension $\cA^{ext}$ of $\cA$ by compact sequences and a strictly increasing sequence $\eta$ such that the restriction $\cA^{ext}_\eta$ is a Silbermann algebra.
\end{theo}
In other words, after extending $\cA$ by adding a suitable set of compact sequences and then passing to a suitable restriction, we arrive at a weakly complete Silbermann pair $(\cA^{ext}_\eta, \, \cA^{ext}_\eta \cap \cK_\eta)$, i.e., the algebra $\cA^{ext}_\eta$ is fractal, and the family of the lifting homomorphisms of its dual ideal $\cA^{ext}_\eta \cap \cK_\eta$ is weakly sufficient for the stability of sequences in $\cA^{ext}_\eta$. \\[3mm]
{\bf Proof.} Let $\cA_0$ be a countable dense subset of $\cA$. We will also assume that the identity sequence is in $\cA_0$. The set $\cA_0^* \cA_0$ is still countable and dense in $\cA$. For each sequence $\bA = (A_n)$ in $\cA_0$, we write
\begin{equation} \label{e55.2}
A_n^*A_n = E_n^* \diag (\lambda_1 (A_n), \, \ldots, \, \lambda_{\delta(n)} (A_n)) E_n
\end{equation}
with a unitary matrix $E_n$ and increasingly ordered eigenvalues  $0 \le
\lambda_1 (A_n) \le \ldots \le \lambda_{\delta(n)} (A_n)$. For $l, \, r \in
\sN$, let $K_{l,r,n}$ be the $\delta(n) \times \delta(n)$-matrix which is zero
if $\max \{l, \, r\} > \delta(n)$ and which has a 1 at the $lr$th entry and
zeros at all other entries if $\max \{l, \, r\} \le \delta(n)$. The sequence
$\bK^{\bA,l,r}$ with entries
\[
K^{\bA,l,r}_n := E_n^* K_{l,r,n} E_n
\]
is a sequence of rank one matrices, hence compact. Let $\cA^{ext}$ stand for the smallest $C^*$-subalgebra of $\cF^\delta$ which contains the algebra $\cA$, the ideal $\cG^\delta$, and all sequences $\bK^{\bA,l,r}$ with $\bA \in \cA_0$ and $l, \, r \in \sN$. It is a simple exercise to show that the algebra $\cA^{ext}$ is still separable. Hence, by Theorems \ref{tf3.20} and \ref{t43.25}, there is a strictly increasing sequence $\eta$ such that the restriction $\cA^{ext}_\eta$ is fractal and essentially fractal. We claim that $\cA^{ext}_\eta$ is a Silbermann algebra. Note that the sequences $\bK^{\bA,r} := \sum_{l=1}^r \bK^{\bA,l,l}$ with entries
\begin{equation} \label{e55.3}
K^{\bA,r}_n := E_n^* \diag (1, \, \ldots, \, 1, \, 0, \, \ldots, \, 0) E_n
\end{equation}
where $r$ ones followed by $\delta(n) - r$ zeros in the diagonal part belong to $\cA^{ext}$.

To simplify notation, we will assume that $\eta$ is the identity mapping (otherwise replace $\delta$ by $\delta \circ \eta$ in what follows), and we denote the (restricted) ideal $\cA^{ext} \cap \cK^\delta$ by $\cJ$. Since $\cA$ is a $C^*$-subalgebra and $\cJ$ is a closed ideal of $\cA^{ext}$, the algebraic sum $\cA + \cJ$ is a $C^*$-subalgebra of $\cA^{ext}$. This subalgebra contains $\cA$, $\cG^\delta$ and all sequences $\bK^{\bA,l,r}$. Thus, $\cA^{ext} = \cA + \cJ$.

Since $\cA^{ext}$ is fractal, the ideal $\cJ/\cG^\delta$ is dual by Corollary \ref{c190309.1}. Let $(I_t)_{t \in T}$ denote the set of its elementary components and, for each $t \in T$, let $W_t : I_t \to L(H_t)$ stand for the associated irreducible representation. As earlier, we will denote an irreducible representation of $I_t$ and its irreducible extensions to $\cA^{ext}/\cG^\delta$ and $\cA^{ext}$ by the same symbol.

Let $\bA \in \cA_0$. We claim that the coset $\bK^{\bA,1,1} + \cG^\delta = \bK^{\bA,1} + \cG^\delta$ is a rank one projection in $(\cA^{ext} \cap \cK^\delta)/\cG^\delta$. Indeed, the entries $K^{\bA,1}_n$ are projection matrices of rank one. Hence, for every positive sequence $(B_n^*B_n) \in \cA^{ext}$, there is a sequence $(\beta_n)$ of complex numbers such that
\[
K^{\bA,1}_n B_n^*B_n K^{\bA,1}_n = \beta_n K^{\bA,1}_n \quad \mbox{for every} \; n \in \sN.
\]
The sequence $(\beta_n K^{\bA,1}_n)_{n \in \sN}$ is fractal, and $\beta_n$ is the largest singular value of $\beta_n K^{\bA,1}_n$. By Proposition \ref{p51.2}, the sequence $(\beta_n)$ is convergent. Since every sequence in $\cA^{ext}$ is a linear combination of four positive sequences, we conclude that, for every sequence $\bC =(C_n) \in \cA^{ext}$, there is a convergent sequence $(\gamma_n)$ of complex numbers such that
\[
K^{\bA,1}_n C_n K^{\bA,1}_n = \gamma_n K^{\bA,1}_n \quad \mbox{for every} \; n \in \sN.
\]
Put $\gamma := \lim_{n \to \infty} \gamma_n$. Then $\bK^{\bA,1} \bC \bK^{\bA,1} - \gamma \bK^{\bA,1} \in \cG^\delta$, which proves the claim.

Since the elementary components of $\cJ/\cG^\delta$ are generated by rank one projections, there is a $t(\bA) \in T$ such that $\bK^{\bA,1} + \cG^\delta \in I_{t(\bA)}$. Since $I_{t(\bA)}$ is an ideal, the equality
\[
\bK^{\bA,l,r} = \bK^{\bA,l,1}\bK^{\bA,1,1}\bK^{\bA,1,r}
\]
implies that $\bK^{\bA,l,r} + \cG^\delta \in I_{t(\bA)}$ for every pair $l, \, r \in \sN$. In particular, all cosets
$\bK^{\bA,r} + \cG^\delta$ belong to $I_{t(\bA)}$. Since the cosets $\bK^{\bA,l,l} + \cG^\delta$ are linearly independent rank one projections and the homomorphism $W_{t(\bA)}$ is irreducible, the operators $W_{t(\bA)} (\bK^{\bA,l,l})$ form a linearly independent set of projection operators of rank one in $L(H_{t(\bA)})$. In particular, the Hilbert space $H_{t(\bA)}$ has infinite dimension.

Since $\cA_0$ is dense in $\cA$, the sequences $\bA + \bK$ with $\bA \in \cA_0$ and $\bK \in \cA^{ext} \cap \cK^\delta$ form a dense subset $\cA_0^{ext}$ of $\cA^{ext}$. Let $\bB := \bA + \bK$ be a sequence of this form, for which $\bB^* \bB$ is not a Fredholm sequence (equivalently, $\bB^* \bB$ is not a $\cJ$-Fredholm sequence, since $\cJ$ contains all compact sequences in $\cA^{ext}$). Then $\bA = (A_n)$ is not a Fredholm sequence, hence
\begin{equation} \label{e55.4}
\lim_{n \to \infty} \lambda_r (A_n) = 0 \quad \mbox{for every} \; r \in \sN
\end{equation}
by (\ref{e43.19}) (see (\ref{e55.2}) and recall that the algebra $\cA^{ext}$ is essentially fractal after restriction). From (\ref{e55.2}) -- (\ref{e55.4}) we conclude that $\bA^*\bA \bK^{\bA,r} \in \cG^\delta$ for every $r \in \sN$, hence
\begin{equation} \label{e55.5}
W_{t(\bA)} (\bA^* \bA) W_{t(\bA)} (\bK^{\bA,r}) = 0 \quad \mbox{for every} \; r \in \sN.
\end{equation}
Since $(W_{t(\bA)} (\bK^{\bA,r}))_{r \ge 1}$ is an increasing sequence of orthogonal projections on $H_{t(\bA)}$, this sequence converges strongly, and its limit, $P$, is the orthogonal projection from $H_{t(\bA)}$ onto the closure of the linear span of the union of the ranges of the $W_{t(\bA)} (\bK^{\bA,r})$ (see, for example, Theorem 4.1.2 in \cite{Mur1}). So we conclude from (\ref{e55.5}) that $W_{t(\bA)} (\bA^* \bA) P = 0$. Thus, and by Theorem \ref{t50.5} $(c)$,
\[
W_{t(\bA)} (\bB^* \bB) P = W_{t(\bA)} (\bA^* \bA) P + W_{t(\bA)} (\bA^*\bK + \bK^* \bA + \bK^*\bK) P
\]
is a compact operator. Then $W_{t(\bA)} (\bB^* \bB)$ cannot be invertible: otherwise, the projection $P$ were compact, but the range of $P$ has infinite dimension, which follows by the same arguments as the infinite dimensionality of $H_{t(\bA)}$.

Thus, whenever $\bB \in \cA_0^{ext}$ and $\bB^* \bB$ is not a Fredholm sequence, then one of the operators $W_t(\bB^* \bB)$ is not invertible. Conversely, if all operators $W_t(\bB^* \bB)$ with $t \in T$ are invertible, then $\bB^* \bB$ is a Fredholm sequence. By Theorem \ref{t50.5} $(e)$ this implies that, whenever all operators $W_t(\bB^* \bB)$ with $t \in T$ are invertible, then the sequence $\bB^* \bB$ is a stable. Since this fact holds for all sequences $\bB$ in the dense subset $\cA_0^{ext}$ of $\cA^{ext}$, the family $(W_t)_{t \in T}$ is weakly sufficient for $\cA^{ext}$, as one easily checks. \hfill \qed \\[3mm]
It is not clear if one can force also a {\em complete} Silbermann pair by a similar construction. The point is that the implication, obtained at the end of the previous proof, does only hold for sequences in a dense subset. Another open question is if there is a version without restriction if one starts with a fractal and essentially fractal separable subalgebra $\cA$.
%

%
%
%\cleardoublepage
%\input{asnuba.ind}
%
%
\end{document}